\documentclass[12pt]{article}
\usepackage{amssymb}
\usepackage{amsthm}
\usepackage{amsgen}
\usepackage{amsfonts}
\usepackage{amsbsy}
\usepackage{float} 
\usepackage[dvips]{graphicx}
\usepackage[dvips]{color}
\usepackage[all]{xy}
\usepackage{makeidx} 

\newtheorem{theorem}{Theorem}[section] 
\newtheorem{corollary}[theorem]{Corollary} 
\newtheorem{conjecture}[theorem]{Conjecture} 
\newtheorem{lemma}[theorem]{Lemma} 
\newtheorem{proposition}[theorem]{Proposition} 
\newtheorem{definition}[theorem]{Definition} 
\newtheorem{example}[theorem]{Example} 
 
\newtheorem{remark}[theorem]{Remark} 
\newtheorem{magma}[theorem]{MAGMA/GAP Example} 
\newcommand{\mthree}[9]{\left(
        \begin{array}{ccc}#1&#2&#3\\#4&#5&#6\\#7&#8&#9
        \end{array}\right)}
\newcommand{\mtwo}[4]{\left(
        \begin{array}{cc}#1&#2\\#3&#4
        \end{array}\right)}
\def\pf{{\bf proof}:\ }
\def\mod{\ {\rm mod}\ }
\def\qed{$\square$}
\def\Spec{{\rm Spec}}
\def\fff{\mathbb{F}}
\def\aaa{\mathbb{A}}
\def\A{\mathbb{A}}
\def\ppp{\mathbb{P}}
\def\P{\mathbb{P}}
\def\qqq{\mathbb{Q}}
\def\Q{\mathbb{Q}}
\def\rrr{\mathbb{R}}
\def\R{\mathbb{R}}
\def\ccc{\mathbb{C}}
\def\cccF{F}
\def\zzz{\mathbb{Z}}
\def\Z{\mathbb{Z}}
\def\ggg{\mathbb{G}}
\def\hhh{{\cal H}}
\def\eee{{\cal E}}
\def\ra{{\rightarrow}}
\newcommand{\hra}{\hookrightarrow}

\makeindex

\begin{document}

\author{
Helena Verrill\thanks{Math Dept, Univ. Hannover, Germany, 
verrill\@@math.uni-hannover.de}
\, and 
\, David Joyner\thanks{Math Dept, USNA, wdj\@@usna.edu}
}
\title{Notes on toric varieties}
\date{7-25-2002}
\maketitle

\begin{abstract}
These notes survey some basic results in toric varieties
over a field $\cccF$, with examples and applications.

A computer algebra 
package (written by the author) is described which
deals with both affine and projective 
toric varieties in any number of dimensions
(written in both MAGMA \cite{M} and GAP \cite{GAP}).
Among other things, the package implements the desingularization
procedure, constructs some error-correcting codes
associated with toric varieties, and computes the
Riemann-Roch space of a divisor on a toric variety.

\end{abstract}

\tableofcontents

\vskip .3in

The main reference for these notes is \cite{F}.

Throughout, let $\cccF$ denote a field.

\section{Introduction}

Let $R=\cccF [x_1,...,x_n]$ be a ring in $n$ variables.
A {\bf binomial relation}
\footnote{Also called a {\bf monomial equation} \cite{ES}.}
in $R$ is one of the form
\[
x_1^{k_1}...x_n^{k_n}=x_1^{\ell_1}...x_n^{\ell_n},
\]
where $k_i\geq 0,\ \ell_j\geq 0$ are integers.
A {\bf binomial variety}
\footnote{Also called a monomial variety, binomial scheme, torus embedding,
...\ . Such a variety need not be normal.} 
is a subvariety of complex affine $n$-space
$\aaa_{\cccF}^n$ defined by a finite set of binomial 
equations. Such varieties arise frequently ``in nature''.
A typical ``toric variety'' (defined more precisely later)
is binomial, though they will be introduced via an {\it a priori}
independent construction
\footnote{In fact, we shall give two constructions
of toric varieties. The basic idea of the first one is to
replace each such binomial equation as above by a relation in a 
semigroup contained in a lattice
and replace $R$ by the ``group algebra'' of this 
semigroup. More details below. By the way, a toric variety is 
always normal (see for example, \cite{F}, page 29).
}.

\subsection{Motivation}

Toric geometry has several interesting facets.
\begin{itemize}
\item
It gives a way of constructing an ambient space (a toric variety)
in which algebraic varieties live, in the same way that 
projective space, and weighted projective space are primarily considered
as ambient spaces.  
\begin{itemize}
\item
Toric varieties are all rational.
\item
Toric varieties are not all projective, (i.e., they cannot all be
embedded in $\P^n$ for some n), so provide 
a wider class of algebraic varieties to work with.
\item
Even if a subvariety of a toric
variety is projective, it generally is embedded
in a toric variety with affine pieces of much lower dimension than the
projective space.
\end{itemize}
\item
In the case of a toric variety, certain algebraic geometry 
computations can be reformulated into simpler
combinatorial problems. For example,
there is a ``simple'' method of resolving certain singularities.
\item 
Toric geometry is used to solve certain compactification problems.
\item 
Toric geometry is a source of many examples in various areas, e.g., 
Batyrev's construction of 
pairs of mirror Calabi-Yau threefolds from ``reflexive polytopes''.
\item 
It has many other applications, e.g.,
to combinatorical geometry \cite{E}, error-correcting
codes \cite{Han}, and connections with 
germs of p-adic orbital integrals (unpublished) and
non-regular continued fractions \cite{F}. 
\end{itemize}

\subsection{Computer programs}

Toric geometry very much lends itself to a computational approach.
There are a few packages available for working with toric geometry, such as
\begin{itemize}
\item
\begin{itemize}
\item TiGERS
\item c program
\item by Birkett Huber and Rekha Thomas 
\item {\sf http://www.math.washington.edu/$\widetilde{\>}$thomas/programs.html}
\item Computes the Gr\"obner fan of a toric ideal
\item Also some maple available from same site:\newline
\centerline{\sf http://www.math.washington.edu/$\widetilde{\>}$thomas/program}
\end{itemize}
\item
\begin{itemize}
\item part of the \emph{singular} computer algebra system
\item
{\sf 
http://www.singular.uni-kl.de/Manual/2-0-0/sing\underline{\hspace{1ex}}376.htm}
\item{\sf http://www.singular.uni-kl.de/Manual/2-0-0/sing\underline{\hspace{1ex}}378.htm}
\item Computes the lattice basis using LLL algorithm, and algorithms by
Conti and Traverso;
Pottier;     
Hosten and Sturmfels;
Di Biase and Urbanke;
 Bigatti, La Scala and Robbiano,
for computing the saturation.
\end{itemize}
\item 
\begin{itemize}
\item ``CY/4d''
\item Kreuzer and Skarke
\item Data base of the classification of all  473,800,776 reflexive polytopes 
in $4$ dimensions.  ($3$ dimensional version also available.)
\item {\sf http://hep.itp.tuwien.ac.at/$\widetilde{\>}$kreuzer/CY/CYcy.html}
\item this gives a method of constructing different examples of mirror pairs of
families of Calabi-Yau threefolds, especially of interest to physicists.
\end{itemize} 
\item
\begin{itemize}
\item ``toric'' (\verb+toric.g+ for GAP 4.3 and \verb+toric.mag+ for MAGMA 2.8)
\item Joyner
\item GAP and MAGMA programs for toric varieties.
\item \cite{J2}
Examples are given in this paper and at the end of each of the
files in \cite{J2} to explain the syntax.
\end{itemize} 
\end{itemize} 
Note that all but the last two of the above programs are restricted to the
case of affine toric varieties, 
(the first two being more from the point of commutative algebra),
and the third one is a fairly specialized
situation, and only available as a data base.  
There is much more to toric geometry that could be programmed 
but so far seems not to have been.

\section{Cones and semigroups}

Let $V=\qqq^n$ having basis
$f_1=(1,0,...,0)$, ..., $f_n=(0,...,0,1)$.
Let $L$ be a lattice in $V$ (i.e., a free rank 
$n$ $\zzz$-module in $V$).
We identify $V$ and $L\otimes_{\zzz}\qqq$.
We use $\cdot$ or $\langle\ ,\ \rangle$ to denote
the (standard) inner product on $V$.
If $e_1,...,e_n$ is a basis for $L$, let
$e^*_1$, ..., $e_n^*$ denote the dual basis
for $L^*$. Let
\[
L^*={\rm Hom}(L,\zzz)=\{ v\in V\ |\ \langle v,w \rangle \in \zzz,
\ \forall w\in L\}
\]
denote the {\bf dual lattice}, so (since the basis 
of $L$ has been choosen to be the standard one)
$L^*$ may be identified with $\zzz^n$.
\index{dual lattice}

\begin{example}
If $L$ is the sublattice of $\zzz^2$ generated by $e_1+2e_2$ and
$2e_1+e_2$ then the dual is not contained in $\zzz^2$.
However, it may be rescaled (multiplied by an integer)
so that it is. MAGMA's \verb+Dual+ command compute's 
this rescaled version.

{\footnotesize{
\begin{verbatim}
> L := LatticeWithBasis(2, [1,2, 2,1]);
> L;
Lattice of rank 2 and degree 2
Basis:
(1 2)
(2 1)
> Lperp:=Dual(L);                      
> Lperp;
Lattice of rank 2 and degree 2
Basis:
( 1  1)
( 2 -1)
> a,b,c:=DualQuotient(L);
> b;
Lattice of rank 2 and degree 2
Basis:
( 1  1)
( 1 -2)
Basis Denominator: 3
> Basis(b);
[
    (1/3 1/3),
    ( 1/3 -2/3)
]
\end{verbatim}
}}
\noindent
This means that 
{\tt Lperp}$=\{a(e_1^*+e_2^*)+ b(2e_1^*-e_2^*)\ |\ a,b\in\zzz\}$
is an integral basis belonging to the dual $L^*$, denoted
{\tt b} in the above MAGMA session. In fact,
\[
\begin{array}{c}
L^*=\{a\frac{e_1^*+e_2^*}{3}+b\frac{e_1^*-2e_2^*}{3}
\ |\ m,n\in\zzz\}\\
=\{(a+\eta)e_1^*+(b+\eta)e_2^*\ |\ a,b\in\zzz,\ 
\eta\in\{0,\frac{1}{3},\frac{2}{3}\}\},
\end{array}
\]
so {\tt Lperp} is ``$3$ times $L^*$''.

The lattice $L$ may be visualized as follows.


\begin{center}
\setlength{\unitlength}{.01cm}

\begin{picture}(256.00,200.00)(-130.00,0.00)
\thicklines
\put(0.00,0.00){\circle*{4}} 
\put(60.00,120.00){\circle*{4}} 
\put(120.00,60.00){\circle*{4}} 
\put(-60.00,60.00){\circle*{4}} 
\put(60.00,-60.00){\circle*{4}} 
\put(120.00,-120.00){\circle*{4}} 
\put(-120.00,-60.00){\circle*{4}} 
\put(-60.00,-120.00){\circle*{4}} 
\put(180.00,0.00){\circle*{4}} 
\put(0.00,180.00){\circle*{4}} 
\put(0.00,60.00){\circle{2}} 
\put(0.00,120.00){\circle{2}} 
\put(-60.00,0.00){\circle{2}} 
\put(-120.00,0.00){\circle{2}} 
\put(120.00,-60.00){\circle{2}} 
\put(0.00,-60.00){\circle{2}} 
\put(0.00,-120.00){\circle{2}} 
\put(-60.00,-60.00){\circle{2}} 
\put(60.00,0.00){\circle{2}} 
\put(120.00,0.00){\circle{2}} 
\put(60.00,120.00){\circle{2}} 
\put(60.00,-120.00){\circle{2}} 
\put(60.00,60.00){\circle{2}} 
\put(0.00,0.00){\vector(2,1){118.00}} 
\put(0.00,0.00){\vector(1,2){58.00}} 
\end{picture}
\end{center}

\end{example}

\vskip 2in

A plain vanilla
\footnote{Sorry, couldn't resist the pun.
The point is that though this is the correct and standard definition of a
cone, later we shall actually reserve the word ``cone'' for a special 
type of cone.}
{\bf cone} in $V$ is a $\sigma$ of the form
\[
\sigma=\{a_1v_1+...+a_mv_m\ |\ a_i\geq 0\},
\]
where $v_1,...,v_m\in V$ is a given collection of 
vectors (so $0\leq m\leq n$), 
called a {\bf basis} (or {\bf generators}) of $\sigma$. 
This cone is also denoted
\[
\sigma = \qqq_{\geq 0}[v_1,...,v_m].
\]
A {\bf rational cone} is one where $v_1,...,v_m\in L$.
A {\bf strongly convex} cone is one which contains
no lines through the origin. A cone generated by linearly independent 
vectors is necessarily strongly convex.
The {\bf dimension of $\sigma$}, denoted dim$(\sigma)$, is
the dimension of the subspace $\sigma+(-\sigma)$
of $V$. 
\index{cone, rational}
\index{cone, strongly convex}

\begin{example}
Let $R$ be a root system of $V$, endowed with the usual 
Euclidean inner product (as in Ch III of 
Humphreys \cite{H}, for example), and
$R_+$ a choice of positive roots. Let $L$ be the root lattice
\footnote{One may replace ``root lattice''
by ``character lattice'' or ``weight lattice''
as well.} of $R$ and $L^*$ its dual lattice (this is {\it not}
the lattice of ``dual roots''). 
Let $\sigma$ be the cone defined by the positive Weyl 
chamber associated to $R_+$. This is a strongly convex
rational cone.

For further examples, see \S\S 2-3 in Verrill \cite{V}.
\end{example}

\begin{example}
\label{ex:3}
Picture of a three dimensional cone and its faces:

\begin{picture}(0,0)%
\includegraphics{cone1.pstex}%
\end{picture}%
\setlength{\unitlength}{1579sp}%
\begingroup\makeatletter\ifx\SetFigFont\undefined%
\gdef\SetFigFont#1#2#3#4#5{%
  \reset@font\fontsize{#1}{#2pt}%
  \fontfamily{#3}\fontseries{#4}\fontshape{#5}%
  \selectfont}%
\fi\endgroup%
\begin{picture}(3600,3510)(-44,-4483)
\put(451,-4411){\makebox(0,0)[lb]{\smash{\SetFigFont{6}{7.2}{\familydefault}{\mddefault}{\updefault}{\color[rgb]{0,0,0}$y$}%
}}}
\put(1051,-1336){\makebox(0,0)[lb]{\smash{\SetFigFont{6}{7.2}{\familydefault}{\mddefault}{\updefault}{\color[rgb]{0,0,0}$\tau_1$}%
}}}
\put(2746,-1141){\makebox(0,0)[lb]{\smash{\SetFigFont{6}{7.2}{\familydefault}{\mddefault}{\updefault}{\color[rgb]{0,0,0}$\tau_2$}%
}}}
\put(1711,-1696){\makebox(0,0)[lb]{\smash{\SetFigFont{6}{7.2}{\familydefault}{\mddefault}{\updefault}{\color[rgb]{0,0,0}$\tau_3$}%
}}}
\put(-44,-1951){\makebox(0,0)[lb]{\smash{\SetFigFont{6}{7.2}{\familydefault}{\mddefault}{\updefault}{\color[rgb]{0,0,0}$\tau_4$}%
}}}
\put(3556,-1561){\makebox(0,0)[lb]{\smash{\SetFigFont{6}{7.2}{\familydefault}{\mddefault}{\updefault}{\color[rgb]{0,0,0}Cone $\sigma$ generated by }%
}}}
\put(3556,-1846){\makebox(0,0)[lb]{\smash{\SetFigFont{6}{7.2}{\familydefault}{\mddefault}{\updefault}{\color[rgb]{0,0,0}$(1,0,0), (1,1,0),(1,1,1)$ and $(1,0,1)$ (circled).}%
}}}
\put(3556,-2701){\makebox(0,0)[lb]{\smash{\SetFigFont{6}{7.2}{\familydefault}{\mddefault}{\updefault}{\color[rgb]{0,0,0}$\tau_2=\R_{\ge 0}(1,1,0)$, $\tau_3=\R_{\ge 0}(1,1,1)$ }%
}}}
\put(3556,-4411){\makebox(0,0)[lb]{\smash{\SetFigFont{6}{7.2}{\familydefault}{\mddefault}{\updefault}{\color[rgb]{0,0,0}The collection of all these cones forms a fan.}%
}}}
\put(3556,-3316){\makebox(0,0)[lb]{\smash{\SetFigFont{6}{7.2}{\familydefault}{\mddefault}{\updefault}{\color[rgb]{0,0,0}Two dimensional faces are given by}%
}}}
\put(3556,-3016){\makebox(0,0)[lb]{\smash{\SetFigFont{6}{7.2}{\familydefault}{\mddefault}{\updefault}{\color[rgb]{0,0,0}and $\tau_4=\R_{\ge0}(1,0,1)$.}%
}}}
\put(3556,-2131){\makebox(0,0)[lb]{\smash{\SetFigFont{6}{7.2}{\familydefault}{\mddefault}{\updefault}{\color[rgb]{0,0,0}The zero dimensional face is $0$, and the}%
}}}
\put(3556,-2416){\makebox(0,0)[lb]{\smash{\SetFigFont{6}{7.2}{\familydefault}{\mddefault}{\updefault}{\color[rgb]{0,0,0}one dimensional faces are $\tau_1=\R_{\ge 0}(1,0,0)$,}%
}}}
\put(3556,-3601){\makebox(0,0)[lb]{\smash{\SetFigFont{6}{7.2}{\familydefault}{\mddefault}{\updefault}{\color[rgb]{0,0,0}$\langle\tau_1,\tau_2 \rangle$, $\langle\tau_2,\tau_3\rangle,\langle\tau_3,\tau_4\rangle$ and $\langle\tau_4,\tau_1\rangle$.}%
}}}
\put(2851,-3961){\makebox(0,0)[lb]{\smash{\SetFigFont{6}{7.2}{\familydefault}{\mddefault}{\updefault}{\color[rgb]{0,0,0}$z$}%
}}}
\put(901,-1636){\makebox(0,0)[lb]{\smash{\SetFigFont{6}{7.2}{\familydefault}{\mddefault}{\updefault}{\color[rgb]{0,0,0}$x$}%
}}}
\end{picture}

\end{example}

\begin{definition}
\index{cone}
\index{dual cone}
A {\bf cone of $L$} is a strongly convex rational cone.

If $\sigma$ is a cone then the {\bf dual cone} is defined by
\[
\sigma^* =\{w\in L^*\otimes \qqq
\ |\ \langle v,w\rangle \geq 0,\ 
\forall v\in \sigma\},
\]
where $L^*$ is the dual lattice.

\begin{magma}
Let $n=2$, $L=L^*=\zzz^2$, and suppose
$\sigma=\qqq_{\geq 0}[e_1,3e_1+4e_2]$. 
To check if $e_1^*-7e_2^*$ or if $4e_1^*-3e_2^*$
belongs to $\sigma^*$ in MAGMA
\footnote{The MAGMA 2.8 code used is available on the internet 
at \cite{J2}.},
type
{\footnotesize{
\begin{verbatim}
load "/home/wdj/magmafiles/toric.mag";
//replace /home/wdj/magmafiles by your path to toric.mag
in_dual_cone([1,-7],[[1,0],[3,4]]);
in_dual_cone([4,-3],[[1,0],[3,4]]);
\end{verbatim}
}}
\noindent
In the first case, MAGMA returns \verb+false+
and in the second case, \verb+true+.
\end{magma}

If $\sigma = \qqq_{\geq 0}[v_1,...,v_m]$ is a cone
of $L=\zzz^n$ and if there are vectors
$v_{m+1},...,v_n$ such that
$\det(v_1,...,v_n)=\pm 1$, then we say that
$\sigma$ is {\bf regular}.

Suppose $\sigma, \sigma'$ are cones
of $L=\zzz^n$.
If dim$(\sigma) = $ dim$(\sigma')$ and 
if there is a $g\in GL_n(\zzz)$ for which 
$\sigma' =g\sigma$ then we say that $\sigma$ is
{\bf isomorphic} to $\sigma'$ and write
$\sigma\cong \sigma'$.

\end{definition}

\subsection{Affine toric varieties}

Associate to the dual cone $\sigma^*$ the semigroup 
\[
S_\sigma =\sigma^*\cap L^*
=\{w\in L^*
\ |\ \langle v,w\rangle \geq 0,\ 
\forall v\in \sigma\}.
\]
Though $L^*$ has $n$ generators {\it as a lattice},
typically $S_\sigma$ will have more than $n$
generators {\it as a semigroup}. 
If $u_1,\dots u_t\in L^*$ are semigroup
generators of $S_\sigma$ then we write
\[
S_\sigma=\zzz_{\geq 0}[u_1,...,u_t],
\]
for brevity.

\begin{remark}
The following question arises: Given a lattice
$L=\zzz [v_1,...,v_m]$ and a cone
$\sigma = \qqq_{\geq 0}[w_1,...,w_m]$,
how do you find $u_1,...,u_n\in L^*$ such that
$S_\sigma=\zzz_{\geq 0}[u_1,...,u_t]$?

First, find a basis for the dual lattice, $L^*$,
say $L^*=\zzz [v_1^*,...,v_m^*]$. 
Next, find generators in $L^*$ of the dual cone, say
$\sigma^* = \qqq_{\geq 0}[w^*_1,...,w^*_m]$,
Let the list of
potential generators be
\[
G=\{w\in L^*\ |\ 
|w|\leq {\rm max}_{1\leq i\leq m} |v_i^*|,
\ w\in \sigma^* \}.
\]
It is conjectured that $G$ is a semigroup basis for
$S_\sigma$.

Another related question arises: If $\sigma\subset \qqq^n$
is a (not necessarily convex) cone, and 
$V=\{v_1,...,v_m\}\subset L^*\subset \qqq^n$,
when is $V$ a set of semigroup generators for $S_\sigma$?

It is conjectured that the following two conditions are
sufficient:
\begin{itemize}
\item
$V\subset \sigma$,

\item
$V$ generates each element in 
\[
\{w\in S_\sigma\ |\ |w|\leq 2\cdot {\rm max}_{v\in V} |v|\}.
\]
\end{itemize}

\end{remark}

\begin{magma}
Let $n=2$, $L=L^*=\zzz^2$, and suppose
$\sigma=\qqq_{\geq 0}[e_1,3e_1+4e_2]$. 
To show that 
\[
S_\sigma=\zzz_{\geq 0}[e_1,e_2,2e_1-e_2, 3e_1-2e_2, 4e_1-3e_2],
\]
type
{\footnotesize{
\begin{verbatim}
load "/home/wdj/magmafiles/toric.mag";
//replace /home/wdj/magmafiles by your path to toric.mag
D := LatticeDatabase();
Lat := Lattice(D, 2, 16);
dual_semigp_gens([[1,0],[3,4]],Lat);
\end{verbatim}
}}
\end{magma}

Let
\[
R_\sigma =\cccF[S_\sigma]
\]
denote the ``group algebra'' of this semigroup. It is a
finitely generated commutative $\cccF$-algebra.
It is in fact integrally closed (\cite{F}, page 29).

We may interprete $R_\sigma$ as a subring of 
$R=\cccF[x_1,...,x_n]$ as follows: First, 
identify each $e_i^*$ with the variable $x_i$.
If $S_\sigma$ is generated as a semigroup by 
vectors of the form $\ell_1 e_1^*+...+\ell_n e_n^*$,
where $\ell_i\in \zzz$, then its image in 
$R$ is generated by monomials of the 
form $x_1^{\ell_1}\dots x_n^{\ell_n}$.

Let
\[
U_\sigma={\rm Spec}\ R_\sigma.
\]
This defines an {\bf affine toric variety}
(associated to $\sigma$).
\index{affine toric variety}
\index{toric variety, affine}

\begin{example}
For the cone in Example \ref{ex:3}, we compute the dual cone by 
taking the intersection of the duals of all the one dimensional cones,
and obtain that $\sigma^\vee$ is spanned by
$(0,0,1), (1,-1,0),(1,0,-1)$ and $(0,1,0)$
\[
\begin{picture}(0,0)%
\includegraphics{cones1dual.pstex}%
\end{picture}%
\setlength{\unitlength}{1579sp}%
\begingroup\makeatletter\ifx\SetFigFont\undefined%
\gdef\SetFigFont#1#2#3#4#5{%
  \reset@font\fontsize{#1}{#2pt}%
  \fontfamily{#3}\fontseries{#4}\fontshape{#5}%
  \selectfont}%
\fi\endgroup%
\begin{picture}(4137,3003)(-674,-4483)
\put(451,-4411){\makebox(0,0)[lb]{\smash{\SetFigFont{6}{7.2}{\familydefault}{\mddefault}{\updefault}{\color[rgb]{0,0,0}$y$}%
}}}
\put(2851,-3961){\makebox(0,0)[lb]{\smash{\SetFigFont{6}{7.2}{\familydefault}{\mddefault}{\updefault}{\color[rgb]{0,0,0}$z$}%
}}}
\put(901,-1636){\makebox(0,0)[lb]{\smash{\SetFigFont{6}{7.2}{\familydefault}{\mddefault}{\updefault}{\color[rgb]{0,0,0}$x$}%
}}}
\put(826,-4186){\makebox(0,0)[lb]{\smash{\SetFigFont{6}{7.2}{\familydefault}{\mddefault}{\updefault}{\color[rgb]{0,0,0}$(0,1,0)$}%
}}}
\put(1726,-3886){\makebox(0,0)[lb]{\smash{\SetFigFont{6}{7.2}{\familydefault}{\mddefault}{\updefault}{\color[rgb]{0,0,0}$(0,0,1)$}%
}}}
\put(1876,-2461){\makebox(0,0)[lb]{\smash{\SetFigFont{6}{7.2}{\familydefault}{\mddefault}{\updefault}{\color[rgb]{0,0,0}$(1,-1,0)$}%
}}}
\put(-674,-2911){\makebox(0,0)[lb]{\smash{\SetFigFont{6}{7.2}{\familydefault}{\mddefault}{\updefault}{\color[rgb]{0,0,0}$(1,0,-1)$}%
}}}
\end{picture}
\]
So $U_\sigma$ is given by ${\rm Spec}(\cccF [y,z,x/y,x/z])\cong
{\rm Spec}(\cccF [X,Y,Z,W]/(XZ-WY))$, 
where the isomorphism of rings is given by
\begin{eqnarray*}
X &\mapsto& y\\
Y &\mapsto& z\\
Z &\mapsto& x/y\\
W &\mapsto& x/z
\end{eqnarray*}
So $U_\sigma$ is an affine threefold in $\A^4$
which is a cone over the projective 
quadric surface $XZ=WY$.  There is a singular point at the origin.
\end{example}

\begin{remark}
(For those less familiar with algebraic geometry:)
Here the spectrum of a ring $R$, 
${\rm Spec}\, R$, is the set of prime ideals in $R$
(with a certain topology).  
A basic result we use is that if $f_1,\dots f_m$ are a collection of
polynomials in variables $x_1,\dots,x_d$ with coefficients in $\ccc$, then
\[
{\rm Spec} (\ccc [x_1,\dots,x_d]/(f_1,f_2,\dots,f_m))
\]
is the affine variety in $d$ dimensional space $\aaa^d=F^d$ defined by the
set of points where $f_1=f_2=\cdots=f_m=0$. 
Generally, we will want to find
a way to write an $\ccc$-algebra in this form, so 
that we can understand its Spec geometrically.

Roughly speaking, 
an algebraic variety is given by a collection of affine
pieces $U_1,U_2,\dots,U_d$
which ``glue'' together.  The affine pieces are given by the zero
sets of polynomial equations in some affine spaces
$\aaa^n(\ccc)=\ccc^n$, and the gluings are given
by maps 
\[
\phi_{i,j}:U_i\ra U_j
\]
which are defined by ratios of polynomials on open subsets of the $U_i$.

A major advantage of the toric geometry description is that the relationships
between the affine pieces $U_\sigma$ are simply described 
by the relationships between the
corresponding cones.  For a toric variety this gives an easy way to
keep track of the ``gluing'' data, as we will see in \S \ref{sec:fan}.
\end{remark}

Equivalently, a toric variety is a normal variety $X$ over $\cccF$
which contains a torus $T=(\cccF^\times )^n$ as an open dense 
(in the Zariski topology)
\footnote{The {\bf Zariski topology} 
on an algebraic variety is given by taking closed sets
to be defined by the zero sets of polynomials.
} 
subset:
\[
T\hookrightarrow X
\]
and such that the natural action of $T$ on $T$ extends to an action on 
$X$. In other words, there is a map:
\[
T\times X\ra X,
\]
which restricts to the above map $T\times T\ra T\subset X$.

\begin{example}
Projective space $\P^n$ over $\cccF$
is a toric variety, containing $(\cccF^\times )^n$:
\begin{eqnarray*}
(\cccF^\times )^n&\hra& \P^n\\
(a_1,a_2,\dots a_n)&\mapsto& (1:a_1:a_2:\dots:a_n)
\end{eqnarray*}
The action of $T=(\cccF^\times )^n$ on $\P^n$ is given by
\[
(a_1,a_2,\dots a_n)\cdot (b_0:b_1:b_2:\dots:b_n)
=(b_0:a_1b_1:a_2b_2:\dots :a_nb_n).
\]

\begin{picture}(0,0)%
\includegraphics{P2R.pstex}%
\end{picture}%
\setlength{\unitlength}{1973sp}%
\begingroup\makeatletter\ifx\SetFigFont\undefined%
\gdef\SetFigFont#1#2#3#4#5{%
  \reset@font\fontsize{#1}{#2pt}%
  \fontfamily{#3}\fontseries{#4}\fontshape{#5}%
  \selectfont}%
\fi\endgroup%
\begin{picture}(5247,3344)(-11,-2593)
\put(3346,-2296){\makebox(0,0)[lb]{\smash{\SetFigFont{8}{9.6}{\familydefault}{\mddefault}{\updefault}{\color[rgb]{0,0,0}$X=0$}%
}}}
\put(3211,-1531){\makebox(0,0)[lb]{\smash{\SetFigFont{8}{9.6}{\familydefault}{\mddefault}{\updefault}{\color[rgb]{0,0,0}$Y=0$}%
}}}
\put(406,-2341){\makebox(0,0)[lb]{\smash{\SetFigFont{8}{9.6}{\familydefault}{\mddefault}{\updefault}{\color[rgb]{0,0,0}$Z=0$}%
}}}
\put(4456,359){\makebox(0,0)[lb]{\smash{\SetFigFont{8}{9.6}{\familydefault}{\mddefault}{\updefault}{\color[rgb]{0,0,0}In $\P^2(\R)$ with projective coordinates }%
}}}
\put(4456, 74){\makebox(0,0)[lb]{\smash{\SetFigFont{8}{9.6}{\familydefault}{\mddefault}{\updefault}{\color[rgb]{0,0,0}$X,Y,Z$, we have that $\P^2(\R)\setminus(\R^\times)^2$}%
}}}
\put(4456,-211){\makebox(0,0)[lb]{\smash{\SetFigFont{8}{9.6}{\familydefault}{\mddefault}{\updefault}{\color[rgb]{0,0,0}is given by the three lines $X=0, Y=0$ and $Z=0$.}%
}}}
\put(4456,-781){\makebox(0,0)[lb]{\smash{\SetFigFont{8}{9.6}{\familydefault}{\mddefault}{\updefault}{\color[rgb]{0,0,0}If $(\R^\times)^2$ is embedded by $(a,b)\mapsto (1,a,b)$,}%
}}}
\put(5206,-1636){\makebox(0,0)[lb]{\smash{\SetFigFont{8}{9.6}{\familydefault}{\mddefault}{\updefault}{\color[rgb]{0,0,0}$(a,b)\cdot(0:Y:Z) =(0:aY:bZ)$}%
}}}
\put(5236,-2401){\makebox(0,0)[lb]{\smash{\SetFigFont{8}{9.6}{\familydefault}{\mddefault}{\updefault}{\color[rgb]{0,0,0}$(a,b)\cdot(X:Y:0) = (X:aY:0)$}%
}}}
\put(5221,-2011){\makebox(0,0)[lb]{\smash{\SetFigFont{8}{9.6}{\familydefault}{\mddefault}{\updefault}{\color[rgb]{0,0,0}$(a,b)\cdot(X:0:Y) = (X:0:bZ)$}%
}}}
\put(4456,-1066){\makebox(0,0)[lb]{\smash{\SetFigFont{8}{9.6}{\familydefault}{\mddefault}{\updefault}{\color[rgb]{0,0,0}then the torus action on these lines is given by}%
}}}
\end{picture}

\end{example}

\begin{lemma} 
\label{lemma:isom}
Let $\sigma\subset V=\qqq^n$ be a cone of 
$L=\zzz^n$.
\begin{itemize}
\item
(Fulton \cite{F}, page 29, or
Ewald \cite{E}, ch VI, \S 3)
$U_\sigma$ is smooth if and only if $\sigma$ is a regular cone.

\item
(Ewald \cite{E}, Theorem 2.11, page 222)
$U_{\sigma}\cong U_{\sigma'}$ (as algebraic
varieties) if and only if $\sigma \cong \sigma'$ (as cones).

\item
(Fulton \cite{F}, page 18, Ewald \cite{E}, Theorem 6.1, page 243)
Let $\phi:L'\rightarrow L$ be a homomorphism of lattices
which maps a cone $\sigma'$ of $L$ to a cone $\sigma$ of
$L$ (i.e., $\phi(\sigma')=\sigma$). Then the dual
$\phi_L^*:L^* \rightarrow (L')^*$ gives rise to a map
$\phi_\sigma^*:S_\sigma\rightarrow S_{\sigma'}$.
This determines a map 
$\phi_*:R_\sigma\rightarrow R_{\sigma'}$, and hence
a morphism $\phi^*:U_{\sigma'}\rightarrow U_{\sigma}$.

\end{itemize}

\end{lemma}

A morphism $\phi^*:U_{\sigma'}\rightarrow U_{\sigma}$
arising as in the above lemma from a homomorphism
$\phi:L'\rightarrow L$ be a homomorphism of lattices
will be called an {\bf affine toric morphism}.
\index{affine toric morphism}

Summarizing the construction, we
have the following sequence to determine an
affine toric variety.
Fix a lattice $L$ in $V$. 
{\tiny{
\[
\begin{array}{ccccccc}
\{{\rm rational\ cones}\}
&\rightarrow 
&\{{\rm comm.\ semigroups}\}
&\rightarrow
&\{{\rm semigroup\ algebras}\}
&\rightarrow 
&\{ {\rm affine \ schemes}\} \\
\sigma
& \longmapsto
&S_\sigma = \sigma^*\cap L^*
& \longmapsto
&\cccF [S_\sigma]
& \longmapsto
&U_\sigma={\rm Spec}\ \cccF [S_\sigma].
\end{array}
\]
}}
Using the notion of a ``fan'',
later we shall see how to, in some cases, patch 
these together into a projective version.

\begin{example}
This is really a non-example.

Let $n=1$, so $L=\zzz$, and consider the semigroup $S$
generated by $2,3$. This corresponds to the 
``coordinate ring'' $\cccF [S]\cong
\cccF[x_1^2,x_1^3]$. This is associated to the
curve $x^3=y^2$, which has a cusp at the origin.
The element $x_1$ in the field of
fractions $\cccF(x_1)$ of $\cccF[S]$ is 
integral over $\cccF[S]$ but is not in
$\cccF[S]$, so $\cccF[S]$ is not integrally closed.
In fact, the curve ${\rm Spec}\  \cccF [S]$ is not normal
and is not a toric variety (although it is
``binomial'').

\end{example}

\begin{example} 
\label{ex:hyperboloid}
Let $n=2$, $L=\zzz^2=L^*$, and 
\[
\sigma =\{ae_2+b(2e_1-e_2)\ |\ a,b\geq 0\}.
\]
This may be visualized as follows:

\vskip .5in

\begin{center}
\setlength{\unitlength}{.01cm}

\begin{picture}(256.00,150.00)(-120.00,0.00)
\thicklines
\put(0.00,64.00){\circle*{4}} 
\put(0.00,128.00){\circle*{4}} 

\put(-64.00,0.00){\circle*{4}} 
\put(0.00,-64.00){\circle*{4}} 

\put(64.00,0.00){\circle*{4}} 
\put(128.00,0.00){\circle*{4}} 
\put(192.00,0.00){\circle*{4}} 

\put(0.00,0.00){\vector(0,1){64.00}} 

\put(64.00,64.00){\circle*{4}} 
\put(128.00,0.00){\circle*{4}} 
\put(0.00,128.00){\circle*{4}} 
\put(128.00,64.00){\circle*{4}} 

\put(128.00,64.00){\circle*{4}} 

\put(0.00,-64.00){\circle*{4}} 
\put(64.00,-64.00){\circle*{4}} 
\put(128.00,-64.00){\circle*{4}} 

\put(0.00,0.00){\vector(2,-1){129.00}} 
\end{picture}
\end{center}

\vskip .5in

In this case, $\sigma^*$ is given by
\[
\sigma^* =\{xe_1^* + ye_2^* \ |\ 2x\geq y\geq 0 \}.
\]
This may be visualized as follows:

\vskip .5in

\begin{center}
\setlength{\unitlength}{.01cm}

\begin{picture}(256.00,150.00)(-100.00,0.00)
\thicklines
\put(0.00,64.00){\circle*{4}} 
\put(0.00,128.00){\circle*{4}} 

\put(-64.00,0.00){\circle*{4}} 
\put(0.00,-64.00){\circle*{4}} 

\put(64.00,0.00){\circle*{4}} 
\put(128.00,0.00){\circle*{4}} 
\put(192.00,0.00){\circle*{4}} 

\put(0.00,0.00){\vector(1,0){64.00}} 

\put(64.00,64.00){\circle*{4}} 
\put(128.00,0.00){\circle*{4}} 
\put(0.00,128.00){\circle*{4}} 
\put(128.00,64.00){\circle*{4}} 

\put(128.00,64.00){\circle*{4}} 

\put(0.00,-64.00){\circle*{4}} 
\put(64.00,-64.00){\circle*{4}} 
\put(128.00,-64.00){\circle*{4}} 

\put(0.00,0.00){\vector(1,2){64.00}} 
\end{picture}
\end{center}
\vskip .5in

Therefore the semigroup $S_\sigma=\sigma^* \cap L^*$ is
generated by
$u_1=e_1^*,u_2=e_1^*+e_2^*,u_3=e_1^*+2e_2^*$.
We associate to these generators the monomials
$x_1,x_1x_2,x_1x_2^2$, respectively, since
$me_1^*+ne_2^* \leftrightarrow x_1^mx_2^n$.
This implies
\footnote{See Lemma \ref{lemma:ucondition}
below for a more rigorous approach.}
\[
\cccF [S_\sigma] =\cccF [x_1,x_1x_2,x_1x_2^2].
\]
In fact
\footnote{Again, see Lemma \ref{lemma:ucondition}
below.}
, this ring is isomorphic to
$\cccF [x,y,z]/(xz-y^2)$, via the map
$x\longmapsto x_1,\ y\longmapsto x_1x_2,\ 
z\longmapsto x_1x_2^2$,
so
\[
U_\sigma = {\rm Spec}(\cccF [S_\sigma])
={\rm Spec}(\cccF[x,y,z]/(xz-y^2)),
\]
which is the surface $xz-y^2=0$.
This identity $xz-y^2=0$ corresponds to the vector identity 
$u_1+u_2=2u_3$, which may be easily verified from the picture
above.
\end{example}

How does one (in general) find the 
equation(s) of the toric variety associated to a cone?
One (algebraic) method is to use the following fact.

\begin{lemma}
(\cite{F}, page 19, Exercise) 
If $S_\sigma$ is generated by $u_1,...,u_t$ then
\[
\cccF [S_\sigma]\cong \cccF [\chi^{u_1},...,\chi^{u_t}]
\cong \cccF [y_1,...,y_t]/I,
\]
where $\chi^{e^*_i}=x_i$, $1\leq i\leq t$, and where
$I$ is the ideal generated by binomials
of the form
\[
y_1^{a_1}...y_t^{a_t}-y_1^{b_1}...y_t^{b_t},
\]
where $a_i\geq 0,\ b_j\geq 0$ are integers satisfying
\begin{equation}
a_1u_1+...+a_tu_t = b_1u_1+...+b_tu_t.
\label{eqn:ucondition}
\end{equation}
\label{lemma:ucondition}
\end{lemma}

Using this lemma, the determination of $I$ may be
reduced to a linear algebra problem (over $\zzz$).
For a detailed example of this, see \S 
\ref{sec:example2} below.

Another way to find the 
equation(s) of the toric variety $U_\sigma$  
is to find a complete (but finite) set of semigroup generators
of $S_\sigma$, since each generator of $I$ corresponds to a relation
(\ref{eqn:ucondition}) for the generators of $S_\sigma$.
This provides, at the very least, a nice geometric interpretation of 
the generators of $I$, as we saw at the end of
Example \ref{ex:hyperboloid} above.

\index{torus}
\begin{definition}
A {\bf torus} (of dimension $n$, defined over the field $F$) 
is a product $T=(F^\times)^n$  of $n$ copies of 
$F^\times = F\setminus\{0\}$ for a field $F$.
There is a natural map $T\times T\ra T$
given by pointwise multiplication.  E.g.,
\[
(F^\times)^2\times (F^\times)^2\ra (F^\times)^2
\]
\[
(a,b)\cdot(c,d)=(ac,bd).
\]

\end{definition}

\begin{example}
For another example, take the cone $\sigma$ spanned by
$e_1$ and $-2e_1+3e_2$.
\[
\begin{picture}(0,0)%
\includegraphics{cones2.pstex}%
\end{picture}%
\setlength{\unitlength}{1579sp}%
\begingroup\makeatletter\ifx\SetFigFont\undefined%
\gdef\SetFigFont#1#2#3#4#5{%
  \reset@font\fontsize{#1}{#2pt}%
  \fontfamily{#3}\fontseries{#4}\fontshape{#5}%
  \selectfont}%
\fi\endgroup%
\begin{picture}(4479,1929)(859,-3598)
\put(1651,-1936){\makebox(0,0)[lb]{\smash{\SetFigFont{6}{7.2}{\familydefault}{\mddefault}{\updefault}{\color[rgb]{0,0,0}$\sigma$}%
}}}
\put(4201,-1861){\makebox(0,0)[lb]{\smash{\SetFigFont{6}{7.2}{\familydefault}{\mddefault}{\updefault}{\color[rgb]{0,0,0}$\sigma^*$}%
}}}
\end{picture}
\]
Then $\sigma^*\cup L^*$ is generated by
$e_2,e_1+e_2$, and $3e_1+2e_2$, corresponding to monomials
$X=y,Y=xy,Z=x^3y^2$, so $U_\sigma$ is isomorphic to an affine variety
in $\aaa^3$ given
by equation
\[
Y^3=XZ.
\]
Note, to verify this we just need to check that the map
\begin{eqnarray*}
\cccF [X,Y,Z]/(Y^3-XZ)&\ra& \cccF [y,xy,x^3y^2]\\
X,Y,Z &\mapsto& y,xy,x^3y^2
\end{eqnarray*}
is an isomorphism.

\end{example}

\begin{example}
\label{ex:Rn}
Let $L=L^*=\zzz^n$ be the standard lattice.

\begin{itemize}
\item
Let $\sigma=\{0\}$ be the trivial cone.
Then $\sigma^*=\qqq^n$ is ``everything''
and so the semigroup
$S_\sigma = \sigma^* \cap L^*$ 
is generated by $\pm e_1^*, ..., \pm e_n^*$.
This means that
\[
R_\sigma =\cccF [S_\sigma]\cong 
\cccF [x_1,x_1^{-1},...,x_n,x_n^{-1}].
\]
This is the coordinate ring of the 
torus $(\cccF^\times )^n$, so $U_\sigma =(\cccF^\times )^n $.
\item
Let $\sigma =\{a_1 e_1+..+a_n e_n\ |\ a_i\geq 0\}$.
The dual cone $\sigma^*$ is equal to $\sigma$,
so 
$S_\sigma = \sigma^* \cap L^* 
$
is generated as a semigroup by $e_1^*,...,e_n^*$.
This implies $R_\sigma \cong \cccF[x_1,...,x_n]$,
so $U_\sigma \cong \cccF^n$.

\end{itemize}

\end{example}

\begin{lemma}
(\cite{F}, Proposition, page 29)
A toric variety $U_\tau$ is smooth if and only if
$\tau$ is generated by a subset of a basis for
$L$.
\label{lemma:smooth}
\end{lemma}

\begin{example}
(\cite{F}, page 35)
Let $\sigma$ be the cone generated by $e_1$ and $-e_1+ne_2$.
Then $U_\sigma\cong \cccF^2/\mu_n$,
where $\mu_n=\{z\in \cccF\ |\ z^n=1\}$ and
where $\mu_n$ acts on $\cccF^2$ by
$\zeta:(z_1,z_2)\longmapsto (\zeta z_1,\zeta z_2)$.
\end{example}

\subsection{Faces and subvarieties}

Suppose $\tau$ is a face
\footnote{
Some, like Ewald \cite{E},
use the convention that $\emptyset,\sigma$ are
(improper) faces of $\sigma$. 
We shall define a {\bf face of $\sigma$}
\index{face of $\sigma$}
to be either $\sigma$ itself or
a subset of the form $H\cap \sigma$,
where $H$ is a ``supporting hyperplane
of $\sigma$''
(i.e., a codimension $1$ subspace of
$V$ for which $\sigma\cap H\not= \emptyset$
and $\sigma$ is contained in exactly one
of the two half-spaces determined by
$H$).
} 
of the cone $\sigma$ of $L$.
How is the associated toric variety $U_\tau$
related to $U_\sigma$?

\begin{lemma} 
\label{lemma:subcone}
If $\tau \subset \sigma$ is a face then there is a 
dominant morphism $U_\tau \rightarrow U_\sigma$.
\end{lemma}

\pf  
If $\tau $ is a face of $\sigma$ then there is a
$u\in S_\sigma =\sigma^* \cap L^*$ such that
$\tau = u^\perp \cap \sigma$, where
$u^\perp =\{v \in V\ |\ u\cdot v =0\}$. We may regard 
$R_\sigma =\cccF [S_\sigma]$ as the coordinate ring
of $U_\sigma$. As such,
\[
U_\tau = \{x\in U_\sigma \ |\ u(x)\not= 0\}.
\]
It is clear that the ``inclusion'' map
$U_\tau\rightarrow U_\sigma$ is surjective 
on a Zariski dense set, so it is
dominant.  \qed 

\begin{example} 
Let $n=2$, $L=\zzz^2=L^*$, and 
\[
\sigma =\{ae_2+b(2e_1-e_2)\ |\ a,b\geq 0\}
\]
and let $\tau =\qqq_{\geq 0}\cdot (2e_1-e_2)$.
In fact, if we take $u=e_1^*+2e_2^*$ then
$\tau =u^\perp \cap \sigma$. If 
we represent $x\in U_\sigma$ by 
the point $x_1e_1+x_2e_2=(x_1,x_2)$ then 
$u(x)=x_1x_2^2$ (or $x_1+2x_2$, but this will instead be written
multiplicatively, as $x_1x_2^2$).
 
The dual cone of $\tau$ is a half-plane,
\[
\tau^* =\{ xe_1^* +ye_2^* \ |\ 2x-y\geq 0\},
\]
so $S_\tau =\tau^* \cap L^*$ is generated
as a semigroup by $e_1^*$, $-e_2^*$, $e_1^*+2e_2^*$,
$-e_1^*-2e_2^*$. The coordinate 
ring associate to $\tau$ is given by
\[
R_\tau
=\cccF[e_1^*, -e_2^*,e_1^*+2e_2^*,-e_1^*-2e_2^*]
\cong \cccF[x_1,x_2^{-1},x_1x_2^2,x_1^{-1}x_2^{-2}].
\]
By inspection, we have an isomorphism
\[
R_\tau \cong \cccF [x,y,z,w]/(wz^2-x,xy-z^2),
\]
via $x=x_1$, $y=x_1^{-1}x_2^{-2}$, $z=x_2^{-1}$,
$w=x_1x_2^2$. This implies that 
$U_\tau$ is the variety in $\cccF^4$ defined
by
\[
wz^2=x,\ \ \ \ \ xy=z^2.
\]
Presumably, we may embed the variety $U_\sigma$
given by $xy=z^2$ determined in Example \ref{ex:hyperboloid}
into $x,y,z,w$-space and regard $U_\tau$ as the
dense subvariety $x\not= 0$ of this.
(This condition implies $w\not= 0$.)

\end{example}

\begin{magma}

Let $n=2$, $L=L^*=\zzz^2$, and suppose
$\sigma=\qqq_{\geq 0}[e_1,3e_1+4e_2]$. 
To find the ideal $I$ defining
the quotient ring $R_\sigma = \cccF [x_1,...,x_5]/I$, type

{\footnotesize{
\begin{verbatim}
load "/home/wdj/magmafiles/toric.mag";
//replace /home/wdj/magmafiles by your path to toric.mag
D := LatticeDatabase();
Lat := Lattice(D, 2, 16);
ideal_affine_toric_variety([[1,0],[3,4]],Lat);
\end{verbatim}
}}
\noindent
MAGMA returns
{\footnotesize{
\begin{verbatim}
Ideal of Polynomial ring of rank 5 over Rational Field
Lexicographical Order
Variables: x1, x2, x3, x4, x5
Basis:
[
    x1*x5^3 - x4^4,
    x2*x5^2 - x4^3
]
\end{verbatim}
}}
\end{magma}

\subsection{The dense torus}

Let $\sigma\subset \qqq^n$ be a cone in $L$.
In the case $\tau=\{0\}\subset \sigma$, the
dominant morphism
$U_\tau\rightarrow U_\sigma$ maps a torus
$U_\tau\cong \cccF^{\times n}$ into a dense
subvariety of $U_\sigma$. We shall briefly
recall another way to view this.

Let $\ggg_m$ denote the multiplicative algebraic group
over $\cccF$. For each integer $k$, the map
$z\longmapsto z^k$ defines an element of
Hom$_{alg.gp.}(\ggg_m,\ggg_m)$. In fact, each element of
Hom$_{alg.gp.}(\ggg_m,\ggg_m)$ arises in this way. 
Given a lattice $L\subset \qqq^n$, let
$T_L=$Hom$_{ab.gp.}(L^*,\ggg_m)$. By Lemma
\ref{lemma:ucondition}, we have
$T_L\subset U_\sigma$. They have the same dimension,
so $T_L$ must be dense in $U_\sigma$.
For an explicit example, see
Example \ref{ex:3.7} below (where $T_L=U_{\sigma_6}$).

\begin{magma}
\label{magma:embed}
To define a torus of dimension (for example) $5$ in MAGMA, as a scheme 
over $\qqq$, type
{\footnotesize{
\begin{verbatim}
load "/home/wdj/magmafiles/toric.mag";
//replace /home/wdj/magmafiles by your path to toric.mag
create_torus(5);
\end{verbatim}
}}
\noindent
The variables \verb+x1, x2, x3, x4, x6, x7, x8, x9, x10+
satisfy $x_i x_{i+5}=1$, for $1\leq i \leq 5$. 

Let $n=2$, $L=L^*=\zzz^2$, and suppose
$\sigma=\qqq_{\geq 0}[e_1,3e_1+4e_2]$. 
To obtain the rational map (of schemes)
$T_L\rightarrow U_\sigma$, type
{\footnotesize{
\begin{verbatim}
embedding_affine_toric_variety([[1,0],[3,4]]);
\end{verbatim}
}}
\noindent
MAGMA returns
{\footnotesize{
\begin{verbatim}
Mapping from: Scheme over Rational Field defined by
x1*x3 - 1
x2*x4 - 1 to Affine Space of dimension 5
Variables : $.1, $.2, $.3, $.4, $.5 %$
with equations : 
x2
x1
x1^2*x4
x1^3*x4^2
x1^4*x4^3
\end{verbatim}
}}
\noindent
This means that $U_\sigma$ is $2$-dimensional
and the dense embedding 
given by
\[
\begin{array}{ccc}
(\cccF^\times )^2&\rightarrow &U_\sigma,\\
(x_1,x_2)&\longmapsto & (x_2,x_1,x_1^2x_4, x_1^3x_4^2, x_1^4x_4^3)
=(x_2,x_1,x_1^2x_2^{-1}, x_1^3x_2^{-2}, x_1^4x_2^{-3}).
\end{array}
\]
\end{magma}

\begin{example}
Another ``toric compactification'' of $(\cccF^\times )^n$ is given by
taking the product of $n$ copies of $\P^1$.  There is a map
\begin{eqnarray*}
(\cccF^\times )^n
& \hra & \overbrace{\P^1\times\P^1\dots \P^1}^{n\ {\rm copies}}\\
(a_1,a_2,\dots,a_n)&\mapsto &
(1:a_1)\times(1:a_2)\times\cdots\times(1:a_n)
\end{eqnarray*}

\begin{picture}(0,0)%
\includegraphics{P1xP1R.pstex}%
\end{picture}%
\setlength{\unitlength}{1973sp}%
\begingroup\makeatletter\ifx\SetFigFont\undefined%
\gdef\SetFigFont#1#2#3#4#5{%
  \reset@font\fontsize{#1}{#2pt}%
  \fontfamily{#3}\fontseries{#4}\fontshape{#5}%
  \selectfont}%
\fi\endgroup%
\begin{picture}(5187,3624)(-11,-2773)
\put(4456,359){\makebox(0,0)[lb]{\smash{\SetFigFont{8}{9.6}{\familydefault}{\mddefault}{\updefault}{\color[rgb]{0,0,0}In $\P^1(\R)\times\P^1(\R)$ with projective coordinates }%
}}}
\put(2776,-886){\makebox(0,0)[lb]{\smash{\SetFigFont{8}{9.6}{\familydefault}{\mddefault}{\updefault}{\color[rgb]{0,0,0}$X=0$}%
}}}
\put(961,-946){\makebox(0,0)[lb]{\smash{\SetFigFont{8}{9.6}{\familydefault}{\mddefault}{\updefault}{\color[rgb]{0,0,0}$Y=0$}%
}}}
\put(1366,239){\makebox(0,0)[lb]{\smash{\SetFigFont{8}{9.6}{\familydefault}{\mddefault}{\updefault}{\color[rgb]{0,0,0}$Z=0$}%
}}}
\put(1396,-1801){\makebox(0,0)[lb]{\smash{\SetFigFont{8}{9.6}{\familydefault}{\mddefault}{\updefault}{\color[rgb]{0,0,0}$W=0$}%
}}}
\put(4456, 74){\makebox(0,0)[lb]{\smash{\SetFigFont{8}{9.6}{\familydefault}{\mddefault}{\updefault}{\color[rgb]{0,0,0}$X,Y$ and $Z,W$ for the two components, we have}%
}}}
\put(4456,-211){\makebox(0,0)[lb]{\smash{\SetFigFont{8}{9.6}{\familydefault}{\mddefault}{\updefault}{\color[rgb]{0,0,0}that $\P^1(\R)\times\P^1(\R)\setminus(\R^\times)^2$ is given by $4$ lines,}%
}}}
\put(4456,-496){\makebox(0,0)[lb]{\smash{\SetFigFont{8}{9.6}{\familydefault}{\mddefault}{\updefault}{\color[rgb]{0,0,0}$X=0, Y=0, Z=0$ and $W=0$}%
}}}
\put(4456,-1066){\makebox(0,0)[lb]{\smash{\SetFigFont{8}{9.6}{\familydefault}{\mddefault}{\updefault}{\color[rgb]{0,0,0}If the embedding is given by $(a,b)\mapsto (1:a)\times (1:b)$}%
}}}
\put(4456,-1351){\makebox(0,0)[lb]{\smash{\SetFigFont{8}{9.6}{\familydefault}{\mddefault}{\updefault}{\color[rgb]{0,0,0}then the torus action on $Y=0$ is given by}%
}}}
\put(5176,-1891){\makebox(0,0)[lb]{\smash{\SetFigFont{8}{9.6}{\familydefault}{\mddefault}{\updefault}{\color[rgb]{0,0,0}$(a,b)\cdot((1:0)\times (Z:W)) = (1:0)\times (Z,bW)$}%
}}}
\end{picture}

\end{example}

The following picture gives a topological view of the situation:

\[
\begin{picture}(0,0)%
\includegraphics{topology.pstex}%
\end{picture}%
\setlength{\unitlength}{1579sp}%
\begingroup\makeatletter\ifx\SetFigFont\undefined%
\gdef\SetFigFont#1#2#3#4#5{%
  \reset@font\fontsize{#1}{#2pt}%
  \fontfamily{#3}\fontseries{#4}\fontshape{#5}%
  \selectfont}%
\fi\endgroup%
\begin{picture}(9341,6000)(301,-5776)
\put(526,-5386){\makebox(0,0)[lb]{\smash{\SetFigFont{6}{7.2}{\familydefault}{\mddefault}{\updefault}{\color[rgb]{0,0,0}(only a finite amount is shaded)}%
}}}
\put(5506,-661){\makebox(0,0)[lb]{\smash{\SetFigFont{6}{7.2}{\familydefault}{\mddefault}{\updefault}{\color[rgb]{0,0,0}$\P^2(\R)$}%
}}}
\put(6136,-4426){\makebox(0,0)[lb]{\smash{\SetFigFont{6}{7.2}{\familydefault}{\mddefault}{\updefault}{\color[rgb]{0,0,0}$\P^1(\R)\times\P^1(\R)$}%
}}}
\put(586,-76){\makebox(0,0)[lb]{\smash{\SetFigFont{6}{7.2}{\familydefault}{\mddefault}{\updefault}{\color[rgb]{0,0,0}Topological picture of}%
}}}
\put(586,-361){\makebox(0,0)[lb]{\smash{\SetFigFont{6}{7.2}{\familydefault}{\mddefault}{\updefault}{\color[rgb]{0,0,0}toric varieties as completions of}%
}}}
\put(586,-646){\makebox(0,0)[lb]{\smash{\SetFigFont{6}{7.2}{\familydefault}{\mddefault}{\updefault}{\color[rgb]{0,0,0}$(\R^\times)^2$}%
}}}
\put(301,-1486){\makebox(0,0)[lb]{\smash{\SetFigFont{6}{7.2}{\familydefault}{\mddefault}{\updefault}{\color[rgb]{0,0,0}Four connected components of $(\R^\times)^2$}%
}}}
\put(5251, 44){\makebox(0,0)[lb]{\smash{\SetFigFont{6}{7.2}{\familydefault}{\mddefault}{\updefault}{\color[rgb]{0,0,0}Two of the possible compactifications:}%
}}}
\put(8251,-1261){\makebox(0,0)[lb]{\smash{\SetFigFont{6}{7.2}{\familydefault}{\mddefault}{\updefault}{\color[rgb]{0,0,0}$Y=0$}%
}}}
\put(7801,-661){\makebox(0,0)[lb]{\smash{\SetFigFont{6}{7.2}{\familydefault}{\mddefault}{\updefault}{\color[rgb]{0,0,0}$X=0$}%
}}}
\put(6451,-2461){\makebox(0,0)[lb]{\smash{\SetFigFont{6}{7.2}{\familydefault}{\mddefault}{\updefault}{\color[rgb]{0,0,0}$Z=0$}%
}}}
\put(9001,-2011){\makebox(0,0)[lb]{\smash{\SetFigFont{6}{7.2}{\familydefault}{\mddefault}{\updefault}{\color[rgb]{0,0,0}edge identified }%
}}}
\put(9001,-2266){\makebox(0,0)[lb]{\smash{\SetFigFont{6}{7.2}{\familydefault}{\mddefault}{\updefault}{\color[rgb]{0,0,0}diagonally}%
}}}
\put(6451,-2716){\makebox(0,0)[lb]{\smash{\SetFigFont{6}{7.2}{\familydefault}{\mddefault}{\updefault}{\color[rgb]{0,0,0}``line at infinity''}%
}}}
\end{picture}
\]

\section{Fans and toric varieties}
\label{sec:fan}

\subsection{Fans}

A fan is a collection of cones which ``fit
together'' well.

\index{fan}
\begin{definition}
A {\bf fan in $L$} is a set $\Delta=\{\sigma\}$
of rational strongly convex
polyhedral cones in $L_{\qqq}=L\otimes \qqq$
such that
\begin{itemize}
\item
if $\sigma\in \Delta$ and $\tau\subset \sigma$ is a face
of $\sigma$ then $\tau \in \Delta$,
\item
if $\sigma_1,\sigma_2\in \Delta$ then 
$\sigma_1\cap \sigma_2$ is a face of 
both $\sigma_1$ and $\sigma_2$ (and hence
belongs to $\Delta$ by the above).
\end{itemize}

If $V=\cup_{\sigma\in \Delta} \sigma$ then
we call the fan {\bf complete}.
\end{definition}
\index{fan, complete}

We shall assume that all fans are finite.

\begin{example}
Some fans in $\R^2$:
\[
\begin{picture}(0,0)%
\includegraphics{fans.pstex}%
\end{picture}%
\setlength{\unitlength}{1579sp}%
\begingroup\makeatletter\ifx\SetFigFont\undefined%
\gdef\SetFigFont#1#2#3#4#5{%
  \reset@font\fontsize{#1}{#2pt}%
  \fontfamily{#3}\fontseries{#4}\fontshape{#5}%
  \selectfont}%
\fi\endgroup%
\begin{picture}(7694,2744)(729,-2783)
\end{picture}
\]
In the picture $2$-dimensional cones are partially shaded.
\end{example}

\begin{example}
Consider the fan given by the two dimensional cone $\sigma$
spanned by $e_1$ and $e_1+2e_2$, and its faces,
where $e_1,e_2$ is the standard basis for the lattice $N\cong\zzz^2$.
\[
\begin{picture}(0,0)%
\includegraphics{sigma.pstex}%
\end{picture}%
\setlength{\unitlength}{1579sp}%
\begingroup\makeatletter\ifx\SetFigFont\undefined%
\gdef\SetFigFont#1#2#3#4#5{%
  \reset@font\fontsize{#1}{#2pt}%
  \fontfamily{#3}\fontseries{#4}\fontshape{#5}%
  \selectfont}%
\fi\endgroup%
\begin{picture}(4542,2470)(1159,-3719)
\put(2701,-2161){\makebox(0,0)[lb]{\smash{\SetFigFont{6}{7.2}{\familydefault}{\mddefault}{\updefault}{\color[rgb]{0,0,0}$\sigma$}%
}}}
\put(4951,-1711){\makebox(0,0)[lb]{\smash{\SetFigFont{6}{7.2}{\familydefault}{\mddefault}{\updefault}{\color[rgb]{0,0,0}$\sigma^*$}%
}}}
\put(3301,-3661){\makebox(0,0)[lb]{\smash{\SetFigFont{6}{7.2}{\familydefault}{\mddefault}{\updefault}{\color[rgb]{0,0,0}$L^*$}%
}}}
\put(1201,-2011){\makebox(0,0)[lb]{\smash{\SetFigFont{6}{7.2}{\familydefault}{\mddefault}{\updefault}{\color[rgb]{0,0,0}$L$}%
}}}
\put(5701,-2161){\makebox(0,0)[lb]{\smash{\SetFigFont{6}{7.2}{\familydefault}{\mddefault}{\updefault}{\color[rgb]{0,0,0}Generators are}%
}}}
\put(5701,-2446){\makebox(0,0)[lb]{\smash{\SetFigFont{6}{7.2}{\familydefault}{\mddefault}{\updefault}{\color[rgb]{0,0,0}circled}%
}}}
\put(4051,-3211){\makebox(0,0)[lb]{\smash{\SetFigFont{5}{6.0}{\familydefault}{\mddefault}{\updefault}{\color[rgb]{0,0,0}$0$}%
}}}
\put(1276,-3286){\makebox(0,0)[lb]{\smash{\SetFigFont{5}{6.0}{\familydefault}{\mddefault}{\updefault}{\color[rgb]{0,0,0}$0$}%
}}}
\end{picture}
\]
Let $e_i^*$ be the dual basis for $L^*$.  Then we have
\[
\sigma^\vee\cup L^*=\langle e_1,e_2,2e_1-e_2\rangle,
\]
So
\[
\cccF [S_\sigma]=\cccF [x,y,x^2/y],
\]
where $x=e^{e_1}$ and $y=e^{e_2}$.
If we set $X=x, Y=y$ and $Z=x^2/y$, then we have
\[
\cccF [S_\sigma]=\cccF [X,Y,Z]/(X^2 - ZY).
\]
So $U_\sigma$ is a quadric cone in $\aaa^2$ given by $X^2-ZY$.
\[
\begin{picture}(0,0)%
\includegraphics{quadric.pstex}%
\end{picture}%
\setlength{\unitlength}{1579sp}%
\begingroup\makeatletter\ifx\SetFigFont\undefined%
\gdef\SetFigFont#1#2#3#4#5{%
  \reset@font\fontsize{#1}{#2pt}%
  \fontfamily{#3}\fontseries{#4}\fontshape{#5}%
  \selectfont}%
\fi\endgroup%
\begin{picture}(7452,3835)(151,-4919)
\put(151,-1651){\makebox(0,0)[lb]{\smash{\SetFigFont{6}{7.2}{\familydefault}{\mddefault}{\updefault}{\color[rgb]{0,0,0}$X^2=YZ$ in $\A^3$}%
}}}
\put(151,-1276){\makebox(0,0)[lb]{\smash{\SetFigFont{6}{7.2}{\familydefault}{\mddefault}{\updefault}{\color[rgb]{0,0,0}$U_\sigma$ is given by}%
}}}
\put(1606,-2191){\makebox(0,0)[lb]{\smash{\SetFigFont{6}{7.2}{\familydefault}{\mddefault}{\updefault}{\color[rgb]{0,0,0}$X$}%
}}}
\put(3511,-3286){\makebox(0,0)[lb]{\smash{\SetFigFont{6}{7.2}{\familydefault}{\mddefault}{\updefault}{\color[rgb]{0,0,0}$Y$}%
}}}
\put(3316,-4861){\makebox(0,0)[lb]{\smash{\SetFigFont{6}{7.2}{\familydefault}{\mddefault}{\updefault}{\color[rgb]{0,0,0}$Z$}%
}}}
\put(5461,-2191){\makebox(0,0)[lb]{\smash{\SetFigFont{6}{7.2}{\familydefault}{\mddefault}{\updefault}{\color[rgb]{0,0,0}$X$}%
}}}
\put(7366,-3286){\makebox(0,0)[lb]{\smash{\SetFigFont{6}{7.2}{\familydefault}{\mddefault}{\updefault}{\color[rgb]{0,0,0}$Y$}%
}}}
\put(7171,-4861){\makebox(0,0)[lb]{\smash{\SetFigFont{6}{7.2}{\familydefault}{\mddefault}{\updefault}{\color[rgb]{0,0,0}$Z$}%
}}}
\put(4006,-1291){\makebox(0,0)[lb]{\smash{\SetFigFont{6}{7.2}{\familydefault}{\mddefault}{\updefault}{\color[rgb]{0,0,0}The torus in $U_\sigma$ is given by}%
}}}
\put(4006,-1576){\makebox(0,0)[lb]{\smash{\SetFigFont{6}{7.2}{\familydefault}{\mddefault}{\updefault}{\color[rgb]{0,0,0}the part where $YZ\not=0$}%
}}}
\end{picture}
\]
For one $1$ dimensional face of $\sigma$ we have
\[
\langle e_1+2e_2\rangle^*\cap L^* = \langle 2e_1-e_2,-2e_1+e_2,e_2\rangle,
\]
so 
\begin{eqnarray*}
U_{\langle e_1+2e_2\rangle}&=&{\rm Spec}(\cccF [y,x^2/y,y/x^2])\\
&=&{\rm Spec}(\cccF [X,Y,Z,Z^{-1}]/(X^2-ZY))
=U_\sigma\setminus\{Z=0\}
\end{eqnarray*}
For the other $1$ dimensional cone we have
\[
U_{\langle e_1\rangle}={\rm Spec}(\cccF [x,y,1/y])
=U_\sigma\setminus\{Y=0\}.
\]
For the zero dimensional cone, we get the embedding
\[
(a,b)\mapsto (a,b,a^2/b)\in U_\sigma .
\]

\end{example}

\begin{example}
For the fan $\Delta$
in $\R$ consisting of $3$ cones, $\tau_0 = \{0\}$, $\tau_1 = \R_{\ge 0}$
and $\tau_2 = \R_{\le0}$, 
\[
\begin{picture}(0,0)%
\includegraphics{P1R.pstex}%
\end{picture}%
\setlength{\unitlength}{1579sp}%
\begingroup\makeatletter\ifx\SetFigFont\undefined%
\gdef\SetFigFont#1#2#3#4#5{%
  \reset@font\fontsize{#1}{#2pt}%
  \fontfamily{#3}\fontseries{#4}\fontshape{#5}%
  \selectfont}%
\fi\endgroup%
\begin{picture}(7875,3625)(-1649,-2894)
\put(1951,539){\makebox(0,0)[lb]{\smash{\SetFigFont{6}{7.2}{\familydefault}{\mddefault}{\updefault}{\color[rgb]{0,0,0}The corresponding variety $X(\Delta):$}%
}}}
\put(-1649,239){\makebox(0,0)[lb]{\smash{\SetFigFont{6}{7.2}{\familydefault}{\mddefault}{\updefault}{\color[rgb]{0,0,0}The fan $\Delta$, consisting}%
}}}
\put(-1649,-46){\makebox(0,0)[lb]{\smash{\SetFigFont{6}{7.2}{\familydefault}{\mddefault}{\updefault}{\color[rgb]{0,0,0}of $3$ cones:}%
}}}
\put(6226,-286){\makebox(0,0)[lb]{\smash{\SetFigFont{6}{7.2}{\familydefault}{\mddefault}{\updefault}{\color[rgb]{0,0,0}projective}%
}}}
\put(6226,-586){\makebox(0,0)[lb]{\smash{\SetFigFont{6}{7.2}{\familydefault}{\mddefault}{\updefault}{\color[rgb]{0,0,0}coordinates}%
}}}
\put(6226,-886){\makebox(0,0)[lb]{\smash{\SetFigFont{6}{7.2}{\familydefault}{\mddefault}{\updefault}{\color[rgb]{0,0,0}$X,Y$}%
}}}
\put(6226, 14){\makebox(0,0)[lb]{\smash{\SetFigFont{6}{7.2}{\familydefault}{\mddefault}{\updefault}{\color[rgb]{0,0,0}$\P^1(\R)\cong S^1$}%
}}}
\put(4651,-736){\makebox(0,0)[lb]{\smash{\SetFigFont{6}{7.2}{\familydefault}{\mddefault}{\updefault}{\color[rgb]{0,0,0}$Y=0$}%
}}}
\put(4576,-2836){\makebox(0,0)[lb]{\smash{\SetFigFont{6}{7.2}{\familydefault}{\mddefault}{\updefault}{\color[rgb]{0,0,0}$X=0$}%
}}}
\put(4351,-1411){\makebox(0,0)[lb]{\smash{\SetFigFont{6}{7.2}{\familydefault}{\mddefault}{\updefault}{\color[rgb]{0,0,0}$U_{\tau_1}$}%
}}}
\put(6226,-2446){\makebox(0,0)[lb]{\smash{\SetFigFont{6}{7.2}{\familydefault}{\mddefault}{\updefault}{\color[rgb]{0,0,0}local variable}%
}}}
\put(6226,-2731){\makebox(0,0)[lb]{\smash{\SetFigFont{6}{7.2}{\familydefault}{\mddefault}{\updefault}{\color[rgb]{0,0,0}$y=Y/X$}%
}}}
\put(6226,-2086){\makebox(0,0)[lb]{\smash{\SetFigFont{6}{7.2}{\familydefault}{\mddefault}{\updefault}{\color[rgb]{0,0,0}$U_{\tau_2}$}%
}}}
\put(751,-1411){\makebox(0,0)[lb]{\smash{\SetFigFont{6}{7.2}{\familydefault}{\mddefault}{\updefault}{\color[rgb]{0,0,0}$\tau_1$}%
}}}
\put(-974,-1411){\makebox(0,0)[lb]{\smash{\SetFigFont{6}{7.2}{\familydefault}{\mddefault}{\updefault}{\color[rgb]{0,0,0}$\tau_2$}%
}}}
\put( 76,-586){\makebox(0,0)[lb]{\smash{\SetFigFont{6}{7.2}{\familydefault}{\mddefault}{\updefault}{\color[rgb]{0,0,0}$\tau_0$}%
}}}
\end{picture}
\]
we have
\begin{eqnarray*}
U_{\tau_1}&=&{\rm Spec}(\cccF [x])\cong\aaa^1\\
U_{\tau_2}&=&{\rm Spec}(\cccF [x^{-1}])\cong\aaa^1\\
U_{\tau_0}&=&{\rm Spec}(\cccF [x,x^{-1}])\cong \cccF^\times
\end{eqnarray*}
So $X(\Delta)\cong \P^1$, with projective coordinates $X,Y$, and
\begin{eqnarray*}
U_{\tau_1}&=&{\rm Spec}(\cccF [x])\cong\P^1\setminus Y=0\\
U_{\tau_2}&=&{\rm Spec}(\cccF [x^{-1}])\cong\P^1\setminus X=0\\
\end{eqnarray*}
Where the map between $\P^1$ and $U_{\tau_1}$ is given by $x=X/Y$,
and the map between $\P^1$ and $U_{\tau_2}$ is given by $x^{-1}=Y/X$.

\end{example}

\begin{example}
$\P^1\times\P^1$ is given by the following fan:
\[
\begin{picture}(0,0)%
\includegraphics{P1P1.pstex}%
\end{picture}%
\setlength{\unitlength}{1579sp}%
\begingroup\makeatletter\ifx\SetFigFont\undefined%
\gdef\SetFigFont#1#2#3#4#5{%
  \reset@font\fontsize{#1}{#2pt}%
  \fontfamily{#3}\fontseries{#4}\fontshape{#5}%
  \selectfont}%
\fi\endgroup%
\begin{picture}(4812,2352)(451,-3598)
\put(451,-1711){\makebox(0,0)[lb]{\smash{\SetFigFont{6}{7.2}{\familydefault}{\mddefault}{\updefault}{\color[rgb]{0,0,0}are the four quadrants:}%
}}}
\put(3811,-1411){\makebox(0,0)[lb]{\smash{\SetFigFont{6}{7.2}{\familydefault}{\mddefault}{\updefault}{\color[rgb]{0,0,0}Cones dual to the}%
}}}
\put(3796,-1696){\makebox(0,0)[lb]{\smash{\SetFigFont{6}{7.2}{\familydefault}{\mddefault}{\updefault}{\color[rgb]{0,0,0}two dimensional cones:}%
}}}
\put(451,-1411){\makebox(0,0)[lb]{\smash{\SetFigFont{6}{7.2}{\familydefault}{\mddefault}{\updefault}{\color[rgb]{0,0,0}The $2$ dimensional cones}%
}}}
\end{picture}
\]
The whole situation is the ``product'' of the situation for
$\P^1$ above.

More generally, for the fan $\Delta_n$
given by taking $n$ dimensional cones
spanned by $(\pm e_1,\pm e_2,\dots,\pm e_n)$ in $\R^n$, the 
toric variety $X(\Delta_n)$ is isomorphic to 
the product of $n$ copies of $\P_1$, i.e., $\P^1\times\dots \P_1$.
\end{example}

\begin{example}
$\P^2$ is given by the following fan:
\[
\begin{picture}(0,0)%
\includegraphics{P2.pstex}%
\end{picture}%
\setlength{\unitlength}{1579sp}%
\begingroup\makeatletter\ifx\SetFigFont\undefined%
\gdef\SetFigFont#1#2#3#4#5{%
  \reset@font\fontsize{#1}{#2pt}%
  \fontfamily{#3}\fontseries{#4}\fontshape{#5}%
  \selectfont}%
\fi\endgroup%
\begin{picture}(5100,2727)(301,-3883)
\put(301,-1336){\makebox(0,0)[lb]{\smash{\SetFigFont{6}{7.2}{\familydefault}{\mddefault}{\updefault}{\color[rgb]{0,0,0}three two dimensional cones}%
}}}
\put(3226,-1561){\makebox(0,0)[lb]{\smash{\SetFigFont{6}{7.2}{\familydefault}{\mddefault}{\updefault}{\color[rgb]{0,0,0}the dual cones}%
}}}
\put(5401,-2461){\makebox(0,0)[lb]{\smash{\SetFigFont{6}{7.2}{\familydefault}{\mddefault}{\updefault}{\color[rgb]{0,0,0}$\Spec(k[x^{-1},yx^{-1}])$}%
}}}
\put(5401,-1711){\makebox(0,0)[lb]{\smash{\SetFigFont{6}{7.2}{\familydefault}{\mddefault}{\updefault}{\color[rgb]{0,0,0}$\Spec(k[x,y])$}%
}}}
\put(5176,-1336){\makebox(0,0)[lb]{\smash{\SetFigFont{6}{7.2}{\familydefault}{\mddefault}{\updefault}{\color[rgb]{0,0,0}Corresponding affine pieces:}%
}}}
\put(5401,-2086){\makebox(0,0)[lb]{\smash{\SetFigFont{6}{7.2}{\familydefault}{\mddefault}{\updefault}{\color[rgb]{0,0,0}$\Spec(k[y^{-1},xy^{-1}])$}%
}}}
\put(601,-2161){\makebox(0,0)[lb]{\smash{\SetFigFont{6}{7.2}{\familydefault}{\mddefault}{\updefault}{\color[rgb]{0,0,0}$\sigma_2$}%
}}}
\put(1351,-3811){\makebox(0,0)[lb]{\smash{\SetFigFont{6}{7.2}{\familydefault}{\mddefault}{\updefault}{\color[rgb]{0,0,0}$\sigma_3$}%
}}}
\put(2026,-1936){\makebox(0,0)[lb]{\smash{\SetFigFont{6}{7.2}{\familydefault}{\mddefault}{\updefault}{\color[rgb]{0,0,0}$\sigma_1$}%
}}}
\end{picture}
\]
If $X,Y,Z$ are the projective coordinates for $\P^2$, then we can
identify the affine pieces with the pieces with local variables
\[
\sigma_1: x=X/Z, y =Y/Z,\]
\[
\sigma_2: x^{-1}=Z/X, yx^{-1}=Y/X,
\]
\[
\sigma_3: y^{-1}=Z/Y, xy^{-1}=X/Y,
\]
which are the natural coordinates for the $Z\not=0$, $X\not=0$ and
$Y\not=0$ parts of $\P^2$ respectively.

More generally, $\P^n$ can be given as a toric variety.  It is
``prettier'' to take the lattice to be given by the sublattice 
of $\Z^{n+1}$ given by
\[
L=\{(a_0,a_1,a_2,\dots,a_n)\in\Z^{n+1}|\sum a_i =0 \},
\]
and to take cones $\sigma_j$ spanned by $e_i - e_j$ for $j=0,\dots n$.
So that $U_{\sigma_j}$ has local coordinates $x_i/x_j$, where
$\prod_{k=0}^n x_k=1$.  If $\P^n$ is a projective variety with
projective coordinates $X_i$, then 
$U_\sigma$ can be identified with the part where $X_j\not=0$, via a map
$x_i/x_j\mapsto X_i/X_j$.
\end{example}

\begin{remark}
Some general remarks on how fans might be implemented on a computer.

A fan is a set of cones satisfying certain conditions. Each cone
$\sigma$ is specified by a set
\footnote{The empty set $\emptyset$ corresponds to the zero cone $\{0\}$.} 
of vectors $\{v_1,...,v_k\}\subset L$:
\[
\begin{array}{c}
\sigma \leftrightarrow \{v_1,...,v_k\}
\\
\sigma =\qqq_{\geq 0}[v_1,...,v_k].
\end{array}
\]
The dimension of $\sigma$ is the dimension of the vector space
span, ${\rm span}_{\qqq}\{v_1,...,v_k\}$.

Suppose
$\sigma \leftrightarrow \{v_1,...,v_k\}$ and
$\sigma' \leftrightarrow \{v'_1,...,v'_{k'}\}$ belong to a fan $\Delta$.
We conjecture that if the generators $v_i,v'_j\in L$ of $\sigma$ 
are choosen to have minimum length then
\[
\sigma\cap \sigma' \leftrightarrow 
\{v_1,...,v_k\}\cap \{v'_1,...,v'_{k'}\}.
\]

We therefore regard a fan $\Delta$ as a set of sets
of vectors $V$ satisfying the following conditions:
\begin{itemize}
\item
any subset of $V$ is in $\Delta$,
\item
if $V$ and $V'$ are in $\Delta$ then $V\cap V'$ is in $\Delta$.
\end{itemize} 

Moreover, $\Delta$ is complete if and only if 
$\qqq^n = \cup_{V\in \Delta}\cup_{v\in V}\qqq_{\geq 0}[v]$.
\end{remark}

\begin{example}
This is really a non-example.
Let $n=2$, $L=\zzz^2$.
Let
\[
\begin{array}{c}
\sigma_1=\qqq_{\geq 0}e_1+\qqq_{\geq 0}(e_1+e_2),\\
\sigma_2=\qqq_{\geq 0}(e_1+e_2)+\qqq_{\geq 0}e_2,\\
\sigma_3 =\qqq_{\geq 0}(e_1+e_2),\\
\sigma_4=\qqq_{\geq 0}e_1,\\
\sigma_5=\qqq_{\geq 0}e_2,\\
\sigma_6=\{0\},\\
\sigma_7= \qqq_{\geq 0}(-e_1)+\qqq_{\geq 0}e_2,\\
\sigma_8=\qqq e_1 +\qqq_{\geq 0}(-e_2).
\end{array}
\]

The set $\Delta=\{\sigma_i \ |\ 1\leq i\leq 8\}$
 may be visualized as follows:

\newpage
\begin{center}
\setlength{\unitlength}{.02cm}

\begin{picture}(256.00,150.00)(-150.00,0.00)
\thicklines
\put(0.00,0.00){\circle*{4}} 
\put(-23.00,12.00){$\sigma_6$} 

\put(0.00,0.00){\vector(0,1){100.00}}
\put(0.00,110.00){$\sigma_5$} 
 
\put(0.00,0.00){\vector(1,0){100.00}}
\put(110.00,0.00){$\sigma_4$} 
 
\put(0.00,0.00){\vector(1,1){100.00}}
\put(100.00,100.00){$\sigma_3$} 
 
\put(50.00,25.00){$\sigma_1$} 
\put(25.00,50.00){$\sigma_2$} 
 
\put(0.00,0.00){\vector(-1,0){100.00}}
\put(-60.00,60.00){$\sigma_7$} 
\put(0.00,-60.00){$\sigma_8$} 

\end{picture}

\end{center}

\vskip .5in

However, it is \underline{not} a fan since
$\sigma_8$ is not a strongly rational cone.
\end{example}

\subsection{Morphisms of fans and toric varieties}

This section is related to the material in \S \ref{sec:autoric}
below.

\index{morphism of fans}
\begin{definition}
\label{defn:autofan}
Let $\Delta$ be a fan of $L$ and $\Delta'$ a fan of $L'$. 
Given a homomorphism of lattices $\phi:L'\rightarrow L$ 
which maps each cone $\sigma'\in \Delta'$ to a cone
$\sigma\in \Delta$, we obtain a map
$\phi_*:\Delta'\rightarrow \Delta$ which we call
the associated {\bf morphism of fans}.
A morphism of a fan to itself which is a bijection as
a mapping of sets, and has an inverse which is a
morphism associated to the inverse of $\phi$,
will be called an {\bf automorphism of a fan}.
\end{definition}

\begin{remark}
In particular, the an automorphism of a fan
is associated to a automorphism of the underlying lattice.
MAGMA has a command to determine the automorphism
group of a lattice, \verb+AutomorphismGroup+.
\end{remark}

\index{equivariant maps between toric varieties}
\begin{definition}
Let $\Delta_1$ and $\Delta_2$ be fans in lattices 
$N_1\cong\Z^{r_1}$ and $N_1\cong\Z^{r_2}$.  Let $T_{N_i} = {N_i}\otimes k$ 
for $i=1,2$.  Suppose there is a map of algebraic tori:
\[
\psi: T_{N_1}\ra T_{N_2}.
\]
A map $f:X(\Delta_1)\ra X(\Delta_2)$ is called
{\bf equivariant} (with respect 
to $\psi$ and the actions of $T_{N_1}$ and $T_{N_2}$
if 
\[
f(\lambda\cdot x)=\psi(\lambda)\cdot f(x).
\]
\end{definition}

\begin{theorem}
Any map of fans $\phi:\Delta_1\ra\Delta_2$ gives rise to a map of
the corresponding toric varieties
\[
\phi_*:X(\Delta_1)\ra X(\Delta_2),
\]
which is equivariant with respect to the map given by
$\phi:N_1\ra N_2$.
Conversely, any equivariant homomorphism between toric varieties 
\[
f:X(\Delta_1)\ra X(\Delta_2)
\]
corresponds to a unique map $\phi:N_1\ra N_2$, giving rise to a map
between the fans $\Delta_1$ and $\Delta_2$, with
\[
f=\phi_* .
\]
\end{theorem}

By the above result, in order to give the equivariant
automorphism group of a toric variety,
we just need to consider the problem of finding the automorphisms of
the corresponding fan.

\begin{example}
\begin{itemize}
\item
The fan for $\P^2$ has automorphism group $D_3$, and the corresponding
automorphisms group of $\P^2$ is generated by
\[
(a:b:c)\mapsto (b:c:a)
\]
\[
(a:b:c)\mapsto (b:a:c)
\]
\item
Is is easy to construct two dimensional toric
varieties with equivariant
automorphism group the dihedral group $D_n$.
For example, $\P^1\times\P^1$ has automorphism group $D_4$, and
for the fan $\Delta$ in $\Z^2$ with cones
$\langle(0,\delta),(\epsilon,\delta)\rangle$ and 
$\langle(\delta,0),(\epsilon,\delta)\rangle$, where $\epsilon,\delta=\pm1$,
the automorphism group is isomorphic to $D_8$.  The variety 
$X(\Delta)$ is isomorphic to $\P^1\times\P^1$ blown up in $4$
points.
\[
\begin{picture}(0,0)%
\includegraphics{automorphism.pstex}%
\end{picture}%
\setlength{\unitlength}{1579sp}%
\begingroup\makeatletter\ifx\SetFigFont\undefined%
\gdef\SetFigFont#1#2#3#4#5{%
  \reset@font\fontsize{#1}{#2pt}%
  \fontfamily{#3}\fontseries{#4}\fontshape{#5}%
  \selectfont}%
\fi\endgroup%
\begin{picture}(2262,2466)(601,-3673)
\put(601,-1411){\makebox(0,0)[lb]{\smash{\SetFigFont{7}{8.4}{\familydefault}{\mddefault}{\updefault}{\color[rgb]{0,0,0}A fan with automorphism }%
}}}
\put(601,-1726){\makebox(0,0)[lb]{\smash{\SetFigFont{7}{8.4}{\familydefault}{\mddefault}{\updefault}{\color[rgb]{0,0,0}group $D_8$}%
}}}
\end{picture}
\]
\end{itemize}
\end{example}

\subsection{Constructing toric varieties from a fan}

We have the following recipe to construct a
quasi-projective toric variety.
\begin{center}
{\tiny{
\begin{tabular}{ccccccc}
$\{$ fans$\}$
&$\rightarrow $
&$\{$ ``fan'' of comm. semigps$\}$
&$\rightarrow$
&$\{$ ``fan'' of semigp algs$\}$
&$\rightarrow$ 
&$\{$  quasi-proj  schemes$\}$ \\
$\Delta$
& $\longmapsto$
&$\{S_\sigma = \sigma^*\cap L^*\ |\ \sigma\in\Delta\}$
& $\longmapsto$
&$\{\cccF [S_\sigma]\ |\ \sigma\in\Delta\}$
&$ \longmapsto$
&$X(\Delta)=\{U_\sigma +\ {\rm gluing\ maps}\ |\ \sigma\in\Delta\}$.\\
\end{tabular}
}}
\end{center}
\vskip .2in

Let $\Delta$ be a fan.
If $\tau$ is a face of both $\sigma_1,\sigma_2\in \Delta$
then Lemma \ref{lemma:subcone}
gives us maps
$\phi_1:U_\tau \rightarrow U_{\sigma_1}$ and
$\phi_2:U_\tau \rightarrow U_{\sigma_2}$. We ``glue''
the affine ``patch'' $U_{\sigma_1}$ to
the affine ``patch'' $U_{\sigma_{12}}$
the overlap using $\phi_1\circ \phi_2^{-1}$.
Define $X(\Delta)$ to be the ``glued scheme''
associated to these maps (as in Iitaka, \S 1.12,
\cite{I}):
\[
X(\Delta)=(\coprod_{\sigma\in\Delta}U_\sigma)/
({\rm gluing}).
\]
This is the {\bf toric variety associated to $\Delta$}.
It is known which toric varieties are complete (see
the lemma below) but, in general, the
classification of projective toric varieties is still
an active area of research
\footnote{See Oda \cite{O}, \S 2.4}.
\index{glue}
\index{toric variety, $X(\Delta)$}.

\begin{lemma}
(Fulton \S 2.4, \cite{F})
 $X(\Delta)$ is complete (as a variety) if and only if 
$\Delta$ is complete (as a fan).
\label{lemma:complete}
\end{lemma}

\index{equivariantly projective toric varieties}
\begin{definition}
Suppose that $\Delta$ is a fan,
$X=X(\Delta)$ is the associated toric variety,
and that there is an 
injective morphism 
$\phi:X\rightarrow \ppp^r$
such that the torus $T=U_{\{0\}}$ in $X$
maps to a subgroup of a torus $T'$ which is dense in
$\ppp^r$ in such a way that the action of $T'$ on
$\ppp^r$ extends that of $\phi(T)$ on $\phi(X(\Delta))$.
In this case, we say that $X$ is {\bf equivariantly
projective}.
\end{definition}

This basically means that there is a projective embedding which is
equivariant with respect to the torus action.

\begin{lemma}
If $\Delta$ is a 2-dimensional complete fan
and if $X(\Delta)$ is smooth then
it is in fact projective. In fact, 
$X(\Delta)$ is equivariantly projective.
\end{lemma}

\pf 
These statements follow from Ewald \cite{E}, Theorem
4.7 (page 160) and Theorem 3.11 (page 277).
\qed

\begin{example}
\label{example:P1}
Let $n=1$. Consider the fan
\[
\Delta =\{\{0\},\ \qqq_{\geq 0},\ \qqq_{\leq 0}
 \}.
\]
We've seen already in Example
\ref{ex:Rn} that
\[
R_\sigma \cong
\left\{
\begin{array}{cc}
\cccF [x],& {\rm if}\ \sigma = \qqq_{\geq 0},\\
\cccF [x^{-1}],& {\rm if}\ \sigma = \qqq_{\leq 0},\\
\cccF [x,x^{-1}],& {\rm if}\ \sigma = \{ 0\},
\end{array}
\right.
\]
and 
\[
U_\sigma \cong
\left\{
\begin{array}{cc}
\cccF,& {\rm if}\ \sigma = \qqq_{\geq 0},\\
\cccF^\times ,& {\rm if}\ \sigma = \{ 0\}.
\end{array}
\right.
\]
Similarly, one can show that if
$\sigma = \qqq_{\leq 0}$ then $U_\sigma \cong \cccF$.
We therefore have dominant maps
\[
U_{\qqq_{\leq 0}} \leftarrow U_{\{ 0\}} \rightarrow U_{\qqq_{\geq 0}},
\]
via the ``obvious'' embeddings
\[
\cccF [x^{-1}]\rightarrow \cccF [x,x^{-1}]  \leftarrow \cccF [x].
\]
The image of $\cccF[x]$ in $\cccF [x,x^{-1}]$ is isomorphic to the 
image of $\cccF[x^{-1}]$ in $\cccF [x,x^{-1}]$ via
the map $x\longmapsto x^{-1}$. This ``gluing'' defines
the variety $X(\Delta)=\ppp^1(\cccF)$. 

This is an example of a projective toric variety.

\end{example}

In general, a toric variety associated to a fan is
not projective. (Recall from Lemma \ref{lemma:complete} that
the fan is complete if and only if $X(\Delta)$ is.)

\begin{example}
\label{ex:3.7}
Let $n=2$, $L=\zzz^2$, $L^*=\zzz^2$.
Let
\[
\begin{array}{c}
\sigma_1=\qqq_{\geq 0}e_1+\qqq_{\geq 0}(e_1+e_2),\\
\sigma_2=\qqq_{\geq 0}(e_1+e_2)+\qqq_{\geq 0}e_2,\\
\sigma_3 =\qqq_{\geq 0}(e_1+e_2),\\
\sigma_4=\qqq_{\geq 0}e_1,\\
\sigma_5=\qqq_{\geq 0}e_2,\\
\sigma_6=\{0\},\\
\sigma_7=\qqq_{\geq 0}(- e_1) +\qqq_{\geq 0}(e_2),\\
\sigma_8= \qqq_{\geq 0}e_1+\qqq_{\geq 0}(-e_2),\\
\sigma_9=\qqq_{\geq 0}(-e_1),\\
\sigma_{10}=\qqq_{\geq 0}(- e_1) +\qqq_{\geq 0}(-e_2),\\
\sigma_{11}=\qqq_{\geq 0}(-e_2).
\end{array}
\]

The set $\Delta=\{\sigma_i \ |\ 1\leq i\leq 11\}$
 may be visualized as follows:

\newpage

\begin{center}
\setlength{\unitlength}{.02cm}

\begin{picture}(256.00,150.00)(-150.00,0.00)
\thicklines
\put(0.00,0.00){\circle*{4}} 
\put(-23.00,-15.00){$\sigma_6$} 

\put(0.00,0.00){\vector(0,1){100.00}}
\put(0.00,110.00){$\sigma_5$} 
 
\put(0.00,0.00){\vector(1,0){100.00}}
\put(110.00,0.00){$\sigma_4$} 
 
\put(0.00,0.00){\vector(1,1){100.00}}
\put(100.00,100.00){$\sigma_3$} 
 
\put(50.00,25.00){$\sigma_1$} 
\put(25.00,50.00){$\sigma_2$} 
 
\put(0.00,0.00){\vector(-1,0){100.00}}
\put(-40.00,50.00){$\sigma_7$} 
\put(50.00,-60.00){$\sigma_8$} 
 
\put(-50.00,-60.00){$\sigma_{10}$} 
 
\put(0.00,0.00){\vector(-1,0){100.00}}
\put(-115.00,5.00){$\sigma_9$} 

\put(0.00,0.00){\vector(0,-1){100.00}}
\put(0.00,-110.00){$\sigma_{11}$} 

\end{picture}
\end{center}

\vskip 1in

The cones $\sigma_7$, $\sigma_8$, and $\sigma_{10}$ are
self-dual. The cone $\sigma_1$ and its dual $\sigma_1^*$ 
(which contains $\sigma_1$) is pictured below.

\begin{center}
\setlength{\unitlength}{.02cm}

\begin{picture}(256.00,150.00)(-150.00,0.00)
\thicklines
\put(0.00,0.00){\circle*{4}} 

\put(0.00,0.00){\line(0,1){100.00}}
\put(0.00,0.00){\line(1,-1){100.00}}

\thinlines
\put(0.00,0.00){\vector(1,0){100.00}}
\put(0.00,0.00){\vector(1,1){70.00}}

\put(70.00,-40.00){$\sigma_1^*$}
\put(30.00,70.00){$\sigma_1^*$}
\put(110.00,70.00){$\sigma_1^*$}
\put(50.00,25.00){$\sigma_1$}

\end{picture}

\end{center}

\vskip 1in

In fact, we have
\[
\begin{array}{c}
\sigma_1^*=\{ae_1^*+be_2^*\ |\ a\geq 0,\ a+b\geq 0\}
=\qqq_{\geq 0}[e_2^*,e_1^*-e_2^*],\\
\sigma_2^*=\{ae_1^*+be_2^*\ |\ b\geq 0,\ a+b\geq 0\}
=\qqq_{\geq 0}[e_1^*,-e_1^*+e_2^*],\\
\sigma_3^*=\{ae_1^*+be_2^*\ |\  a+b\geq 0\}
=\qqq_{\geq 0}[e_1^*-e_2^*,e_1^*+e_2^*],\\
\sigma_4^*=\{ae_1^*+be_2^*\ |\  a\geq 0\}
=\qqq_{\geq 0}[e_1^*,e_2^*,-e_2^*],\\
\sigma_5^*=\{ae_1^*+be_2^*\ |\  b\geq 0\}
=\qqq_{\geq 0}[e_1^*,-e_1^*,e_2^*],\\
\sigma_6^*=\qqq e_1^* + \qqq e_2^*,\\
\sigma_7^*=\qqq_{\geq 0}[-e_1^*, e_2^*],\\
\sigma_8^*=\qqq_{\geq 0}[e_1^*,-e_2^*],\\
\sigma_9^*=\{ae_1^*+be_2^*\ |\  a\leq 0\}
=\qqq_{\geq 0}[-e_1^*,e_2^*,-e_2^*],\\
\sigma_{10}^*=\qqq_{\geq 0}[-e_1^*,-e_2^*],\\
\sigma_{11}^*=\{ae_1^*+be_2^*\ |\  b\leq 0\}
=\qqq_{\geq 0}[e_1^*,-e_1^*,-e_2^*].
\end{array}
\]
The associated semigroups are
\[
\begin{array}{c}
S_{\sigma_1} =\sigma_1^*\cap L^* 
= \zzz_{\geq 0} [e_1^*,e_2^*,e_1^*- e_2^*],\\
S_{\sigma_2} =\sigma_2^*\cap L^* 
= \zzz_{\geq 0} [e_1^*,e_2^*,-e_1^*+ e_2^*],\\
S_{\sigma_3} =\sigma_3^*\cap L^* 
= \zzz_{\geq 0} [e_1^*,e_2^*,e_1^*- e_2^*,-e_1^*+ e_2^*],\\
S_{\sigma_4} =\sigma_4^*\cap L^* 
= \zzz_{\geq 0} [e_1^*,e_2^*,- e_2^*],\\
S_{\sigma_5} =\sigma_5^*\cap L^* 
= \zzz_{\geq 0} [e_1^*,-e_1^*, e_2^*],\\
S_{\sigma_6} =\sigma_6^*\cap L^* 
= \zzz_{\geq 0} [e_1^*,-e_1^*, e_2^*,-e_2^*],\\
S_{\sigma_7} =\sigma_7^*\cap L^* 
= \zzz_{\geq 0} [-e_1^*,e_2^*],\\
S_{\sigma_8} =\sigma_8^*\cap L^* 
= \zzz_{\geq 0} [e_1^*,- e_2^*],\\
S_{\sigma_9} =\sigma_9^*\cap L^* 
= \zzz_{\geq 0} [-e_1^*,e_2^*,-e_2^*],\\
S_{\sigma_{10}} =\sigma_{10}^*\cap L^* 
= \zzz_{\geq 0} [-e_1^*,-e_2^*],\\
S_{\sigma_{11}} =\sigma_{11}^*\cap L^* 
= \zzz_{\geq 0} [e_1^*,-e_1^*,-e_2^*].
\end{array}
\]
(The element $e_1^*$ in the set of generators
of $S_{\sigma_1}$ and $S_{\sigma_3} $ are extraneous. 
However, their inclusion makes
some comparisons mentioned later a little bit easier.)
The associated coordinate rings are

\[
\begin{array}{c}
R_{\sigma_1} \cong \cccF[x_1,x_2,x_1x_2^{-1}]
\cong \cccF[x,y,w]/(wy-x)\\
\cong \cccF[x_2,x_1x_2^{-1}]
\cong \cccF[y,w],\\
R_{\sigma_2} \cong \cccF[x_1,x_2,x_1^{-1}x_2]
\cong \cccF[x,y,z]/(zx-y),\\
R_{\sigma_3} \cong \cccF[x_1,x_2,x_1^{-1}x_2,x_1x_2^{-1}]
\cong \cccF[x,y,z,w]/(xz-y,zw-1)\\
\cong \cccF[x_2,x_1^{-1}x_2,x_1x_2^{-1}]
\cong \cccF[y,z,w]/(zw-1),\\
R_{\sigma_4} \cong \cccF[x_1,x_2,x_2^{-1}]
\cong \cccF[x,y,v]/(vy-1),\\
R_{\sigma_5} \cong \cccF[x_1,x_1^{-1},x_2]
\cong \cccF[x,u,y]/(ux-1),\\
R_{\sigma_6} \cong \cccF[x_1,x_2,x_1^{-1},x_2^{-1}]
\cong \cccF[x,y,u,v]/(ux-1,vy-1),\\
R_{\sigma_7} \cong \cccF[x_1^{-1},x_2]\cong \cccF[u,y],\\
R_{\sigma_8} \cong \cccF[x_1,x_2^{-1}]\cong \cccF[x,v],\\
R_{\sigma_{9}} \cong \cccF[x_1^{-1},x_2,x_2^{-1}]\cong \cccF[u,y,v],\\
R_{\sigma_{10}} \cong \cccF[x_1^{-1},x_2^{-1}]\cong \cccF[u,v],\\
R_{\sigma_{11}} \cong \cccF[x_1,x_1^{-1},x_2^{-1}]\cong \cccF[x,u,v].
\end{array}
\]
From these, one can read off the equations for
$U_{\sigma_i}$, $1\leq i\leq 11$. For example, $U_{\sigma_6}$ is 
$ux-1,vy-1$, the torus $\cccF^\times \times \cccF^\times$.
We have
\newpage

\[
\begin{array}{c}
U_{\sigma_1} \cong \{(x,y,w)\ |\ wy-x=0\} ,\\
U_{\sigma_2} \cong \{(x,y,z)\ |\ zx-y=0\},\\
U_{\sigma_3} \cong \{(y,z,w)\ |\ zw=1\},\\
U_{\sigma_4} \cong \cccF\times \cccF^\times,\\
U_{\sigma_5} \cong \cccF\times \cccF^\times,\\
U_{\sigma_6} \cong \cccF^\times\times \cccF^\times,\\
U_{\sigma_7} \cong \cccF\times \cccF,\\
U_{\sigma_8} \cong \cccF\times \cccF,\\
U_{\sigma_{9}} \cong \cccF\times \cccF^\times,\\
U_{\sigma_{10}} \cong \cccF\times \cccF,\\
U_{\sigma_{11}} \cong \cccF\times \cccF^\times.
\end{array}
\]
In particular, all these affine ``patches''
are non-singular.

We have $\sigma_3\subset \sigma_1$ and
$\sigma_3\subset \sigma_2$. By Lemma
\ref{lemma:subcone}, there should be 
corresponding morphisms between the
associated affine toric varieties.
(These morphisms are used to ``glue'' these affine
pieces together.)
We have $U_{\sigma_3}$ is 
$zw-1=0$, which is the equation needed to
compare $U_{\sigma_2}$, which is $xz-y=0$, with
the equation for
$U_{\sigma_1}$, which is $wy-x=0$, on the intersection of the two.
We have maps
\[
U_{\sigma_1}
\leftarrow
U_{\sigma_3}
\rightarrow
U_{\sigma_2}.
\]
The map $z\longmapsto w^{-1}$ is the ``gluing''
map on the overlap.

As another example, we have $\sigma_6\subset \sigma_i$, 
$1\leq i \leq 11$. We have birational maps
\[
U_{\sigma_6}
\rightarrow
U_{\sigma_i},\ \ \ 
1\leq i \leq 11.
\]
This gives the dense embedding of the torus
in $X(\Delta)$.

It is known that any non-singular complete surface
is projective (see for example,
Theorem 8.9 in Iitaka \cite{I}), so
\[
X(\Delta)=(\coprod_{\sigma\in\Delta}U_\sigma)/
({\rm gluing}),
\]
is projective (it is complete since the fan $\Delta$ is
a complete fan).

It is in fact, equivariantly projective
since the fan $\Delta$ is strongly polytopal
\footnote{In the sense of Ewald \cite{E}, 
Definition 43 in \S V4, page 159, and
Theorem 3.11, \S VII, page 277. In two dimensions, this basically means
that the points in each wall of the polytope generate
(as a cone) the cones in the fan and the 
vertices  of the polytope generate
(as a cone) the separating rays between the
cones of the fan.}
 with supporting polytope
pictured below.

\begin{center}
\setlength{\unitlength}{.01cm}

\begin{picture}(256.00,110.00)(-150.00,0.00)
\thinlines
\put(0.00,0.00){\circle*{4}} 

\put(0.00,0.00){\vector(0,1){100.00}}
 
\put(0.00,0.00){\vector(1,0){100.00}}

\put(0.00,0.00){\vector(1,1){100.00}}
 
 
\put(0.00,0.00){\vector(-1,0){100.00}}
 
 
\put(0.00,0.00){\vector(-1,0){100.00}}

\put(0.00,0.00){\vector(0,-1){100.00}}

\thicklines

\put(50.00,50.00){\line(0,-1){50.00}}
\put(50.00,50.00){\line(-1,0){50.00}}

\put(0.00,-50.00){\line(1,1){50.00}}
\put(0.00,-50.00){\line(-1,1){50.00}}

\put(-50.00,0.00){\line(1,1){50.00}}

\end{picture}
\end{center}

\vskip .5in

\end{example}

Let $\phi:L'\rightarrow L$ be a homomorphism of lattices
and let $\phi_\qqq$ denote its extension to 
$V'\rightarrow V$ by linearity. Let
$\Delta$ be a fan of $L$ and $\Delta'$ a fan of $L'$ such that
for each cone $\sigma'\in \Delta'$ there is
a cone $\sigma\in \Delta$ such that 
$\phi_\qqq (\sigma')\subset \sigma$. Recall
from Definition \ref{defn:autofan} that
such a $\phi_\qqq$
is called a morphism of fans.

The next result recalls a generalization of 
Lemma \ref{lemma:isom}

\begin{lemma}
(Fulton \cite{F}, page 23 and Ewald \cite{E}, Theorem 6.1, 
page 243)
\label{lemma:functor}
Let $\phi_\qqq$ be a mapping of fans as above. Then 
the morphisms 
\[
\phi^*:{\rm Spec}(\cccF [S_\sigma])\rightarrow 
{\rm Spec}(\cccF[S_{\sigma'}]),
\]
of Lemma \ref{lemma:isom}
piece together to give a morphism 
$\phi^* : X(\Delta)\rightarrow  X(\Delta')$.
\end{lemma}

A morphism $\phi^*:X(\Delta)\rightarrow  X(\Delta')$
arising as in the above lemma from a homomorphism
$\phi:L'\rightarrow L$ be a homomorphism of lattices
will be called a {\bf toric morphism}.

Let $T$ be a torus which embeds as a dense subvariety 
of $X=X(\Delta)$. Identify $T$ with its image in $X$.
It turns out that a toric morphism $\psi:X\rightarrow X'$ is
$T$-equivariant
\footnote{In the sense that $\psi(t)\psi(x)=\psi(\mu(t)x)$,
where $\mu:T\rightarrow {\rm Aut}(X)$ is the natural action
of $T$ on $X$. See Ewald \cite{E}, Theorem 6.4, for a proof.}.

\subsection{The complement of the torus---stars}

As mentioned earlier, all toric varieties over $\cccF$ of dimension $n$ 
contain $(\cccF^\times )^n$ as a dense open subset, 
so it is the the complement 
$X(\Delta)\setminus (\cccF^\times )^n$ that distinguishes them (and how it
``glues'' to $(\cccF^\times )^n$).

\index{orbit}
\begin{definition}
The {\bf orbit} of a point $p$ in an $n$ dimensional
 toric variety $X(\Delta)$ 
is the set of points $a\cdot x$ where $a\in(\cccF^\times )^n$.
\end{definition}

A toric variety decomposes into a disjoint 
set of orbits.  One of the orbits is 
$(\cccF^\times )^n$ as a subset of $X(\Delta)$, and the other orbits form the
complement.  For example, for $\P^2$ we have the decomposition
\[
\P^2
=
\left\{
\renewcommand\arraystretch{1.5}
\begin{array}{l}
\{(a:b:c) | abc\not=0\}\cong(\cccF^\times )^2\\
\cup \{(0:b:c) | bc\not=0\}
\cup \{(a:0:c) | ac\not=0\}
\cup\{(a:b:0) | ab\not=0\}\\
\cup\{(1:0:0)\}
\cup\{(0:1:0)\}
\cup\{(0:0:1)\}
\end{array}
\right..
\]
The orbits are in one to one correspondence with the cones.  For each
$d$ dimensional cone there is a corresponding $n-d$ dimensional orbit.
If $F$ is algebraically closed, 
the (Zariski) closures of the orbits are themselves toric varieties
(see Lemma \ref{prop:3.22} below).

\index{star}
\begin{definition}
For a $d$ dimensional cone $\sigma$ in a fan $\Delta$
in $N\cong\Z^n$, the {\bf star} of $\sigma$ is
defined as follows.
Take any projection $\pi$
along $\sigma$ to a lattice $N_\sigma$ of dimension
$n-d$. In other words, there is a linear map
\[
\pi:N\ra N_\sigma,
\]
with $\ker\pi=\Z\sigma$.  This map extends to a map from $N_\Q$ to
$(N_\sigma)_\Q$ which I also denote by $\pi$.
Then define
\footnote{Note, though the concept of fan is standard, the notation 
$\Delta_\sigma$ used here may not be.
}
the fan $\Delta_\sigma$ in $N_\sigma$ to be the set
of images of all cones $\tau$ with $\sigma\prec\tau$ under the projection.
In other words, we have
\[
\Delta_\sigma  =  \{\pi(\tau) | \tau \in\Delta, \> \sigma\prec\tau\}.
\]
\end{definition}

\begin{magma}
To compute the star of a cone $\sigma$ in a fan $\Delta$,
use the \verb+star+ command in the {\bf toric} package.

Consider the fan $\Delta$ determined by the cones
\[
\{
\qqq_{\geq 0}[(2,-1),(1,0)],
\qqq_{\geq 0}[(1,0),(1,1)],
\qqq_{\geq 0}[(1,1),(2,0)]\}
\]
and let $\sigma=\qqq_{\geq 0}[(2,-1),(1,0)]$.
The star of $\sigma$ is the set $\{ \sigma\}$, since it is a 
cone of maximal dimension in $\Delta$. This, as well as some other
examples, are computed below using GAP.

{\footnotesize{
\begin{verbatim}
gap> RequirePackage("guava");
true
gap> Read("c:/gap/gapfiles/toric.g"); 
 toric.g for GAP 4.3, wdj, version 7-11-2002
available functions: dual_semigp_gens, cart_prod_lists,
 in_dual_cone, max_vectors, inner_product, toric_points,
 ideal_affine_toric_variety, embedding_affine_toric_variety, 
 toric_code, toric_codewords, divisor_polytope,
divisor_polytope_lattice_points, riemann_roch, flatten, 
faces, in_cone, normal_to_hyperplane, number_of_cones_dim , 
subcones_of_fan, star, betti_number, cardinality_of_X, 
euler_characteristic, ... 
gap> 
gap> Cones2:=[[[2,-1],[1,0]],[[1,0],[1,1]],[[1,1],[2,0]]];
[ [ [ 2, -1 ], [ 1, 0 ] ], [ [ 1, 0 ], [ 1, 1 ] ], [ [ 1, 1 ], [ 2, 0 ] ] ]
gap> star([[1,0]],Cones2);
[ [ [ 1, 0 ] ], [ [ 2, -1 ], [ 1, 0 ] ], [ [ 1, 0 ], [ 1, 1 ] ] ]
gap> star([[1,0],[2,-1]],Cones2);
[ [ [ 2, -1 ], [ 1, 0 ] ] ]
gap> 
gap> Cones3:=[ [ [2,0,0],[0,2,0],[0,0,2] ], [ [2,0,0],[0,2,0],[1,1,-2] ] ];
[ [ [ 2, 0, 0 ], [ 0, 2, 0 ], [ 0, 0, 2 ] ], 
  [ [ 2, 0, 0 ], [ 0, 2, 0 ], [ 1, 1, -2 ] ] ]
gap> star([[2,0,0]],Cones3);
[ [ [ 2, 0, 0 ] ], [ [ 0, 0, 2 ], [ 2, 0, 0 ] ], [ [ 0, 2, 0 ], [ 2, 0, 0 ] ],
  [ [ 1, 1, -2 ], [ 2, 0, 0 ] ], [ [ 2, 0, 0 ], [ 0, 2, 0 ], [ 0, 0, 2 ] ], 
  [ [ 2, 0, 0 ], [ 0, 2, 0 ], [ 1, 1, -2 ] ] ]
gap> star([[2,0,0],[0,2,0]],Cones3);
[ [ [ 0, 2, 0 ], [ 2, 0, 0 ] ], [ [ 2, 0, 0 ], [ 0, 2, 0 ], [ 0, 0, 2 ] ], 
  [ [ 2, 0, 0 ], [ 0, 2, 0 ], [ 1, 1, -2 ] ] ]
gap> 
\end{verbatim}
}}
\end{magma}

The following results may be found in \S 3.1 in \cite{F}.

\begin{proposition}
\label{prop:3.22}
Assume $F$ is algebraically closed.
The closures of the orbits of points in $X(\Delta)$ are in correspondence
with the cones of $\Delta$, and are given by $X(\Delta_\sigma)$ for
$\sigma\in\Delta$.
\end{proposition}

\begin{corollary}
Assume $F$ is algebraically closed.
The closures of the components of the complement of $(\cccF^\times )^n$ in a 
complete toric
variety $X(\Delta)$ are isomorphic to the toric varieties $X(\Delta_\sigma)$
where $\sigma$ runs through all the one dimensional cones in $\Delta$.
\end{corollary}

\begin{example}
Consider the fan $\Delta$ in $\Z^3$ with 
one dimensional cones are generated by $(\pm 1,\pm1,\pm1)$.
These points can be pictured as the vertices of a cube with center at the 
origin.  For each face of the cube we define a
three-dimensional cone generated by the
vectors at the vertices of the face.  So there are $6$ three-dimensional
faces, and $X(\Delta)$ is complete, and has $6$ singularities.
The star of each one dimensional cone contains $3$ two-dimensional cones,
and gives a toric variety isomorphic to $\P^2$. 

\vskip .1in

\begin{picture}(0,0)%
\includegraphics{star.pstex}%
\end{picture}%
\setlength{\unitlength}{1184sp}%
\begingroup\makeatletter\ifx\SetFigFont\undefined%
\gdef\SetFigFont#1#2#3#4#5{%
  \reset@font\fontsize{#1}{#2pt}%
  \fontfamily{#3}\fontseries{#4}\fontshape{#5}%
  \selectfont}%
\fi\endgroup%
\begin{picture}(14366,6903)(-302,-5248)
\put(-299,1439){\makebox(0,0)[lb]{\smash{\SetFigFont{5}{6.0}{\familydefault}{\mddefault}{\updefault}{\color[rgb]{0,0,0}The $8$ one dimensional cones of $\Sigma$}%
}}}
\put(11116,854){\makebox(0,0)[lb]{\smash{\SetFigFont{5}{6.0}{\familydefault}{\mddefault}{\updefault}{\color[rgb]{0,0,0}This picture represents the union of}%
}}}
\put(11116,539){\makebox(0,0)[lb]{\smash{\SetFigFont{5}{6.0}{\familydefault}{\mddefault}{\updefault}{\color[rgb]{0,0,0}$8$ copies of $P^2$ which form}%
}}}
\put(11116,224){\makebox(0,0)[lb]{\smash{\SetFigFont{5}{6.0}{\familydefault}{\mddefault}{\updefault}{\color[rgb]{0,0,0}$X(\Sigma)\setminus (k^*)^3$.}%
}}}
\put(5401,914){\makebox(0,0)[lb]{\smash{\SetFigFont{5}{6.0}{\familydefault}{\mddefault}{\updefault}{\color[rgb]{0,0,0}Intersection of the fan with}%
}}}
\put(5401,599){\makebox(0,0)[lb]{\smash{\SetFigFont{5}{6.0}{\familydefault}{\mddefault}{\updefault}{\color[rgb]{0,0,0}a cube and a hyperplane:}%
}}}
\put(7951,-3061){\makebox(0,0)[lb]{\smash{\SetFigFont{5}{6.0}{\familydefault}{\mddefault}{\updefault}{\color[rgb]{0,0,0}The star of a one dimensional}%
}}}
\put(7951,-3376){\makebox(0,0)[lb]{\smash{\SetFigFont{5}{6.0}{\familydefault}{\mddefault}{\updefault}{\color[rgb]{0,0,0}cone:}%
}}}
\put(10276,-4861){\makebox(0,0)[lb]{\smash{\SetFigFont{5}{6.0}{\familydefault}{\mddefault}{\updefault}{\color[rgb]{0,0,0}This gives the toric variety $\P^2$}%
}}}
\end{picture}

The picture shows how these copies of $\P^2$ intersect.  The intersections
correspond to the $2$ dimensional cones of $\Delta$.

\end{example}

\section{Polyhedra and support functions}

In this section we describe a useful way of constructing 
toric varieties using polyhedra.

\subsection{Toric varieties from polyhedra}

An alternative way to define a toric variety is, instead of starting with
a fan in a lattice $N$, to start with a polyhedron in the dual lattice $M$.
This method has limitations, in that it can only produce complete
projective toric varieties, but nevertheless, it is worth considering.
In particular 
there are some nice results about certain subvarieties of toric varieties,
which come down to counting points in faces of the corresponding polyhedra.

\index{fan $\Delta(P)$ for a polyhedron $P$}
\begin{definition}
Let $N$ be a lattice of rank $n$, and $P$ a convex polyhedron in $N_\Q$
with vertices at points in $N$.  

Then if $F$ is a face of $P$ of dimension $d$, we define a cone 
$\widehat\sigma(F)$ in $L^*$ by 
\[
\widehat\sigma(F)=
\left\{\lambda(x-y)\ |\ x\in P, y\in F,\lambda \in \Q
\right\}
\]
Then we define a cone $\sigma(F)$ in $N_\Q$ by taking the dual:
\[
\sigma(F)=(\widehat\sigma(F))^* .
\]
\index{$\Delta(P)$}
The fan $\Delta(P)$ is given by
\[
\Delta(P) = \{\sigma(F) : F\ {\rm a\ face\ of }\ P\}.
\]
\end{definition}

\[
\begin{picture}(0,0)%
\includegraphics{poly1.pstex}%
\end{picture}%
\setlength{\unitlength}{1579sp}%
\begingroup\makeatletter\ifx\SetFigFont\undefined%
\gdef\SetFigFont#1#2#3#4#5{%
  \reset@font\fontsize{#1}{#2pt}%
  \fontfamily{#3}\fontseries{#4}\fontshape{#5}%
  \selectfont}%
\fi\endgroup%
\begin{picture}(10672,7982)(-1649,-6702)
\put(3376,389){\makebox(0,0)[lb]{\smash{\SetFigFont{7}{8.4}{\familydefault}{\mddefault}{\updefault}{\color[rgb]{0,0,0}$p_1$}%
}}}
\put(4576,-511){\makebox(0,0)[lb]{\smash{\SetFigFont{7}{8.4}{\familydefault}{\mddefault}{\updefault}{\color[rgb]{0,0,0}$p_2$}%
}}}
\put(3151,-1261){\makebox(0,0)[lb]{\smash{\SetFigFont{7}{8.4}{\familydefault}{\mddefault}{\updefault}{\color[rgb]{0,0,0}$p_3$}%
}}}
\put(4051, 14){\makebox(0,0)[lb]{\smash{\SetFigFont{7}{8.4}{\familydefault}{\mddefault}{\updefault}{\color[rgb]{0,0,0}$e_3$}%
}}}
\put(3901,-886){\makebox(0,0)[lb]{\smash{\SetFigFont{7}{8.4}{\familydefault}{\mddefault}{\updefault}{\color[rgb]{0,0,0}$e_1$}%
}}}
\put( 76,-2161){\makebox(0,0)[lb]{\smash{\SetFigFont{7}{8.4}{\familydefault}{\mddefault}{\updefault}{\color[rgb]{0,0,0}Points in $P-e_1$}%
}}}
\put( 76,-2476){\makebox(0,0)[lb]{\smash{\SetFigFont{7}{8.4}{\familydefault}{\mddefault}{\updefault}{\color[rgb]{0,0,0}circled:}%
}}}
\put(-1274,839){\makebox(0,0)[lb]{\smash{\SetFigFont{7}{8.4}{\familydefault}{\mddefault}{\updefault}{\color[rgb]{0,0,0}Polyhedron $P$  in a lattice:}%
}}}
\put(-1649,-1711){\makebox(0,0)[lb]{\smash{\SetFigFont{7}{8.4}{\familydefault}{\mddefault}{\updefault}{\color[rgb]{0,0,0}Example of construction of $\sigma(F)$:}%
}}}
\put(2776,-436){\makebox(0,0)[lb]{\smash{\SetFigFont{7}{8.4}{\familydefault}{\mddefault}{\updefault}{\color[rgb]{0,0,0}$e_2$}%
}}}
\put(3676,-436){\makebox(0,0)[lb]{\smash{\SetFigFont{7}{8.4}{\familydefault}{\mddefault}{\updefault}{\color[rgb]{0,0,0}$f$}%
}}}
\put(3151,839){\makebox(0,0)[lb]{\smash{\SetFigFont{7}{8.4}{\familydefault}{\mddefault}{\updefault}{\color[rgb]{0,0,0}Labels for the faces:}%
}}}
\put(4276,-4186){\makebox(0,0)[lb]{\smash{\SetFigFont{5}{6.0}{\familydefault}{\mddefault}{\updefault}{\color[rgb]{0,0,0}$0$}%
}}}
\put(6676,-2536){\makebox(0,0)[lb]{\smash{\SetFigFont{7}{8.4}{\familydefault}{\mddefault}{\updefault}{\color[rgb]{0,0,0}${\sigma(e_1)}:$}%
}}}
\put(6451,-4861){\makebox(0,0)[lb]{\smash{\SetFigFont{7}{8.4}{\familydefault}{\mddefault}{\updefault}{\color[rgb]{0,0,0}${\sigma(p_3)}:$}%
}}}
\put(901,-6511){\makebox(0,0)[lb]{\smash{\SetFigFont{5}{6.0}{\familydefault}{\mddefault}{\updefault}{\color[rgb]{0,0,0}$0$}%
}}}
\put(3901,-6511){\makebox(0,0)[lb]{\smash{\SetFigFont{5}{6.0}{\familydefault}{\mddefault}{\updefault}{\color[rgb]{0,0,0}$0$}%
}}}
\put(5926,1064){\makebox(0,0)[lb]{\smash{\SetFigFont{7}{8.4}{\familydefault}{\mddefault}{\updefault}{\color[rgb]{0,0,0}The fan $\Delta(P)$}%
}}}
\put(3601,-2536){\makebox(0,0)[lb]{\smash{\SetFigFont{7}{8.4}{\familydefault}{\mddefault}{\updefault}{\color[rgb]{0,0,0}$\widehat{\sigma}(e_1):$}%
}}}
\put(3226,-4861){\makebox(0,0)[lb]{\smash{\SetFigFont{7}{8.4}{\familydefault}{\mddefault}{\updefault}{\color[rgb]{0,0,0}$\widehat{\sigma}(p_3):$}%
}}}
\put(1501,-3886){\makebox(0,0)[lb]{\smash{\SetFigFont{5}{6.0}{\familydefault}{\mddefault}{\updefault}{\color[rgb]{0,0,0}$0$}%
}}}
\put(151,-4861){\makebox(0,0)[lb]{\smash{\SetFigFont{7}{8.4}{\familydefault}{\mddefault}{\updefault}{\color[rgb]{0,0,0}Points in $P - p_3$:}%
}}}
\end{picture}
\]

Conversely, for a fan we can find a corresponding polyhedron.  However,
this is not unique, as shown in the following examples:
\[
\begin{picture}(0,0)%
\includegraphics{poly2.pstex}%
\end{picture}%
\setlength{\unitlength}{1579sp}%
\begingroup\makeatletter\ifx\SetFigFont\undefined%
\gdef\SetFigFont#1#2#3#4#5{%
  \reset@font\fontsize{#1}{#2pt}%
  \fontfamily{#3}\fontseries{#4}\fontshape{#5}%
  \selectfont}%
\fi\endgroup%
\begin{picture}(8442,8997)(1,-9442)
\put(  1,-661){\makebox(0,0)[lb]{\smash{\SetFigFont{7}{8.4}{\familydefault}{\mddefault}{\updefault}{\color[rgb]{0,0,0}A fan $\Delta $ with}%
}}}
\put(  1,-976){\makebox(0,0)[lb]{\smash{\SetFigFont{7}{8.4}{\familydefault}{\mddefault}{\updefault}{\color[rgb]{0,0,0}four cones:}%
}}}
\put(3901,-661){\makebox(0,0)[lb]{\smash{\SetFigFont{7}{8.4}{\familydefault}{\mddefault}{\updefault}{\color[rgb]{0,0,0}The cones each define an}%
}}}
\put(3901,-976){\makebox(0,0)[lb]{\smash{\SetFigFont{7}{8.4}{\familydefault}{\mddefault}{\updefault}{\color[rgb]{0,0,0}``angle'':}%
}}}
\put(1801,-9361){\makebox(0,0)[lb]{\smash{\SetFigFont{7}{8.4}{\familydefault}{\mddefault}{\updefault}{\color[rgb]{0,0,0}Three polyhedra $\Delta$ with $\Delta (P)=\Delta $}%
}}}
\put(601,-5161){\makebox(0,0)[lb]{\smash{\SetFigFont{7}{8.4}{\familydefault}{\mddefault}{\updefault}{\color[rgb]{0,0,0}The ``angles'' glue together to give polyhedra:}%
}}}
\end{picture}
\]
A polyhedra gives more structure, and can be used to define a ``polarization''
of the toric variety.  This corresponds to taking a support function on a
fan, which we will see later.  A fan with support function gives
rise to a unique polyhedra. 

\begin{remark}
 We do not have to define a fan to get a
toric variety from a polyhedron in a lattice.  We can go directly from
a polyhedron to a toric variety as follows:

\index{$X(P)$}

\begin{definition}
For a polyhedron $P$ in a lattice $L^*$, define a cone $C_P$  on
$P$ by
\[
C_P=\{0\}\cup 
\{(u,\boldsymbol{v})\in \Q\oplus L^*\ |\ u > 0, (\boldsymbol{v}/u)\in P\}.
\]
Then let 
\[
S_P =\cccF [C_P]
\]
be the $\cccF$-algebra generated by $C_P$.  
As before, we use variables $y_i$ when we 
consider elements of $L^*$ multiplicatively, writing
\[
e^{a_1\boldsymbol{e}_1+\dots+a_n\boldsymbol{e}_n}=
y_1^{a_1}y_2^{a_2}\dots y_n^{a_n}.
\]
The algebra $S_P$ has a grading 
inherited from the grading on $k[y_1,y_1^{-1},\dots y_n,y_n^{-1}]$.
We define the toric variety associated to $P$ by
\[
X(P)={\rm Proj}(\cccF [S_P]),
\]
where ${\rm Proj}$ of a ring means the collection 
of homogenous prime ideals in the ring.
\end{definition}
\end{remark}
\index{$C_P$}
\index{$S_P$}

To go from a fan to a polyhedron, we need to define support functions.

\subsection{Support functions}

This section is based on Chapter 2 of Oda's book \cite{O}.

\index{support functions}
\begin{definition}
Let $\Delta$ be a fan in a lattice $L$.
A {\bf $\Delta$-linear support function} is a real valued function
\[
h:\vert\Delta\vert\ra\R,
\]
such that 
\[
h:\vert\Delta\vert\cap N\ra\Z,
\]
and $h$ is linear on all cones $\sigma\in\Delta$.
\end{definition}

This means that if $\sigma,\tau\in \Delta$, then 
$h|_\sigma$ and $h|_\tau$
are linear functions, which are equal on $\sigma\cap\tau$.
To illustrate a support function on a fan $\Delta$ we 
can draw lines where the fan has constant value $0$, $1$, $2$,
and so on, as in the examples in the following pictures:
\[
\begin{picture}(0,0)%
\includegraphics{support.pstex}%
\end{picture}%
\setlength{\unitlength}{1579sp}%
\begingroup\makeatletter\ifx\SetFigFont\undefined%
\gdef\SetFigFont#1#2#3#4#5{%
  \reset@font\fontsize{#1}{#2pt}%
  \fontfamily{#3}\fontseries{#4}\fontshape{#5}%
  \selectfont}%
\fi\endgroup%
\begin{picture}(10362,3528)(-2399,-3673)
\put(-2399,-361){\makebox(0,0)[lb]{\smash{\SetFigFont{7}{8.4}{\familydefault}{\mddefault}{\updefault}{\color[rgb]{0,0,0}Lines of constant integer values for $3$ different support functions for the same fan}%
}}}
\end{picture}
\]
In the above examples the value of the support function is always positive,
but this is not necessary.  Any element of the dual, $L^*$ gives a support
function on any fan $\Delta$ in $L$.

A support function can be defined by specifying
for each cone $\sigma$ in $\Delta$ an element $h_\sigma\in L^*$
such that
\begin{equation}
\label{eqn:hn}
h(n)=\langle n,h_\sigma\rangle
\end{equation}
for all $n\in\sigma$, and such that 
$\langle n,h_\sigma\rangle=\langle n,h_\tau\rangle$ for $n\in\tau\prec\sigma$.
\index{$h_\sigma$}

\index{$SF(L,\Delta)$}
\index{$\Delta$-linear support functions}
\begin{definition}
\label{def:SF}
For a fan $\Delta$, the set $SF(L,\Delta)$ is the set of all 
$\Delta$-linear support functions.
\end{definition}

A support function is determined by its values on the one dimensional cones
of a fan.
If we let
\[
\Delta(1) = \{\sigma\in\Delta | \dim\sigma =1\},
\]
then we have
\begin{equation}
\label{eqn:4.2}
SF(\Delta)\cong\Z^{\Delta(1)}.
\end{equation}
(See \S 2.1 in \cite{O}.)

\subsection{From a fan with support function to a polyhedron}

\index{$P(\Delta,h)$}
\begin{definition}
Suppose $\Delta$ is a complete fan in $L^*\cong\Z^n$, 
with support function $h$.  Let $\Delta(n)$ be the set of
$n$ dimensional cones in $\Delta$.
\index{$\Delta(n)$}
Then we define a polyhedron $P(\Delta,h)$ by
\[
P(\Delta,h)=\{m\in L^*_\Q\ |\ \langle m,n\rangle \ge -h(n),\ 
\forall n\in L_\Q\}
=\bigcap_{\sigma\in\Delta(n)}(-h_\sigma + \sigma^\vee),
\]
where $h_\sigma$ satisfies (\ref{eqn:hn}).
\end{definition}

To see the second equality above, 
note that if $m = -h_\sigma + s$ for some $s\in\sigma^*$,
we have for $n\in\sigma$ that 
\[
\langle m,n\rangle = \langle -h_\sigma + s,n\rangle
=
\langle -h_\sigma,n \rangle + \langle s,n\rangle
\ge \langle h_\sigma,n \rangle=h(n),
\]
So
$\bigcap_{\sigma\in\Delta(n)}( -h_\sigma + \sigma^\vee)\subset
P(\Delta,h)$.
Conversely, if $\langle m,n\rangle \ge -h(n)$ for all $n\in\sigma\in\Delta$,
then for all $n\in\sigma$ we have
\[
\langle m+h_\sigma,n\rangle=\langle m,n\rangle + h(n)\ge 0
\]
so $m+h_\sigma\in\sigma^\vee$, i.e., $m\in -h_\sigma+\sigma^{\vee}$,
which gives the reverse inclusion.


\begin{example}
In this diagram we show how to construct the polyhedron $P (\Delta,h)$
for a fan with support.
\[
\begin{picture}(0,0)%
\includegraphics{poly3.pstex}%
\end{picture}%
\setlength{\unitlength}{1579sp}%
\begingroup\makeatletter\ifx\SetFigFont\undefined%
\gdef\SetFigFont#1#2#3#4#5{%
  \reset@font\fontsize{#1}{#2pt}%
  \fontfamily{#3}\fontseries{#4}\fontshape{#5}%
  \selectfont}%
\fi\endgroup%
\begin{picture}(11236,9528)(-1949,-9073)
\put(3616,-2356){\makebox(0,0)[lb]{\smash{\SetFigFont{5}{6.0}{\familydefault}{\mddefault}{\updefault}{\color[rgb]{0,0,0}$0$}%
}}}
\put(6751,-2626){\makebox(0,0)[lb]{\smash{\SetFigFont{7}{8.4}{\familydefault}{\mddefault}{\updefault}{\color[rgb]{0,0,0}part of the edge of}%
}}}
\put(6751,-2941){\makebox(0,0)[lb]{\smash{\SetFigFont{7}{8.4}{\familydefault}{\mddefault}{\updefault}{\color[rgb]{0,0,0}the polyhedron}%
}}}
\put(6751,-2311){\makebox(0,0)[lb]{\smash{\SetFigFont{7}{8.4}{\familydefault}{\mddefault}{\updefault}{\color[rgb]{0,0,0}this outline becomes}%
}}}
\put(3601,-5011){\makebox(0,0)[lb]{\smash{\SetFigFont{7}{8.4}{\familydefault}{\mddefault}{\updefault}{\color[rgb]{0,0,0}are marked round the sides}%
}}}
\put(-1649,239){\makebox(0,0)[lb]{\smash{\SetFigFont{7}{8.4}{\familydefault}{\mddefault}{\updefault}{\color[rgb]{0,0,0}$L$ with level sets of a support}%
}}}
\put(2101,-526){\makebox(0,0)[lb]{\smash{\SetFigFont{7}{8.4}{\familydefault}{\mddefault}{\updefault}{\color[rgb]{0,0,0}and $h_\sigma$ circled:}%
}}}
\put(2101,-211){\makebox(0,0)[lb]{\smash{\SetFigFont{7}{8.4}{\familydefault}{\mddefault}{\updefault}{\color[rgb]{0,0,0}$L^*$, with $\sigma^*$ shaded}%
}}}
\put(-1649,-61){\makebox(0,0)[lb]{\smash{\SetFigFont{7}{8.4}{\familydefault}{\mddefault}{\updefault}{\color[rgb]{0,0,0}function $h$ marked}%
}}}
\put(-1649,-361){\makebox(0,0)[lb]{\smash{\SetFigFont{7}{8.4}{\familydefault}{\mddefault}{\updefault}{\color[rgb]{0,0,0}and $\sigma$ is shaded}%
}}}
\put(-1949,-1411){\makebox(0,0)[lb]{\smash{\SetFigFont{7}{8.4}{\familydefault}{\mddefault}{\updefault}{\color[rgb]{0,0,0}$\sigma$}%
}}}
\put(1501,-1411){\makebox(0,0)[lb]{\smash{\SetFigFont{7}{8.4}{\familydefault}{\mddefault}{\updefault}{\color[rgb]{0,0,0}$\sigma^*$}%
}}}
\put(-1649,-4261){\makebox(0,0)[lb]{\smash{\SetFigFont{7}{8.4}{\familydefault}{\mddefault}{\updefault}{\color[rgb]{0,0,0}In this picture $L^*$ and $L$ are}%
}}}
\put(-1649,-4576){\makebox(0,0)[lb]{\smash{\SetFigFont{7}{8.4}{\familydefault}{\mddefault}{\updefault}{\color[rgb]{0,0,0}identified via the standard inner }%
}}}
\put(-1649,-4891){\makebox(0,0)[lb]{\smash{\SetFigFont{7}{8.4}{\familydefault}{\mddefault}{\updefault}{\color[rgb]{0,0,0}product on $\Z^2$, and $h_\sigma$ }%
}}}
\put(-1649,-5206){\makebox(0,0)[lb]{\smash{\SetFigFont{7}{8.4}{\familydefault}{\mddefault}{\updefault}{\color[rgb]{0,0,0}is marked with a point in $\sigma$ for}%
}}}
\put(-1649,-5521){\makebox(0,0)[lb]{\smash{\SetFigFont{7}{8.4}{\familydefault}{\mddefault}{\updefault}{\color[rgb]{0,0,0}all $\sigma\in\Delta$}%
}}}
\put(3601,-4111){\makebox(0,0)[lb]{\smash{\SetFigFont{7}{8.4}{\familydefault}{\mddefault}{\updefault}{\color[rgb]{0,0,0}In this figure the points $-h_\sigma$ are}%
}}}
\put(3601,-4411){\makebox(0,0)[lb]{\smash{\SetFigFont{7}{8.4}{\familydefault}{\mddefault}{\updefault}{\color[rgb]{0,0,0}marked, and the polygon $P (\Delta,h)$}%
}}}
\put(4111,-811){\makebox(0,0)[lb]{\smash{\SetFigFont{7}{8.4}{\familydefault}{\mddefault}{\updefault}{\color[rgb]{0,0,0}outline of $\sigma^* - h_\sigma$}%
}}}
\put(3601,-4711){\makebox(0,0)[lb]{\smash{\SetFigFont{7}{8.4}{\familydefault}{\mddefault}{\updefault}{\color[rgb]{0,0,0}is drawn. Dual ``angles'' of the cones}%
}}}
\end{picture}
\]
Here are $3$ more examples:
\[
\begin{picture}(0,0)%
\includegraphics{poly4.pstex}%
\end{picture}%
\setlength{\unitlength}{1579sp}%
\begingroup\makeatletter\ifx\SetFigFont\undefined%
\gdef\SetFigFont#1#2#3#4#5{%
  \reset@font\fontsize{#1}{#2pt}%
  \fontfamily{#3}\fontseries{#4}\fontshape{#5}%
  \selectfont}%
\fi\endgroup%
\begin{picture}(10512,6832)(-2549,-6077)
\put(571,-6001){\makebox(0,0)[lb]{\smash{\SetFigFont{7}{8.4}{\familydefault}{\mddefault}{\updefault}{\color[rgb]{0,0,0}The corresponding polygons}%
}}}
\put(-2399, 89){\makebox(0,0)[lb]{\smash{\SetFigFont{7}{8.4}{\familydefault}{\mddefault}{\updefault}{\color[rgb]{0,0,0}Angles of dual cones indicated, and points $h_\sigma$ marked ($L$ and $L^*$ are identified).}%
}}}
\put(-2549,539){\makebox(0,0)[lb]{\smash{\SetFigFont{7}{8.4}{\familydefault}{\mddefault}{\updefault}{\color[rgb]{0,0,0}Lines of constant integer values for $3$ different support functions for the same fan}%
}}}
\end{picture}
\]
In the above four examples we see that in only one case of $h$
do we get a polyhedron which will satisfy
$\Delta(P (\Delta,h)) = \Delta$.  In a theorem below we will give
the condition we need on $h$ to have this property.
\end{example}

\subsection{From polyhedra to fans with support function}

\index{$h_P$}

We saw how to get from a polyhedron to a fan.  
Here we try to reverse that direction.

\begin{definition}
\index{support function, upper convex}
\index{convex, upper}
\index{convex, support function}
\index{convex, strictly upper}
\index{strictly convex}
A support function $h$ on $N\cong\Z^n$ is {\bf upper convex} if 
\[
h(n_1+n_2)\le h(n_1)+h(n_2).
\]
If in addition 
\[
h_{\sigma_1}\not=h_{\sigma_2},
\]
for $\sigma_1,\sigma_2$
two $n$ dimensional cones in $\Delta$, then $h$ is called 
{\bf strictly upper convex}.
\end{definition}

Fulton \cite{F} uses the same terminology but omits the 
adjective ``upper''.

A polyhedron $P $
in $L^*$ also gives a strictly upper convex 
support function on the fan $\Delta(P )$ in $L$, 
which is defined on $L$ by:
\[
h_P (n)=-\inf\{\langle m,n\rangle\ |\ m\in P \}.
\]

\begin{theorem}
For a support function $h$ on a fan $\Delta$, we have 
\[
\Delta(P (\Delta,h)) = \Delta
\]
if and only if $h$ is strictly upper convex.
\end{theorem}

For further details on the above result, see Lemma 2.12 and
Theorem 2.22 in \cite{O}.

We also have the following important result:

\begin{theorem}
For a (finite) complete fan $\Delta$, the toric variety $X(\Delta)$
is projective (i.e., can be embedded by rational functions 
in $\P^m$ for some $m$) if and only if there is a strictly upper convex 
support function on $\Delta$.
\end{theorem}

For further details, see Theorem 2.13 and
Corollary 2.16 in \cite{O}.

Note that a polyhedron $P $ in $L^*$ defines a finite complete fan with a 
strictly upper convex 
support function, $h_P $, so $X(\Delta(P ))$ is projective.
If we want to restrict attention to projective toric varieties, we can
restrict ourselves to toric varieties given by polyhedra.
To understand where the above result comes from, we need to look at linear
systems (or line bundles or sheaves) 
on $X(\Delta)$ defined by support functions.

\section{Resolution of singularities for surfaces}
\index{resolution of singularities for surfaces}

In this section, we describe in some detail how to resolve
singularities of an affine toric variety.

\subsection{Singular examples}

First, we shall give some examples of singular toric
surfaces. 

\subsubsection{Example 1}

Consider the cone $\sigma$ generated by
$e_1$ and $3e_1+2e_2$. This is a ``cone in $L$'', where
$L=\zzz e_1+\zzz e_2=\zzz^2$.
The angle between the generating vectors is 
less than $\pi/2$, so we shall call this an ``acute cone.''

\begin{center}
\setlength{\unitlength}{.02cm}

\begin{picture}(256.00,100.00)(-100.00,0.00)
\thicklines
\put(0.00,0.00){\circle*{4}} 
\put(0.00,0.00){\vector(1,0){50.00}}
\put(60.00,0.00){$e_1$} 
 
\put(0.00,0.00){\vector(3,2){100.00}}
\put(110.00,70.00){$3e_1+2e_2$} 
 
\put(80.00,30.00){$\sigma$} 

\end{picture}
\end{center}

\vskip .8in

We have 
$\sigma=\qqq_{\geq 0}[e_1,3e_1+2e_2]$,
$\sigma^*=\qqq_{\geq 0}[2e_1^*-3e_2^*,e_2^*]$, so
\[
S_\sigma=\zzz_{\geq 0}[e_2^*,e_1^*-e_2^*,2e_1^*-3e_2^*].
\]
This implies that
$\cccF[S_\sigma] = \cccF[x_2,x_1x_2^{-1},x_1^2x_2^{-3}]
\cong \cccF[x,y,z]/(y^2-xz)$,
so $U_\sigma$ is the (singular) cone
$y^2=xz$.

To resolve this, we try to choose a subcone
(more precisely, a refinement of the fan
associated to $\sigma$) $\tau\subset \sigma$
and we try to find a morphism
$\pi:U_\tau\rightarrow U_\sigma$
such that $U_\tau$ is smooth and
$\pi$ is birational.

Let $\tau=\qqq_{\geq 0} [e_1,2e_1+e_2]$,
so $\tau^*=\qqq_{\geq 0} [e_1^*-2e_2^*,e_2^*]$, 
and
\[
S_{\tau}=\tau^*\cap L^*
=\zzz_{\geq 0}[   e_1^*-2e_2^*,e_2^*].
\]
This implies $\cccF[S_\tau]=\cccF[x_2,x_1x_2^{-2}]$,
so $U_\tau \cong \cccF^2$. Since
$U_\sigma = \cccF[x_2,x_1x_2^{-1},x_1^2x_2^{-3}]
\cong \cccF[x,y,z]/(y^2-xz)$, $(X,Y)\longmapsto (X,XY,XY^2)$
defines a morphism $U_\tau\rightarrow U_\sigma$.
This is the resolution $\pi$.

\subsubsection{Example 2}
\label{sec:example2}

Consider the cone $\sigma$ generated by
$e_2$ and $3e_1-2e_2$. This is a ``cone in $L$'', where
$L=\zzz e_1+\zzz e_2=\zzz^2$.

The angle between the generating vectors is 
greater than $\pi/2$, so we shall call this an ``obtuse cone.''

\begin{center}
\setlength{\unitlength}{.02cm}

\begin{picture}(256.00,120.00)(-100.00,0.00)
\thicklines
\put(0.00,0.00){\circle*{4}} 
\put(0.00,0.00){\vector(0,1){50.00}}
\put(0.00,60.00){$e_2$} 
 
\put(0.00,0.00){\vector(3,-2){100.00}}
\put(120.00,-80.00){$3e_1-2e_2$} 
 
\put(80.00,60.00){$\sigma$} 

\end{picture}

\end{center}

\vskip .8in

The dual cone $\sigma^*$ is given by
\[
\sigma^*=\{ae^*_1+be^*_2\ | \ b\geq 0, \ \frac{3}{2}a\geq b\},
\]
so it is generated by $e^*_1$ and $2e^*_1+3e^*_2$. 

\vskip .5in

\begin{center}
\setlength{\unitlength}{.02cm}

\begin{picture}(256.00,120.00)(-100.00,0.00)
\thicklines
\put(0.00,0.00){\circle*{4}} 
\put(0.00,0.00){\vector(1,0){50.00}}
\put(60.00,0.00){$e_1^*$} 
 
\put(0.00,0.00){\vector(2,3){100.00}}
\put(105.00,135.00){$2e_1^*+3e_2^*$} 
 
\put(80.00,70.00){$\sigma^*$} 

\end{picture}

\end{center}
\vskip .4in

The
associated semigroup $S_\sigma$ is given by
\[
S_\sigma = \sigma^*\cap L^*,
\]
and is generated (as a semigroup) by
\[
u_1=e^*_1,\ \ u_2=2e^*_1+3e^*_2,\ \  u_3=e^*_1+e^*_2.
\]
What are the equations defining
the affine variety $U_\sigma= {\rm Spec}\, \cccF [S_\sigma]$
associated to the cone $\sigma$?

To determine these equations, we use Lemma 
\ref{lemma:ucondition} above.

For our example, let $c_i=a_i-b_i$, $i=1,2,3$,
in the notation of Lemma 
\ref{lemma:ucondition}.
In this case, (\ref{eqn:ucondition}) becomes
\[
\vec{0}=c_1u_1+c_2u_2+c_3u_3=
c_1e_1^*+c_2(2e_1^*+3e_2^*)+c_3(e_1^*+e_2^*)
=(c_1+2c_2+c_3,3c_2+c_3).
\]
This implies $3c_2+c_3=0,\ c_1+2c_2+c_3=0$,
so (subtracting) $c_1=c_2$. In other words,
\[
a_1-b_1=a_2-b_2,\ \ 3(a_2-b_2)=b_3-a_3.
\]
In particular, if $a_1=b_1$ then we must also have
$a_2=b_2$ and $a_3=b_3$, so
$y_1^{a_1}y_2^{a_2}y_3^{a_3}-y_1^{b_1}y_2^{b_2}y_3^{b_3}=0$.
Moreover, if $a_1>b_1$ then $a_2>b_2$ and
$a_3<b_3$.
For example, if we take
\[
a_1=1,\ b_1=0,\ 
a_2=1,\ b_2=0,\ 
a_3=0,\ b_3=3,
\]
then $y_1y_2-y_3^3$ is a generator of $I$.
Since the value of $c_1$ determines those of
$c_2$ and $c_3$, we have
\[
I=(y_1y_2-y_3^3),
\]
so
\[
U_\sigma = 
\{ (y_1,y_2,y_3)\ |\ y_1y_2=y_3^3 \}.
\]
Moreover, in the notation of Fulton \cite{F}, \S 1.3, we have
\[
\chi^{u_1}=x_1,\ 
\chi^{u_2}=x_1^2x_3^3,\ 
\chi^{u_3}=x_1x_2,
\]
so the coordinate ring of $U_\sigma$ is given by
\[
\cccF [S_\sigma] =
\cccF [x_1,x_1^2x_3^3,x_1x_2],
\]
by Lemma \ref{lemma:ucondition}.

From this equation for $U_\sigma$, we see that
it has a singularity at $(0,0,0)$ (the Jacobian
of $y_1y_2-y_3^3$ vanishes at the origin).
How can the surface $U_\sigma$ be desingularized?

To this end, consider the two subcones of $\sigma$,
\[
\sigma_1= \qqq_{\geq 0} e_1+\qqq_{\geq 0} e_2,
\]
\[\sigma_2= \qqq_{\geq 0} (3e_1-2e_2)+\qqq_{\geq 0} e_1.
\]


\begin{center}
\setlength{\unitlength}{.02cm}

\begin{picture}(256.00,100.00)(-100.00,0.00)
\thicklines
\put(0.00,0.00){\circle*{4}} 
\put(0.00,0.00){\vector(1,0){50.00}}
\put(60.00,0.00){$e_1$} 
\put(0.00,0.00){\vector(0,1){50.00}}
\put(0.00,60.00){$e_2$} 
\put(60.00,60.00){$\sigma_1$} 

\thinlines
\put(0.00,0.00){\vector(3,-2){100.00}}
\put(110.00,-70.00){$3e_1-2e_2$} 
 
\put(80.00,-30.00){$\sigma_2$} 

\end{picture}

\end{center}
\vskip .7in

\noindent
We shall see that the varieties associated to these subcones 
have birational morphisms satisfying
\[
U_{\sigma_1}\rightarrow U_\sigma,\ \ 
U_{\sigma_2}\rightarrow U_\sigma.
\]

First, we determine the equations defining
$U_{\sigma_1}, U_{\sigma_2}$ as affine varieties.
To do this, we compute the dual cones and then
use Lemma \ref{lemma:ucondition} as we did for
$U_\sigma$ above.
The dual cone for $\sigma_1$ is easy
\[
\sigma_1^*=\qqq_{\geq 0} e^*_1
+\qqq_{\geq 0} e^*_2.
\]
The dual cone for $\sigma_2$ is given by
\[
\begin{array}{c}
\sigma_2^*=
\{ a e^*_1+b e^*_2\ |\ 
(a e^*_1+b e^*_2)e_1\geq 0,\ 
(a e^*_1+b e^*_2)(3e_1-2e_2)\geq 0\}\\
=\{ a e^*_1+b e^*_2\ |\ 
a\geq 0,\ 
\frac{3}{2}a-b\geq 0\}\\
=\qqq_{\geq 0} (-e^*_1)+
\qqq_{\geq 0} (2e^*_1+3e^*_2).
\end{array}
\]

Next, we determine the semigroups $S_{\sigma_1}$ and
$S_{\sigma_2}$.
We have
\[
S_{\sigma_1}=\sigma^*_1\cap L^*
=\zzz_{\geq 0} e^*_1
+\zzz_{\geq 0} e^*_2,
\]
and
\[
S_{\sigma_2}=\sigma^*_2\cap L^*
=\zzz_{\geq 0} (-e^*_2)
+\zzz_{\geq 0} (2e^*_1+3e^*_2)
+\zzz_{\geq 0} (e^*_1+e^*_2).
\]
In the first case, 
\[
\cccF [S_{\sigma_1}]
\cong \cccF [x_1,x_2],
\]
which is the coordinate ring of our variety $U_{\sigma_1}$.
Thus
\[
U_{\sigma_1}\cong \cccF^2,
\]
which is of course non-singular.

In the second case, $S_{\sigma_2}$ is generated by
$u_1=-e^*_2$,
$u_2=2e^*_1+3e^*_2$,
$u_3=e^*_1+e^*_2$, so
\[
\cccF [S_{\sigma_2}]
\cong \cccF [x_2^{-1},x_1^2x_2^3,x_1x_2]
\cong \cccF [y_1,y_2,y_3]/I,
\]
where $I$ is generated by
\[
y_1^{a_1}y_2^{a_2}y_3^{a_3}-y_1^{b_1}y_2^{b_2}y_3^{b_3},
\]
and $a_i\geq 0$, $b_j\geq 0$ satisfy
\[
c_1(-e^*_2)+c_2(2e^*_1+3e^*_2)+c_3(e^*_1+e^*_2)=\vec{0},
\]
where $c_i=a_i-b_i$, $i=1,2,3$. This implies
$2c_2+c_3=0$, $-c_1+3c_2+c_3=0$, 
so (subtracting) we have
$c_1=c_2$ and $c_3=-2c_2$, or
$a_1-b_1=a_2-b_2$ and $a_3-b_3=-2(a_2-b_2$.
If 
\[
a_1=1,\ b_1=0,\ 
a_2=1,\ b_2=0,\ 
a_3=0,\ b_3=2,
\]
then $y_1y_2-y_3^2$ is a generator of $I$.
As in the case of $U_\sigma$, we have that 
$U_{\sigma_1}$ is given by $y_1y_2=y_3^2$.
This is singular.

What is the relationship, if any (and there is one),
between $U_{\sigma}$, $U_{\sigma_1}$, and 
$U_{\sigma_2}$?

Since $\sigma^*$ is a subcone of $\sigma^*_1$,
$S_\sigma$ is a sub-semigroup of $S_{\sigma_1}$,
so there is an inclusion homomorphism
\[
\cccF [S_{\sigma}]
\rightarrow 
\cccF [S_{\sigma_1}].
\]
Indeed, the obvious map
\[
\cccF [S_{\sigma}]
=\cccF [x_1,x_1^2x_2^3,x_1x_2]
\rightarrow
\cccF [S_{\sigma_1}]
=\cccF [x_1,x_2]
\]
fits the bill.
This implies that there is a morphism
\[
U_{\sigma_1} = {\rm Spec}(\cccF [S_{\sigma_1}])
\rightarrow
U_\sigma = {\rm Spec}(\cccF [S_\sigma]).
\]

The semigroups $S_{\sigma_1}$ and $S_{\sigma}$
generate the same sub{\it group} of $L^*$
(namely $L^*$ itself). Therefore (by
\cite{F}, page 18, Exercise), we see that the
morphism $U_{\sigma_1} \rightarrow
U_\sigma$ is birational. In fact, we can
determine this morphism explicitly.
To this end, recall
$U_\sigma$ is given by $y_1y_2=y_3^3$.
The map
\[
(y_1,y_3)\longmapsto (y_1,y_3^3/y_1,y_3),
\]
sends a point of
$U_{\sigma_1}-\{(y_1,y_3)\ |\ y_1\not= 0\}$
to a point of
$U_{\sigma}-\{(y_1,y_2,y_3)\in U_\sigma \ |\ y_1\not= 0\}$.
This is the birational morphism,
$U_{\sigma_1}\rightarrow U_{\sigma}$.

Since $\sigma^*$ is a subcone of $\sigma^*_2$,
$S_\sigma$ is a sub-semigroup of $S_{\sigma_2}$,
so there is an inclusion homomorphism
\[
\cccF [S_{\sigma}]
\rightarrow 
\cccF [S_{\sigma_2}].
\]
Indeed, $\cccF [S_{\sigma}]
=\cccF [x_1,x_1^2x_2^3,x_1x_2]$ is contained in
$\cccF [S_{\sigma_2}]= \cccF [x_2^{-1},x_1^2x_2^3,x_1x_2]$,
since $x_1=(x_2^{-1})(x_1x_2)$.
There is a morphism
\[
U_{\sigma_2} \rightarrow U_\sigma
\]
which is birational (for the same reason as in the
$U_{\sigma_1}$ case above).
In fact, we can
determine this morphism explicitly.
To this end, recall
$U_{\sigma_2}$ is given by $y_1y_2=y_3^2$.
Send a point
\[
(Y_1,Y_2,Y_3)\in U_{\sigma_2}
\]
to the point
\[
(y_1,y_2,y_3)=(Y_1,Y_2/Y_3,Y_3)\in U_{\sigma}.
\]
This defines a birational morphism (defined on the subset
$Y_3\not= 0$).

Next, consider the two subcones of $\sigma_2$,
\[
\sigma_{21}= \qqq_{\geq 0} e_1+\qqq_{\geq 0} (2e_1-e_2),
\]
\[
\sigma_{22}= \qqq_{\geq 0} (2e_1-e_2)+\qqq_{\geq 0} (3e_1-2e_2).
\]


\newpage

\begin{center}
\setlength{\unitlength}{.02cm}

\begin{picture}(256.00,80.00)(-80.00,0.00)
\thicklines
\put(0.00,0.00){\circle*{4}} 
\put(0.00,0.00){\vector(1,0){50.00}}
\put(60.00,0.00){$e_1$} 
\put(0.00,0.00){\vector(2,-1){100.00}}
\put(110.00,-50.00){$2e_1-e_2$} 
\put(80.00,-20.00){$\sigma_{21}$} 

\thinlines
\put(0.00,0.00){\vector(3,-2){140.00}}
\put(140.00,-90.00){$3e_1-2e_2$} 
 
\put(110.00,-70.00){$\sigma_{22}$} 

\end{picture}

\end{center}
\vskip 1in

\noindent
We shall see that the varieties associated to these subcones 
have birational morphisms satisfying
\[
U_{\sigma_{21}}\rightarrow U_{\sigma_2},\ \ 
U_{\sigma_{22}}\rightarrow U_{\sigma_2}.
\]

To compute the equations defining $U_{\sigma_{21}}$, $U_{\sigma_{22}}$, 
we first determine the dual cones $\sigma^*_{21}$, $\sigma^*_{21}$.
We have
\[
\sigma^*_{21}=
\{ a e^*_1+b e^*_2\ |\ a\geq 0,\ 2a\geq b\}
=\zzz_{\geq 0}(-e^*_2)+\zzz_{\geq 0}(e^*_1+2e^*_2)
\]
and
\[
\begin{array}{c}
\sigma^*_{22}=
\{ a e^*_1+b e^*_2\ |\ 3a\geq 2b,\ 2a\geq b\}\\
=\zzz_{\geq 0}(-e^*_1-2e^*_2)+\zzz_{\geq 0}(2e^*_1+3e^*_2).
\end{array}
\]
These give
\[
\cccF [S_{\sigma^*_{21}}]=\cccF[x_2^{-1},x_1x_2^2],
\]
which contains $x_1$, and
\[
\cccF [S_{\sigma^*_{22}}]=\cccF[x_1^{-1}x_2^{-2},x_1^2x_2^3].
\]

Since $\sigma_{21}$ is generated by $e_1$, $2e_1-e_2$ (which is a basis for
$L$), $U_{\sigma_{21}}$ is smooth, by the above lemma. In fact,
$U_{\sigma_{21}}\cong \cccF^2$.
(Indeed, the ideal $I$ associated to
$U_{\sigma_{21}}$ is generated by 
$y_1^{a_1}y_2^{a_2}y_3^{a_3}-y_1^{b_1}y_2^{b_2}y_3^{b_3}$
where $c_1(-e^*_2)+c_2(e^*_1+2e^*_2)=\vec{0},$
$c_i=a_i-b_i$, $i=1,2,3$. These imply $c_1=c_2=0$, so $I$ is trivial.)

The ideal $I$ for $U_{\sigma_{22}}$ is generated by the binomials
$y_1^{a_1}y_2^{a_2}y_3^{a_3}-y_1^{b_1}y_2^{b_2}y_3^{b_3}$
where $c_1(-e^*_1-2e^*_2)+c_2(2e^*_1+3e^*_2)=\vec{0},$
$c_i=a_i-b_i$, $i=1,2$. This implies 
$c_1=2c_2$, $2c_1=3c_2$. These imply $c_1=c_2=0$, 
so $I$ is trivial.

What is the relationship between
\[
U_{\sigma_{2}}, U_{\sigma_{21}}, U_{\sigma_{22}}?
\]
Recall $U_{\sigma_{2}}$ was given by $y_1y_2=y_3^3$. The morphism
\[
(y_1,y_3)\longmapsto
(y_1,y_3^3/y_1,y_3)
\]
defines a birational morphism
\[
U_{\sigma_{21}}
\rightarrow
U_{\sigma_{2}} ,
\]
and similarly for $U_{\sigma_{22}}
\rightarrow U_{\sigma_{2}}$.

The composite morphism
\[
U_{\sigma_{21}}
\rightarrow
U_{\sigma_{2}} \rightarrow
U_{\sigma} 
\]
is a resolution of singularities.
This corresponds to the refinement 
pictured below.



\begin{center}
\setlength{\unitlength}{.02cm}

\begin{picture}(256.00,80.00)(-100.00,0.00)
\thinlines
\put(0.00,0.00){\circle*{4}} 
\put(0.00,0.00){\vector(1,0){50.00}}
\put(60.00,0.00){$e_1$} 
\put(0.00,0.00){\vector(2,-1){100.00}}
\put(110.00,-50.00){$2e_1-e_2$} 
\put(80.00,-20.00){$\sigma_{21}$} 

\put(60.00,60.00){$\sigma_{1}$} 

\thicklines
\put(0.00,0.00){\vector(0,1){50.00}}
\put(0.00,60.00){$e_2$} 

\put(0.00,0.00){\vector(3,-2){140.00}}
\put(155.00,-100.00){$3e_1-2e_2$} 
 
\put(115.00,-70.00){$\sigma_{22}$} 

\end{picture}

\end{center}
\vskip .9in



\subsection{Resolution of singularities}

Let $\Delta$ be a fan of $L=\zzz^2$. A refinement
$\Delta'$ of $\Delta$ yields an equivariant
birational map
$X(\Delta')\rightarrow X(\Delta)$.
We seek here two things.
\begin{itemize}
\item
A refinement $\Delta'$ for which $X(\Delta')$ is
smooth.
\item
An explicit expression for the morphism
$X(\Delta')\rightarrow X(\Delta)$.
\end{itemize}
To simplify matters, we shall assume that 
$\Delta$ is the fan associated to a cone
of $L$. Furthermore, we shall
simply replace $X(\Delta)$ by $U_\sigma$.
The examples in the previous section
provide a good blueprint of what to expect.

\begin{magma}

Let $n=2$, $L=L^*=\zzz^2$, and suppose
$\sigma=\qqq_{\geq 0}[e_1,3e_1+4e_2]$. 
To obtain the birational desingularization (as schemes)
$\cccF^2 \rightarrow U_\sigma$ in MAGMA, type

{\footnotesize{
\begin{verbatim}
load "/home/wdj/magmafiles/toric.mag";
//replace /home/wdj/magmafiles by your path to toric.mag
D := LatticeDatabase();
Lat := Lattice(D, 2, 16);
desing_affine_toric_variety([[1,0],[3,4]],Lat);
\end{verbatim}
}}
\noindent
MAGMA returns
{\footnotesize{
\begin{verbatim}
Mapping from: Affine Space of dimension 2
Variables : x1, x2 to Affine Space of dimension 5
Variables : $.1, $.2, $.3, $.4, $.5
with equations : 
$.2
$.1
$.1^2/$.2
$.1^3/$.2^2
$.1^4/$.2^3 %$
\end{verbatim}
}}
\noindent
Basically, this is the same as the dense embedding map,
given in MAGMA \ref{magma:embed}, except that the 
notation for the variables returned by MAGMA is a
little different (a MAGMA ``feature'', 
since the mapping does not belong to the 
coordinate ring of functions).
This means that $U_\sigma$ is $2$-dimensional
and the mapping is given by
given by
\[
\begin{array}{ccc}
\cccF^2&\rightarrow &U_\sigma,\\
(x_1,x_2)&\longmapsto & (x_2,x_1,x_1^2x_4, x_1^3x_4^2, x_1^4x_4^3)
=(x_2,x_1,x_1^2x_2^{-1}, x_1^3x_2^{-2}, x_1^4x_2^{-3}).
\end{array}
\]
\end{magma}

Let $\sigma_{a,b}$ denote the following two-dimensional cone
\[
\sigma_{a,b}=
\qqq_{\geq 0}[e_1,ae_1+be_2],
\]
where $(a,b)\in \zzz^2$, $b\not= 0$. Let $S_{a,b}=
S_{\sigma_{a,b}}$ and $U_{a,b}= U_{\sigma_{a,b}}$.
By Oda \cite{O}, Proposition 1.24,
$U_{a,b}\cong \cccF^2/\mu_b$
(see also Exercise 2.10 above).

\begin{lemma}
\label{lemma:5.1}
(a)
$\sigma_{a,b}\cong \sigma_{a,-b}$ and
$\sigma_{a,b}\cong \sigma_{a\pm b,b}$.

(b) $U_{a,b}\cong U_{a,-b}$ and
$U_{a,b}\cong U_{a\pm b,b}$. 
\end{lemma}

\pf  
(a) The map 
$A=
\left(
\begin{array}{cc}
1 & 0\\
0 & -1 
\end{array}
\right)\in GL_2(\zzz)$
defines an isomorphism
$\sigma_{a,b}\rightarrow \sigma_{a,-b}$
and $A=
\left(
\begin{array}{cc}
1 & 0\\
\pm 1 & 1 
\end{array}
\right)\in GL_2(\zzz)$
defines an isomorphism
$\sigma_{a,b}\rightarrow \sigma_{a\pm b,b}$.

(b) This follows from Lemma 
\ref{lemma:isom} using (a).
\qed 

The following result shows that, in dimension 2, ``most'' affine 
toric varieties are
of the above form (up to isomorphism).

\begin{lemma}
Let $a,b,c,d\in \zzz$ with $gcd(a,b)=1$
(in particular, $a\not=0,\ b\not=0$),
and let
$\sigma=\qqq_{\geq 0}[(a,b),(c,d)]$
be a 2-dimensional
cone. There are $e,f\in \zzz$ such that $\sigma\cong \sigma_{e,f}$,
hence
$U_\sigma\cong U_{e,f}$.
\end{lemma}

\pf 
First, if $gcd(a,b)=r$ then we may replace $a$ by $a/r$ 
and $b$ by $b/r$, since
$\qqq_{\geq 0}[(a,b),(c,d)]=\qqq_{\geq 0}[(a/r,b/r),(c,d)]$
(they are equal, not just isomorphic).
So, we may assume without loss of generality
that $gcd(a,b)=1$. In this case
(by a standard result in elementary number theory), 
there exist $x,y\in \zzz$ such that
$ax+by=1$. We have
\[
\left(
\begin{array}{cc}
a & b\\
c & d 
\end{array}
\right)
\left(
\begin{array}{cc}
x & -b\\
y & a 
\end{array}
\right)
=
\left(
\begin{array}{cc}
1 & b'\\
c' & d' 
\end{array}
\right),
\]
for some $b',c',d'\in\zzz$. Note
$\left(
\begin{array}{cc}
x & -b\\
y & a 
\end{array}
\right)\in GL_2(\zzz)$.
Finally, note 
\[
\left(
\begin{array}{cc}
1 & b'\\
c' & d' 
\end{array}
\right)
\left(
\begin{array}{cc}
1 & b'\\
0 & -1 
\end{array}
\right)
=
\left(
\begin{array}{cc}
1 & 0\\
c' & d'' 
\end{array}
\right),
\]
for some $d''\in \zzz$.
Let $e=c', f=d''$.
\qed

Because of the following lemma, we may dispose of the 
case $b=\pm 1$. 

\begin{lemma}
If $\sigma=\sigma_{a,1}$ or $\sigma=\sigma_{a,-1}$ then
$U_\sigma\cong \cccF^2$. (In particular,
$U_\sigma$ is smooth.
\end{lemma}

\pf 
Exercise (using Lemma \ref{lemma:ucondition}). \qed 

Therefore, we may assume $|n|>1$.

Let $\sigma_0$ denote the cone generated by
$e_1,e_2$. Let $\sigma=\sigma_{m,n}$,
where $m\not= 0$ and $n\not= 0$.
In this case, $U_{\sigma_0}\cong \cccF^2$ is smooth and there
is a birational resolution 
$\pi:U_{\sigma_0}\rightarrow U_{\sigma}$. 

We explicitly (as is possible) determine $\pi$. 

It suffices to consider the case $n>0$ and $gcd(|m|,n)=1$.
(By Lemma \ref{lemma:5.1}
above, we may assume we are in this case.)
Note
$\sigma^*=\qqq_{\geq 0}[e_2^*,ne_1^*-me_2^*]$.

\begin{lemma}
$S_\sigma$ is generated by $e_2^*$, $ne_1^*-me_2^*$,
and elements of the form
$ae_1^*+be_2^*$, $(a,b)\in J$,
where $J$ is a subset of lattice points in
$L$ contained in the parallelogram 
from $(0,0)$ to
$(0,1)$ to $(n,1-m)$ to $(n,-m)$.
\end{lemma}

\newpage

\begin{center}
\setlength{\unitlength}{.01cm}

\begin{picture}(256.00,210.00)(-130.00,0.00)
\put(0.00,0.00){\circle{2}} 
\put(60.00,120.00){\circle{2}} 
\put(120.00,60.00){\circle{2}} 
\put(0.00,0.00){\circle*{6}} 
\put(60.00,-60.00){\circle{2}} 
\put(-120.00,-60.00){\circle{2}} 
\put(-60.00,-120.00){\circle{2}} 
\put(180.00,0.00){\circle{2}} 
\put(0.00,60.00){\circle{2}} 
\put(-60.00,60.00){\circle{2}} 
\put(-60.00,0.00){\circle*{6}} 
\put(-120.00,0.00){\circle{2}} 
\put(120.00,-60.00){\circle{2}} 
\put(0.00,-60.00){\circle{2}} 
\put(180.00,-60.00){\circle{2}} 
\put(180.00,0.00){\circle{2}} 
\put(180.00,60.00){\circle{2}} 
\put(180.00,120.00){\circle*{6}} 
\put(120.00,120.00){\circle{2}} 
\put(180.00,180.00){\circle{2}} 
\put(120.00,180.00){\circle{2}} 
\put(60.00,0.00){\circle{2}} 
\put(120.00,0.00){\circle{2}} 
\put(60.00,120.00){\circle{2}} 
\put(60.00,60.00){\circle{2}} 
\qbezier[250](-60.00,0.00)(180.0,120.0)(180.0,120.0)
\qbezier[250](180.00,120.00)(-60.0,0.0)(-60.0,0.0)
\thicklines
\put(-60.00,-60.00){\line(0,1){60.00}} 
\put(180.00,120.00){\line(0,1){60.00}} 
\put(-60.00,-60.00){\line(4,3){238.00}} 
\put(-60.00,0.00){\line(4,3){238.00}} 
\put(-450,120){generators for $S_{\sigma_{-3,4}}$}
\end{picture}

\end{center}
\vskip .5in

Before proving this lemma, note it implies 
\[
\cccF[S_\sigma]=\cccF[x_2,x_2^nx_2^{-m},\{x_1^ax_2^b\}_{(a,b)\in J}].
\]
If we identify $U_{\sigma_0}$ with
$\cccF^2$, hence $\cccF[S_{\sigma_0}]$
with $\cccF[x_1,x_2]$, then the
map
\begin{equation}
\label{eqn:resol}
(x_1,x_2)\longmapsto (x_2,x_2^nx_2^{-m},\{x_1^ax_2^b\}_{(a,b)\in J})
\end{equation}
defines a morphism $\pi:
U_{\sigma_0}\rightarrow U_{\sigma}$.
This is the explicit resolution.

\pf
The semigroup generators of $S_\sigma$ must belong to
$L\cap \sigma^*$. Each such generator must in fact 
belong to the parallelogram spanned by the vectors 
$(0,1)$ and $(n,-m)$ since the parallelogram forms a 
fundamental domain for $\qqq^2/\zzz[ e_2^*,ne_1^*-me_2^*]$.
\qed

Here's a stronger version of the above lemma.

\begin{theorem} 
$S_\sigma$ is generated by $e_2^*$, $ne_1^*-me_2^*$,
and elements of the form
$ae_1^*+be_2^*$, $(a,b)\in J$, where $J$ is 
contained in the triangle from $(0,0)$ to
$(0,1)$ to $(n,-m)$.
\end{theorem}

\pf
By the above lemma, we know that the generators of 
$S_\sigma$ belong to the parallelogram generated by
the vectors $e_2^*$, $ne_1^*-me_2^*$.
Let $P$ be this parallelogram, let $T$ be the 
triangle of the theorem, and let $w$ denote the
vertex $e_2^*+ne_1^*-me_2^*$.
To prove the theorem, we need to show 
that: if $v\in P\cap L$ then there are $t_i\in T\cap L$,
$i=1,2$, such that $v=t_1+t_2$. Symbolically, we need
to show: if $v\not= w$ then
\[
T\cap L\cap (v- T\cap L)\not= \emptyset .
\]
By symmetry, this is equivalent to the following
{\it claim}: if $v\in w-T\cap L$, $v\not= w$, then
\[
(w-T\cap L)\cap (w-v+ T\cap L)\not= \emptyset .
\]
This is geometrically ``obvious''. To prove it, suppose not.
Assume 
\[
(w-T\cap L)\cap (w-v+ T\cap L) = \emptyset .
\]
Let $v$ be a vector satisfying this and where 
$w-v$ is as small as possible
(lexicographically). But 
\[
(w-T\cap L)\cap (v_0+ T\cap L)\not= \emptyset,
\]
where $v_0=(1,1)$. Since $v_0$ is the smallest vector in 
$T\cap L$, our assumption must be false.
\qed

\section{Riemann-Roch spaces}

In this section we will assume $X$ is non singular, so that Weil diviors
and Cartier divisors are the same.  (The space of
Cartier divisors is the subspace of locally principal Weil divisors.)

\subsection{Divisors and Linear Systems}

\index{Weil divisors, group of}
\begin{definition}
For an algebraic variety $X$ the group of {\bf Weil divisors} on $X$ is
(roughly) given by
\[
{\rm Div}(X) = \Z[\{{\rm irreducible\ subvarieties\ of }\ X
\ {\rm of\ codimension }\ 1\}].
\]
\end{definition}

\index{${\rm orb}(\sigma)$}
Previously we have mentioned orbits of points in toric varieties.
For a cone $\sigma\in\Delta$ we define
\[
{\rm orb}(\sigma) = \{ {\rm unique}\ T-{\rm orbit\ in\ }U_\sigma
\ {\rm which\ is\ closed\ in}\ U_\sigma 
\},
\]
where as usual closed is with respect to the Zariski topology.

For example, for $\P^2$, with homogenous coordinates $X,Y,Z$, one
of the affine pieces corresponding to a cone is $X\not=0$.  In this piece
the orbit $XYZ\not=0$ is not closed, since it is not defined by the
vanishing of any finite set of polynomials.  But the subvariety given
by the point $(1:0:0)$ is closed and $T$ invariant.
In the affine piece $XY\not=0$, the unique $T$ invariant divisor is
given by the line $Z=0$ (restricted to this piece).

Define a Weil divisor corresponding to $\sigma$ by:
\[
V(\sigma)={\rm closure\ of\ orb}(\sigma).
\]
\index{$V(\sigma)$}
\index{divisor, Weil}

\index{$T$-invariant\\ divisors}
\begin{definition}
For a toric variety $X$ with dense open torus $T$, a Weil divisor
$D$ is {\bf $T$ invariant} if $D=T\cdot D$.
The space of $T$ invarant Weil divisors is denoted $T{\rm Div}(X)$.
\end{definition}
\index{$T{\rm Div}(X)$}

\index{$D_h$}
\begin{definition}
For any support function $h$ on a fan $\Delta$ we have a corresponding
$T$ invariant divisor given by
\[
D_h=\sum_{\sigma\in\Delta(1)} h(n(\sigma))V(\sigma),
\]
where $n(\sigma)$ is an element of $N$ which generates $\sigma$.
(Note, $n(\sigma)$ is only defined when $\sigma$ is one dimensional.)
\end{definition}

\begin{theorem}
We have
\[
T{\rm Div}(X) = \sum_{\sigma\in\Delta(1)}\Z V(\sigma)\cong SF(\Delta),
\]
where $SF$ is as in Definition \ref{def:SF}.
\end{theorem}

For further details, see \S 4.2 above, \S 3.3 of \cite{F}, or 
\S 2.1 in \cite{O}.

\index{valuation $v_f$}
\begin{definition}
For any  rational function $f$ on a variety $X$ (i.e., defined by
polynomials), there is a corresponding valuations on divisors given
by 
\[
v_f(D) = {\rm order\ of\ vanishing\ of } f\ {\rm along }\ D.
\]
If $f$ has a pole along $D$ then the valuation is negative.
\end{definition}

Rather than give a proper definition of ``order of vanishing'',
we give an example below.

\index{$T$-invariant\\ divisors}
\begin{definition}
If $f$ is a rational function on a variety $X$ (i.e., defined by
polynomials), then there is a corresponding Weil divisor, given by
the free abelian group
\[
(f) = \sum_{D\in {\rm Div}(X)} v_f(D)D.
\]
A divisor of this form is called a {\bf principal divisor}.
\end{definition}
\index{principal divisor}

\begin{example}
For $\P^2$, with projective coordinates $X,Y,Z$, consider the function
\[
f(X,Y,Z) = \frac{X^2(Y-3Z)}{(XY-Z^2)(X-Z)}.
\]
This funtion is zero on the lines $X=0$, $Y=3Z$.  It vanishes to order
$2$ on $X=0$.  It has poles along the curve $XY=Z^2$ and the line $X-Z$,
so the corresponding divisor is:
\[
(f)=2(X=0) + (Y=3Z) - (XY=Z^2) - (X=Z).
\]
\end{example}

\index{linear equivalence}
\begin{definition}
Two Weil divisors $D_1$ and $D_2$
are {\bf linearly equivalent},
written $D_1\sim D_2$, if there is some function $f$ with
\[
(f) = D_1 - D_2.
\]
\end{definition}

\begin{example}
Any two lines on $\P^2$ are linearly equivalent, since if 
$l_1(X,Y,Z)$ and $l_2(X,Y,Z)$ are linear functions defining lines
$L_1$ and $L_2$, then
\[
\left(\frac{l_1}{l_2}\right)=L_1 - L_2.
\]
\end{example}

\index{complete linear system}
\begin{definition}
If $D$ is a divisor on $X$, then the {\bf complete linear system}
defined by $D$ is given by
\[
|D|=\{D'\in {\rm Div}(X) | D \sim D'\}.
\]
$|D|$ has the structure of a projective space, since any $D'\in|D|$
corresponds to some function $f$, with $(f)=D'-D$, 
and we  also have $f\mapsto (f)-D\in|D|$.
Functions can
be added together, and multiplied by elements of $k$, but $(cf)=(f)$
for any non zero constant $c$, so we quotient the vector space of non-zero
functions on $X$ by $\cccF^\times$ to obtain the desired
projective space.  
\end{definition}
\index{$|D|$}

Linear systems are important because they can be used to define
functions from varieties to projective space.  If $D_1\sim D_2\sim D_3$
are divisors, locally defined by the vanishing of some polynomials
$f_1,f_2,f_3$ on some affine piece of a variety, then there is a function
given (locally) by
\[
x\mapsto (f_1(x) : f_2(x) : f_3(x)).
\]
The fact that $D_1\sim D_2$ means that this map can be patched together
globally, with the ratio $f_1(x)/f_2(x)$ giving the value of the function
locally.

\begin{theorem}
Let $X$ be a smooth toric variety.
For a support function $h$ on a fan $\Delta$,
the linear system $|D_h|$ defines a smooth embedding of $X(\Delta)$ in
projective space if and only if $h$ is strictly upper convex.
\end{theorem}

This is equivalent to the lemma on page 69 of Fulton \cite{F}.

\subsection{An explicit basis for $L(D)$}

Let $\Delta$ be a fan in a lattice $L$. Denote by
$\tau_1$, ..., $\tau_n$ the edges or rays of
the fan and let $v_i$ denote the first (smallest)
lattice point along the ray $\tau_i$.
Let $D_i$ denote the Weil divisor 
\[
D_i={\rm Hom}(\tau_i^*\cap L^*,\cccF^\times),
\]
which may be regarded as the closure of the orbit of $T$ acting on 
the edge $\tau_i$. Let 
\[
P_D=\{x=(x_1,...,x_n)\ |\ \langle x,v_i\rangle \geq -d_i, 
\ \forall 1\leq i\leq k\}
\]
denote the polytope associated to the Weil
divisor $D=d_1D_1+...+d_kD_k$,
where $D_i$ is as above. 
\index{$P_D$}
There is a fairly simple condition which determines whether or not
$D$ is Cartier - see the exercise on page 62 of \cite{F}. 
Moreover, there is a fairly simple condition which tells you
whether or not
$D$ is ample - see the proof of the proposition on page
68 of \cite{F}. 

\index{line bundle}
If $F$ is a topological field, a {\bf line bundle} on an 
$F$-variety $X$ is given by a $F$-manifold 
$L$ with a surjective map
\[
\pi: L\ra X,
\]
such that the inverse image $\pi^{-1}(x)$
is a one dimensional vector space over the underlying field $F$.

\index{equivariant line bundle}
For a toric variety $X(\Delta)$, containing an 
open dense torus $(k^\times )^n$, a line bundle 
$\pi: L\ra X(\Delta)$ is called an {\bf equivariant line bundle} if
$(\cccF^\times)^n$ acts on $L$, and for all 
$a\in (\cccF^\times)^n$ and all $v\in L$
we have
\[
\pi(a\cdot v)=a\cdot\pi(v).
\]

We won't go into sheaves in detail, but all invertible sheaves on (a smooth) 
toric variety $X$ are defined by $T$-invariant Weil divisors.
Moreover, certain computations of the cohomology of invertible
sheaves on $X$ boils down to combinatorial computations.
For example, we have results like the following.

\begin{theorem}
For a polyhedron $P$ in a lattice $L^*$, we have a support function
$h=h_P$ on the toric variety $X(\Delta(P))$, and a corresponding
sheaf ${\cal O}(D_h)$.  The space of global sections of the sheaf,
$H^0(X(\Delta(P)),{\cal O}(D_h))$
is a finite dimensional vector space with basis given by the set of
lattice points in $P\cap L^*$.
\end{theorem}

For further details, see \S 3.4 in \cite{F} or
Lemma 2.3 in \cite{O}.

\index{$L(D)$}
We define the Riemann-Roch space $L(D)$ by
\[
L(D)= \Gamma(X,{\cal O}(D))=H^0(X,{\cal O}(D))
\]
(see for example Griffiths and Harris \cite{GH}, page 136,
for a natural isomorphism between this space and the
``usual definition'').
By Fulton \cite{F}, page 66, we have
\[
L(D)=\oplus_{u\in P_D} \cccF\cdot \chi^u,
\]
where $\chi=(x_1,...,x_n)$ and $\chi^u$ is the
associated monomial in multi-index notation.

\begin{example}
We continue example \ref{ex:toric} above.

Let $\Delta$ be the fan generated by 
\[
v_1=2e_1-e_2,\ \ v_2=-e_1+2e_2,\ \ v_3=-e_1-e_2.
\]

In the notation above, the divisor $D=d_1D_1+d_2D_2+d_3D_3$
is a Cartier divisor
\footnote{This is an Exercise on page 65 of \cite{F},
the solution of which is an easy calculation using the
Exercise on page 62, which is in turn, basically
solved in the back of the book.}
if and only if
$d_1\equiv d_2\equiv d_3\ ({\rm mod}\ 3)$.
Let 
\[
\begin{array}{c}
P_D=\{(x,y)\ |\ \langle (x,y),v_i\rangle \geq -d_i, \ \forall i\}\\
=\{(x,y)\ |\ 2x-y \geq -d_1, -x+2y \geq -d_2, -x-y \geq -d_3 
\}
\end{array}
\]
denote the polytope associated to $D$. The 
\verb+divisor_polytope+ command in \verb+toric.mag+
implements an algorithm which determines the inequalities 
describing $P_D$ in general (at the end of the 
file \verb+toric.mag+ \cite{J2} there are examples of
this).

If $d_1=d_2=6$ and $d_3=0$ then $P_D$ is a triangle in the plane 
with vertices at $(-6,-6)$, $(-2,2)$, and $(2,-2)$. Note that it remains
invariant under the action of $G$. Moreover, the $G$ action
on $P_D\cap L^* $ has $7$ singleton orbits
(the lattice points $(-i,-i)$, where $0\leq i\leq 6$)
and $12$ orbits of size $2$.
A basis for the Riemann-Roch space $L(D)$ is returned by
the \verb+toric.mag+ commands:

{\footnotesize{
\begin{verbatim}
> load "/home/wdj/magmafiles/toric.mag"; 
Loading "/home/wdj/magmafiles/toric.mag"

toric.mag for MAGMA 2.8, wdj, version 3-2-2002
available functions: dual_semigp_gens, cart_prod_lists, 
in_dual_cone, max_vectors, create_affine_space, toric_points, 
ideal_affine_toric_variety, affine_toric_variety, 
embedding_affine_toric_variety, desing_affine_toric_variety, 
divisor_polytope, divisor_polytope_lattice_points, riemann_roch,  
toric_code, toric_codewords, ... 

> 
> DB := LatticeDatabase();      
> Lat := Lattice(DB, 2, 16);Lat;
Standard Lattice of rank 2 and degree 2
> Cones:=[[[2,-1],[-1,2]],[[-1,2],[-1,-1]],[[-1,-1],[2,-1]]];
> Div:=[6,6,0];
> RR:=riemann_roch(Div,Cones,Lat);          
> RR;                                       
[
    1/(x1^6*x2^6),
    1/(x1^5*x2^5),
    1/(x1^5*x2^4),
    1/(x1^4*x2^5),
    1/(x1^4*x2^4),
    1/(x1^4*x2^3),
    1/(x1^4*x2^2),
    1/(x1^3*x2^4),
    1/(x1^3*x2^3),
    1/(x1^3*x2^2),
    1/(x1^3*x2),
    1/x1^3,
    1/(x1^2*x2^4),
    1/(x1^2*x2^3),
    1/(x1^2*x2^2),
    1/(x1^2*x2),
    1/x1^2,
    x2/x1^2,
    x2^2/x1^2,
    1/(x1*x2^3),
    1/(x1*x2^2),
    1/(x1*x2),
    1/x1,
    x2/x1,
    1/x2^3,
    1/x2^2,
    1/x2,
    1,
    x1/x2^2,
    x1/x2,
    x1^2/x2^2
]

\end{verbatim}
}}
\end{example}

\section{Application to error-correcting codes}
\label{sec:codes}

Error-correcting codes associated to a toric
variety were introduced by J. Hansen \cite{H}
in the case of surfaces. In many cases, he found
good estimates for the parameters $n$ (the length),
$k$ (the dimension), and $d$ (the minimum distance)
of the code. The estimates on $d$
were based on more general techniques by his student
S. Hansen \cite{Han2} (in fact, this
paper was originally part of his PhD thesis \cite{Han1}). 

The next section briefly recalls his J. Hansen's construction.

The section after that gives a construction which is a little more 
general than that in \cite{H}, though it still
falls in the framework of the general class of
codes constructed in \cite{Han2}.

\subsection{Hansen codes}
\index{Hansen codes}

We recall briefly some codes associated to a toric surface,
constructed by J. Hansen \cite{H}.

Let $F=\fff_q$ be a finite field with $q$ elements
and let $\overline{F}$ denote a separable algebraic closure. 
Let $L$ be a lattice in $\qqq^2$ generated by $v_1, v_2\in \zzz^2$,
$P$ a polytope in $\qqq^2$,
and $X(P)$ the associated toric surface. Let $P_L=P\cap \zzz^2$.

There is a dense embedding of $GL(1)\times GL(1)$ into $X(P)$
given as follows.
Let $T_L=Hom_{\zzz}(L,GL(1))$ (which is $\cong GL(1)\times GL(1)$
by sending $t=(t_1,t_2)$ to $m_1v_1+m_2v_2\longmapsto
e(\ell)(t)=t_1^{m_1}t_2^{m_2}$)
and let $e(\ell):T_L\rightarrow \overline{F}$
be defined by $e(\ell)(t)=t(\ell)$, $t\in T_L$.

Impose an ordering on the set $T_L(F)$ (changing the ordering
leads to an equivalent code).
Define the code $C=C_P\subset F^n$ to be the linear code generated
by the vectors
\begin{equation}
\label{eqn:code}
B=\{(e(\ell)(t))_{t\in T_L(F)}\ |\ \ell \in L\cap P_L\},
\end{equation}
where $n=(q-1)^2$. In some special cases, the dimension if $C$ is known
and an estimate of its minimum distance can be given
(see Theorem \ref{thrm:hansen} below).

\begin{example}
The case $q=2$ is trivial.

The first non-trivial example is for $q=3$. 
Let $L=\zzz^2$, so 
\[
T_L(F)=\{(1,1),(1,2),(2,1),(2,2)\},
\]
which we use as our ordering.
Let $P_L$ be the polytope
with vertices $(0,0)$, $(1,0)$, and $(0,1)$. 

In this case,

\begin{center}
\begin{tabular}{|c|c|c|c|c|} \hline
$t$ & (1,1)&(1,2)&(2,1)&(2,2) \\ \hline
$e(0,0)(t)$ & 1 & 1 & 1 & 1 \\ \hline
$e(1,0)(t)$ & 1 & 1 & 2 & 2 \\ \hline
$e(0,1)(t)$ & 1 & 2 & 1 & 2 \\ \hline
\end{tabular}
\end{center}
Thus, 
\[
\begin{array}{c}
C={\rm span}_F\{(1,1,1,1),(1,1,2,2),(1,2,1,2)\}\\
=\{(0,0,0,0), (1,1,1,1),(2,2,2,2),(0,1,0,1), \\
(1,2,1,2), (0,2,0,2),(2,1,2,1),(1,0,1,0), \\
(2,1,2,1), (2,0,2,0),(1,2,1,2), (0,0,1,1),\\
 (0,0,2,2), (1,1,2,2), (1,1,0,0), (2,2,0,0),\\
 (2,2,1,1), (2,0,0,1),(1,0,0,2), (0,2,1,0),\\
(0,1,2,0), (0,2,2,1), (0,1,1,2),  (1,2,2,0),\\ 
(2,1,1,0), (2,1,0,2), (1,2,0,1)
\}.
\end{array}
\]
It's minimum distance is $2$.

\end{example}

Hansen gives lower bounds on the minimum distance $d$ of 
such codes in the cases:

(a) $P$ is an isoceles triangle with vertices $(0,0)$,$(a,a)$,$(0,2a)$,

(b) $P$ is an isoceles triangle with vertices $(0,0)$,$(a,0)$,$(0,a)$, or

(c) $P$ is a rectangle with vertices $(0,0)$,$(a,0)$,$(0,b)$,$(a,b)$,

\noindent
provided $q$ is ``sufficient large''. His precise result is 
recalled below.

\begin{theorem}
(Hansen \cite{Han})
\label{thrm:hansen}
Let $a$, $b$ be positive integers. Let $P$ be the polytope 
defined in (a)-(c) above. 

\begin{itemize}
\item[(a)] Assume $q>2a+1$. The code $C=C_P$ has
\[
n=(q-1)^2,\ \ k=(a+1)^2,\ \ d\geq n-2a(q-1).
\]

\item[(b)] Assume $q>a+1$. The code $C=C_P$ has
\[
n=(q-1)^2,\ \ k=(a+1)(a+2)/2,\ \ d\geq n-a(q-1).
\]

\item[(c)] Assume $q>{\rm max}(a,b)+1$. The code $C=C_P$ has
\[
n=(q-1)^2,\ \ k=(a+1)(b+1),\ \ d\geq n-a(q-1)-b(q-1)+ab.
\]

\end{itemize}
\end{theorem}

\begin{magma}
We examine an example of part (b) of Hansen's theorem stated above.
First load \verb+toric.mag+ as follows (replace 
\verb+/home/wdj/magmafiles+ by your path to 
\verb+toric.mag+). 

{\footnotesize{
\begin{verbatim}

> load "/home/wdj/magmafiles/toric.mag"; 
Loading "/home/wdj/magmafiles/toric.mag"

toric.mag for MAGMA 2.8, wdj, version 8-28-2001
available functions: dual_semigp_gens, cart_prod_lists, 
in_dual_cone, max_vectors, create_affine_space, toric_points, 
ideal_affine_toric_variety, affine_toric_variety, 
embedding_affine_toric_variety, desing_affine_toric_variety, 
toric_code, toric_codewords, ... 

> Polyb:=[[0,0],[0,1],[1,0],[1,1],[0,2],[1,2],[2,0],[2,1],[0,3],[3,0]];
> C:=toric_code(Polyb,GF(3),Lat);C;
[4, 4, 1] Linear Code over GF(3)
\end{verbatim}
}}

The command \verb+toric_codewords(Polyb,GF(3),Lat);+ returns the list
of all codewords:
\[
\begin{array}{c}
\{
  ( 0, 0, 0, 0 ),
    ( 1, 0, 0, 0 ),
    ( 2, 0, 0, 0 ),
    ( 2, 1, 0, 0 ),
    ( 0, 1, 0, 0 ),\\
\    ( 1, 1, 0, 0 ),
    ( 1, 2, 0, 0 ),
    ( 2, 2, 0, 0 ),
    ( 0, 2, 0, 0 ),
    ( 0, 2, 1, 0 ),
    ( 1, 2, 1, 0 ),\\
\    ( 2, 2, 1, 0 ),
    ( 2, 0, 1, 0 ),
    ( 0, 2, 0, 0 ),
    ( 0, 2, 1, 0 ),
    ( 1, 2, 1, 0 ),
    ( 2, 2, 1, 0 ),\\
\    ( 2, 0, 1, 0 ),
    ( 0, 0, 1, 0 ),
    ( 1, 0, 1, 0 ),
    ( 1, 1, 1, 0 ),
    ( 2, 1, 1, 0 ),
    ( 0, 1, 1, 0 ),\\
\    ( 0, 1, 2, 0 ),
    ( 1, 1, 2, 0 ),
    ( 2, 1, 2, 0 ),
    ( 2, 2, 2, 0 ),
    ( 0, 2, 2, 0 ),
    ( 1, 2, 2, 0 ),
    ( 1, 0, 2, 0 ),\\
\    ( 2, 0, 2, 0 ),
    ( 0, 0, 2, 0 ),
    ( 0, 0, 2, 1 ),
    ( 1, 0, 2, 1 ),
    ( 2, 0, 2, 1 ),
    ( 2, 1, 2, 1 ),
    ( 0, 1, 2, 1 ),\\
\    ( 1, 1, 2, 1 ),
    ( 1, 2, 2, 1 ),
    ( 2, 2, 2, 1 ),
    ( 0, 2, 2, 1 ),
    ( 0, 2, 0, 1 ),
    ( 1, 2, 0, 1 ),
    ( 2, 2, 0, 1 ),\\
\    ( 2, 0, 0, 1 ),
    ( 0, 0, 0, 1 ),
    ( 1, 0, 0, 1 ),
    ( 1, 1, 0, 1 ),
    ( 2, 1, 0, 1 ),
    ( 0, 1, 0, 1 ),
    ( 0, 1, 1, 1 ),\\
\    ( 1, 1, 1, 1 ),
    ( 2, 1, 1, 1 ),
    ( 2, 2, 1, 1 ),
    ( 0, 2, 1, 1 ),
    ( 1, 2, 1, 1 ),
    ( 1, 0, 1, 1 ),
    ( 2, 0, 1, 1 ),\\
\    ( 0, 0, 1, 1 ),
    ( 0, 0, 1, 2 ),
    ( 1, 0, 1, 2 ),
    ( 2, 0, 1, 2 ),
    ( 2, 1, 1, 2 ),
    ( 0, 1, 1, 2 ),\\
\    ( 1, 1, 1, 2 ),
    ( 1, 2, 1, 2 ),
    ( 2, 2, 1, 2 ),
    ( 0, 2, 1, 2 ),
    ( 0, 2, 2, 2 ),
    ( 1, 2, 2, 2 ),
    ( 2, 2, 2, 2 ),\\
\    ( 2, 0, 2, 2 ),
    ( 0, 0, 2, 2 ),
    ( 1, 0, 2, 2 ),
    ( 1, 1, 2, 2 ),
    ( 2, 1, 2, 2 ),
    ( 0, 1, 2, 2 ),\\
\    ( 0, 1, 0, 2 ),
    ( 1, 1, 0, 2 ),
    ( 2, 1, 0, 2 ),
    ( 2, 2, 0, 2 ),
    ( 0, 2, 0, 2 ),\\
\    ( 1, 2, 0, 2 ),
    ( 1, 0, 0, 2 ),
    ( 2, 0, 0, 2 ),
    ( 0, 0, 0, 2 )\}.
\end{array}{c}
\]

\end{magma}

Other examples (we invite the reader to experiment her/himself)
support the following conjecture.

\begin{conjecture}
Hansen's theorem holds if one replaces the symbol
$\geq$ in the lower estimate for $d$ 
by $=$.
\end{conjecture}

\subsection{Other toric codes}

Though the results in \cite{F} apply to
toric varieties over $\ccc$, we shall work
over a finite field $\fff_q$ having $q=p^k$ elements,
where $p$ is a prime and $k\geq 1$. 
We assume that the results of \cite{F} have 
analogs over $\fff_q$.

Let $M\cong \zzz^n$ be a lattice in $V=\rrr^n$
and let $N\cong \zzz^n$ denote its dual.
Let $\Delta$ be a fan (of rational cones, with
respect to $M$) in $V$ and let $X=X(\Delta)$ denote the
toric variety associated to $\Delta$. 
Let $T$ denote a dense torus in $X$.

Let $D=P_1+...+P_n$ be a positive $1$-cycle on $X$,
where the points $P_i\in X(\fff_q)$ are distinct.
Let $G$ be a $T$-invariant divisor on $X$ which does not ``meet'' $D$,
in the sense that no element of the support of $D$ 
intersects any element in the support of $G$. We
write this as
\[
{\rm supp}(G)\cap {\rm supp}(D)=\emptyset.
\]
Some additional assumptions on $G$ and $D$ shall be made later.
Let 
\[
L(G)=\{0\}\cup \{f\in \fff_q(X)^\times \ |\ div(f)+G\geq 0\}
\]
denote the Riemann-Roch space associated to $G$.
According to \cite{F}, \S 3.4, there is a polytope 
$P_G$ in $V$ such that $L(G)$ is spanned by the
monomials $x^a$ (in multi-index notation), for $a\in P_G\cap N$.
Let $C_L=C_L(D,G)$ denote the code define by
\[
C_L=\{(f(P_1),...,f(P_n))\ |\ f\in L(G)\}.
\]
This is the {\bf Goppa code associated to $X$, $D$, and $G$}. 
The dual code is denoted
\[
C=\{(c_1,...,c_n)\in \fff_q^n\ |\ \sum_{i=1}^n c_if(P_i)=0,
\forall f\in L(G)\}.
\]
\index{Goppa code associated to $X$, $D$, and $G$}

\begin{example}
\label{ex:7.3}

Let $\Delta$ be the fan generated by 
\[
v_1=2e_1-e_2,\ \ v_2=-e_1+2e_2,\ \ v_3=-e_1-e_2.
\]
Let $X$ be the toric variety associated to $\Delta$.

In the notation above, the divisor $D=d_1D_1+d_2D_2+d_3D_3$
is a Cartier divisor if and only if
$d_1\equiv d_2\equiv d_3\ ({\rm mod}\ 3)$.

\begin{center}
\setlength{\unitlength}{.01cm}

\begin{picture}(256.00,150.00)(-120.00,0.00)
\thicklines
\put(0.00,64.00){\circle*{4}} 
\put(0.00,128.00){\circle*{4}} 

\put(-64.00,0.00){\circle*{4}} 
\put(0.00,-64.00){\circle*{4}} 

\put(-64.00,-64.00){\circle*{4}} 
\put(-64.00,64.00){\circle*{4}} 
\put(-64.00,128.0){\circle*{4}} 

\put(64.00,0.00){\circle*{4}} 
\put(128.00,0.00){\circle*{4}} 
\put(192.00,0.00){\circle*{4}} 

\put(0.00,0.00){\vector(-1,2){64.00}} 
\put(0.00,0.00){\vector(-1,-1){64.00}} 

\put(64.00,64.00){\circle*{4}} 
\put(128.00,0.00){\circle*{4}} 
\put(0.00,128.00){\circle*{4}} 
\put(-128.00,0.00){\circle*{4}} 

\put(128.00,80.00){$\sigma_1$} 
\put(-128.00,64.00){$\sigma_2$} 
\put(50.00,-100.00){$\sigma_3$} 

\put(135.00,-75.00){$v_1$} 
\put(-75.00,135.00){$v_2$} 
\put(-85.00,-85.00){$v_3$} 

\put(0.00,-64.00){\circle*{4}} 
\put(64.00,-64.00){\circle*{4}} 
\put(128.00,-64.00){\circle*{4}} 

\put(0.00,0.00){\vector(2,-1){129.00}} 
\end{picture}

\end{center}
\vskip .5in

Let 
\[
\begin{array}{c}
P_D=\{(x,y)\ |\ \langle (x,y),v_i\rangle \geq -d_i, \ \forall i\}\\
=\{(x,y)\ |\ 2x-y \geq -d_1, -x+2y \geq -d_2, -x-y \geq -d_3 
\}
\end{array}
\]
denote the polytope associated to the Weil
divisor $D=d_1D_1+d_2D_2+d_3D_3$,
where $D_i$ is as above. 

Let 
\[
G=10D_3,\ \ \ 
D=D_1+D_2+D_3.
\]
Then $P_G$ is a triangle in the plane 
with vertices at $(0,0)$, $(-10/3,20/3)$, and $(20/3,10/3)$,
see Figure \ref{fig:polytopeG}. It's area is $50/3$.

\begin{figure}[h]
\begin{minipage}{\textwidth}
\begin{center}
\vspace{1.0 cm}
\includegraphics[height=5cm,width=5cm]{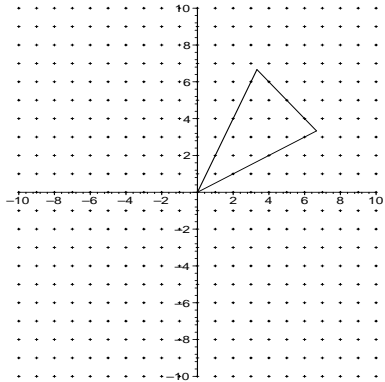}
\end{center}
\end{minipage}
\caption{The polytope associated to the divisor $G$.}
\label{fig:polytopeG}
\end{figure}

Let $T\subset X$ denote the dense torus of $X$.
In this example, the patch $U_{\sigma_1}$ is an affine variety with 
coordinates $x_1,x_2,x_3$ given by $x_1^3-x_2x_3=0$. 
The torus embedding
$T\hookrightarrow U_{\sigma_1}$ is given by
sending $(t_1,t_2)$ to 
$(x_1,x_2,x_3)=(t_1t_2,t_1t_2^2,t_1^2t_2)$.
The patch $U_{\sigma_2}$ is an affine variety with 
coordinates $y_1,y_2,y_3$ given by $y_2^2-y_1y_3=0$. 
The torus embedding
$T\hookrightarrow U_{\sigma_2}$ is given by
sending $(t_1,t_2)$ to 
$(y_1,y_2,y_3)=(t_1^{-2}t_2^{-1},t_1^{-1},t_1^{-1}t_2)$.
The patch $U_{\sigma_3}$ is an affine variety with 
coordinates $z_1,z_2,z_3$ given by $z_2^2-z_1z_3=0$. 
The torus embedding
$T\hookrightarrow U_{\sigma_3}$ is given by
sending $(t_1,t_2)$ to 
$(x_1,x_2,x_3)=(t_1^{-1}t_2^{-2},t_2^{-1},t_1t_2^{-1})$.

In the local coordinates of $U_{\sigma_1}$, the
space $L(G)$ has as a basis,
\[
\begin{array}{c}
\{h_i\}=
\{ 1, x_1 x_2, x_1^2 x_2, x_1 x_2^2, x_1^2 x_2^2, 
x_1^3 x_2^2, x_1^4 x_2^2, x_1^2 x_2^3, x_1^3 x_2^3,\\ 
x_1^4 x_2^3, x_1^5 x_2^3, x_1^6 x_2^3, x_1^2 x_2^4, 
x_1^3 x_2^4, x_1^4 x_2^4, x_1^5 x_2^4, x_1^6 x_2^4, \\
x_1^3 x_2^5, x_1^4 x_2^5, x_1^5 x_2^5, x_1^3 x_2^6, 
x_1^4 x_2^6\}\, .
\end{array}
\]
In particular, it is 22-dimensional.

Choose points $P_1,...,P_n\in X(\fff_q)$ and let
\[
C=\{(c_1,...,c_n)\in \fff_q^n\ |\ \sum_{i=1}^n c_i f(P_i)=0,
\forall f\in L(G)\}.
\]
For example, we know
\[
\begin{array}{c}
|X(\fff_2)|=7,\ \ |X(\fff_3)|=13,\ \ |X(\fff_4)|=21,\ \\  
|X(\fff_5)|=31,\ \ |X(\fff_7)|=57,\ \  |X(\fff_8)|=73,...
\end{array}
\]

\end{example}

One can (though we do not do so here)
use \cite{Han2} to
estimate the parameters $n,k,d$ for 
$C_L=C_L(G,D,X)$.
The result, roughly speaking, says that
if $X$ is a ``nice'' toric variety of dimension $N$,
if $G$ is a ``nice'' divisor on $X$, and if
$D$ is a 1-cycle of ``sufficiently large degree'',
then
\[
n={\rm deg}(D),
\]
\[
k={\rm dim}\, H^0(X,{\cal O}(G))=|P_G\cap M|,
\]
and $d$ is bounded from below (roughly) by $n-N!\cdot k$.

\section{The toric automorphism group}
\label{sec:autoric}

Since a (dense) torus $T$ acts on $X=X(\Delta)$, 
$T$ forms a subgroup of the automorphism group of $X$.
From Lemma \ref{lemma:functor}, we see that another source of 
automorphisms comes from the automorphisms of the fan $\Delta$.

Let $\Delta$ be a fan in $V=L\otimes \qqq$,
$T_L=$Hom$_{ab.gp.}(L^*,\ggg_m)$, and let
$\Delta'$ be a fan in $V'=L'\otimes \qqq$,
$T_{L'}=$Hom$_{ab.gp.}((L')^*,\ggg_m)$.

\begin{lemma} (Theorem 1.13, Oda \cite{O})
Let $f:T_{L'}\rightarrow T_{L'}$
be a homomorphism. If $\psi : X(\Delta)\rightarrow  X(\Delta')$
is a morphism which is equivariant with respect to $f$
then there is a mapping of fans $\phi_\qqq:\Delta'\rightarrow 
\Delta$ such that $\psi=\phi^*$
(in the sense of Lemma \ref{lemma:functor}).
\end{lemma}

In particular, an equivariant automorphism of a toric variety 
is a toric automorphism.

\begin{example}
Consider the fan
\[
\Delta =\{\{0\},\ \qqq_{\geq 0},\ \qqq_{\leq 0}
 \}
\]
in Example \ref{example:P1} above. This corresponds
to the toric variety $X(\Delta)=\ppp^1$.
The auomorphism $\phi$ of $\Delta$ which 
swaps $\qqq_{\geq 0}$ with $\qqq_{\leq 0}$ and leaves 
$\{0\}$ alone induces a toric morphism $\phi^*$
on $X(\Delta)$. Locally, this morphism $\phi^*$ is simply
the ``gluing'' map discussed in Example \ref{example:P1}.

\end{example}

\begin{remark}
We observe a connection with the previous section 
\ref{sec:codes} on codes.

Note that if $G$ is a subgroup of the toric automorphism group of
$X(\Delta)$ and $G$ induces an automorphism of the lattice
$L$ which leaves the polytope $P_L$ invariant then 
$G$ induces an automorphism on the code $C$ defined
in (\ref{eqn:code}).
\end{remark}

\begin{remark}
Here's an algorithm to find the automorphism group of a fan
(hence the toric automorphism group of a toric variety):

Recall all the cones in a fan are polyhedral.
It follows from this that the 
automorphism group of a fan
\footnote{Recall Definition \ref{defn:autofan}.}, 
is contained in the
automorphism group of its set of rays
(all the one-dimensional cones in the fan)
\footnote{The rays do not uniquely determine the fan, 
nor even the cones themselves.}.
To find the automorphism group of the rays,
we 
\begin{enumerate}
\item
 write each ray as $\tau_i=\qqq_{\geq 0}\cdot v_i$,
for some $v_i\in L$,

\item
 for each $g\in {\rm Aut}(L)\subset {\rm Aut}(\qqq^n)$ do:

\begin{enumerate}
\item
for each $i$, compute $g(v_i)$, and ask: is $g(v_i)\in \tau_j$,
for some $j$? 

\item
if no for some $i$, go to the next $g$;
if yes for all $i$, add $g$ to ${\rm Aut}(Rays)$.

\item
Go to the next $g$, if one exists. If no more $g$
exist, then return ${\rm Aut}(Rays)$.
\end{enumerate}
\item
Now that ${\rm Aut}(Rays)$ has been determined,
perform the following double for loop.
For each $g$ in ${\rm Aut}(Rays)$ do:
for each cone $\sigma\in \Delta$ do:
does $g$ send $\sigma$ to some cone
$\sigma'\in \Delta$. If the answer is ``yes'',
for all $\sigma$ then add $g$ to ${\rm Aut}(\Delta)$. 
Go to the next $g$, if one exists. If no more $g$
exist, then return ${\rm Aut}(\Delta)$.
\end{enumerate}
\end{remark}

\begin{example}
\label{ex:toric}
Let $\Delta$ be the fan generated by 
\[
v_1=2e_1-e_2,\ \ v_2=-e_1+2e_2,\ \ v_3=-e_1-e_2.
\]
There is a picture of this in Example \ref{ex:7.3} above.

Let $G$ denote the automorphism group of $\Delta$ 
generated by the map $g$ sending $(x,y)$ to $(y,x)$, swapping 
$v_1$ and $v_2$ and leaving $v_3$ fixed. By Theorem
1.13 in Oda \cite{O}, 
this corresponds to a $T$-equivariant automorphism of
$X(\Delta)$.

In this example, the patch $U_{\sigma_1}$ is an affine variety with 
coordinates $x_1,x_2,x_3$ given by $x^3-x_2x_3=0$. 
The automorphism $g$ acts on $U_{\sigma_1}$ sending
$(x_1,x_2,x_3)$ to $(x_1,x_3,x_2)$. The torus embedding
$T\hookrightarrow U_{\sigma_1}$ is given by
sending $(t_1,t_2)$ to 
$(x_1,x_2,x_3)=(t_1t_2,t_1t_2^2,t_1^2t_2)$.

The patch $U_{\sigma_2}$ is an affine variety with 
coordinates $y_1,y_2,y_3$ given by $y^2-y_1y_3=0$. 
The automorphism $g$ does not act on $U_{\sigma_2}$.
The torus embedding
$T\hookrightarrow U_{\sigma_2}$ is given by
sending $(t_1,t_2)$ to 
$(y_1,y_2,y_3)=(t_1^{-2}t_2^{-1},t_1^{-1},t_1^{-1}t_2)$.

The patch $U_{\sigma_3}$ is an affine variety with 
coordinates $z_1,z_2,z_3$ given by $z^2-z_1z_3=0$. 
The automorphism $g$ does not act on $U_{\sigma_3}$.
The torus embedding
$T\hookrightarrow U_{\sigma_3}$ is given by
sending $(t_1,t_2)$ to 
$(x_1,x_2,x_3)=(t_1^{-1}t_2^{-2},t_2^{-1},t_1t_2^{-1})$.
The automorphism $g$ sends $U_{\sigma_2}$ to $U_{\sigma_3}$
by sending $(y_1,y_2,y_3)$ to $(z_1,z_2,z_3)$.

\end{example}

\section{Betti numbers}

\index{$b_k$}
\index{$d_k$}
\index{Betti number}
Let $X=X(\Delta)$ be an $n$-dimensional 
smooth toric variety associated to a fan $\Delta$.
Let $d_k$ denote the number of distinct $k$-dimensional
cones in $\Delta$. 
The {\bf $k^{th}$ Betti number} of $X$ is the rank of $H^k(X,\zzz)$. 
Let $b_k$ denote the the $k^{th}$ Betti number of $X$.
It is known (see Fulton, \S 4.5) that 
\[
b_{2k}=\sum_{j=k}^n (-1)^{j-k} 
\left(
\begin{array}{c}
j \\ 
k
\end{array}
\right)
d_{n-j}
.
\]
The {\tt toric} package can compute this.

\begin{example}
Let $\Delta$ be the fan whose cones of maximal dimension
are defined by 
\[
\sigma_1=\rrr_{\geq 0}(1,0)+\rrr_{\geq 0}(0,1),
\ \ \ \ 
\sigma_2=\rrr_{\geq 0}(0,1)+\rrr_{\geq 0}(-1,0),
\]
\[ 
\sigma_3=\rrr_{\geq 0}(-1,0)+\rrr_{\geq 0}(0,-1),
\ \ \ \ 
\sigma_4=\rrr_{\geq 0}(0,-1)+\rrr_{\geq 0}(1,0).
\]
The toric variety associated to this fan is $X=\ppp^1\times \ppp^1$.
The {\tt toric} package commands to compute the
Betti numbers $b_{2k}$ are as follows
\footnote{We will use the GAP version to compute this example.}.

{\footnotesize{
\begin{verbatim}
gap> RequirePackage("guava");
gap> Read("c:/gap/gapfiles/toric.g"); 
gap> Cones:=[[[1,0],[0,1]],[[0,1],[-1,0]],[[-1,0],[0,-1]],[[0,-1],[1,0]]];
gap> betti_number(Cones,1);
gap> betti_number(Cones,2);
gap> euler_characteristic(Cones);
4
\end{verbatim}
}}

\noindent
GAP returns $b_1=0$ to the first command, 
\verb+betti_number(Cones,1);+, and $b_2=2$ 
to the second command. The last command tells us that the Euler
characteristic of $X$ is $4$.
\end{example}

\section{Counting points over a finite field}

Let $X=X(\Delta)$ be a smooth toric variety associated to a fan 
$\Delta$. Let $q$ be a prime power and $\fff=GF(q)$ denote a
field with $q$ elements. It is known (see Fulton,
\S 4.5) that 
\[
|X(\fff)|=\sum_{k=0}^n (q-1)^k d_{n-k}.
\]
The {\tt toric} package can compute this.

\begin{example}
Let $\Delta$ be the fan whose cones of maximal dimension
are defined by 
\[
\sigma_1=\rrr_{\geq 0}(1,0)+\rrr_{\geq 0}(0,1),
\ \ \ \ 
\sigma_2=\rrr_{\geq 0}(0,1)+\rrr_{\geq 0}(-1,0),
\]
\[
\sigma_3=\rrr_{\geq 0}(-1,0)+\rrr_{\geq 0}(0,-1),
\ \ \ \ 
\sigma_4=\rrr_{\geq 0}(0,-1)+\rrr_{\geq 0}(1,0),
\]
so $X(\Delta)=\ppp^1\times \ppp^1$.
The {\tt toric} package commands to compute the
number of points mod $p$, $|X(\fff_p)|$ are as follows
\footnote{We will use the GAP version to compute this example.}.

{\footnotesize{
\begin{verbatim}
gap> RequirePackage("guava");
gap> Read("c:/gap/gapfiles/toric.g"); 
gap> Cones:=[[[1,0],[0,1]],[[0,1],[-1,0]],[[-1,0],[0,-1]],[[0,-1],[1,0]]];
gap> cardinality_of_X(Cones,2);
gap> cardinality_of_X(Cones,3);
gap> cardinality_of_X(Cones,5);
\end{verbatim}
}}

\noindent
GAP returns $|X(\fff_2)|=9$ to the first command, 
\verb+cardinality_of_X(Cones,2);+, 
$|X(\fff_3)|=16$ to the second command, 
$|X(\fff_5)|=36$ to the last command. 
\end{example}

\section{The compactification of a universal elliptic
curve}

In this section, we discuss the compactification of the 
universal family of elliptic
curves with level $N$ structure.

An elliptic curve is given by $\ccc/L$ for some lattice $L$ in $\ccc$.
Two elliptic curves given by lattices $L_1$ and $L_2$ are 
{\bf homothetic} 
\index{homothetic}
if there is some $\alpha\in\ccc$ with
$L_1=\alpha L_2.$
We have
\[
L_1,L_2\ {\rm are\ homothetic\ }\iff\ E_1,E_2\ {\rm are\ isomorphic}.
\]
By suitable scaling, an elliptic curve is isomorphic to 
\[
E_\tau = \ccc/\langle 1,\tau\rangle.
\]
The subgroup of points of order $N$ on $E_\tau$ is given by
\[
E[N]=\{\overline{(n + m\tau)/N} | n,m\in\Z\}\cong(\Z/N\Z)^2
\]
An {\bf elliptic curve with full level $N$ structure} is an elliptic curve
\index{elliptic curve with full level $N$ structure}
together with two points spanning $E[N]$.  E.g., $(E_\tau,1/N,\tau/N).$

Let 
\[
\Gamma(N) = \left\{
\mtwo{a}{b}{c}{d}\in SL_2(\Z)\Big\vert \mtwo{a}{b}{c}{d} \equiv
\mtwo{1}{0}{0}{1} \mod N \right\}.
\]

Any elliptic curve over $\ccc$
can be described as the quotient of the complex plane by
a lattice generated by $1$ and $\tau\in\hhh = \{x\in\ccc |{\rm Im}(x)>0\}$.
We have
\[
E_\tau = \ccc/\langle 1,\tau\rangle .
\]
The Weierstrass $\wp$ function gives the algebraic equations for
$E_\tau$ \cite{Har}.

We can show that 
\[
E_\tau\cong E_{\tau'}\ \iff\ \tau'=g\tau \ {\rm for\ some\ }g\in SL_2(\Z).
\]
So, the quotient $ SL_2(\Z)\setminus\hhh$ parameterises elliptic curves
over $\ccc$.  For a congruence subgroup $\Gamma_0(N)\subset  SL_2(\Z)$,
the space 
\[
Y_0(N) = \Gamma_0(N)\setminus\hhh
\]
parameterizes elliptic curves with
level $N$ structure, i.e., with a subgroup of order $N$.
The problem with $Y_0(N)$ is that it is not compact.  The compactification
is given by 
\[
X_0(N) = \Gamma(N)\setminus\hhh^*
\]
where  $\hhh^* = \hhh\cup\Q\cup\{\infty\}$.
\index{$X_0(N)$}
It is a bit harder to find the compactification of $\eee^\circ(N)$; we
need to know what is the appropriate preimage of the cusps,
(which are the finite number of points 
$\Gamma(N)\setminus(\Q\cup\infty)$), and how these ``glue'' into the rest of
the space.

First of all, we can describe the whole of $\eee ^\circ$ as a quotient:
Let 
\[
\begin{array}{c}
H(N)=\left\{
\mthree{1}{Nr}{Ns}{0}{Na+1}{Nb}{0}{Nc}{Nd+1}
\in SL_3(\Z) 
\Bigg\vert a,b,c,d,r,s\in\Z\right\}\\
\cong (N\Z\oplus N\Z)\rtimes\Gamma(N).
\end{array}
\]
$H(N)$ acts on $\ccc \times\hhh$ by
\[
\mthree{1}{Nr}{Ns}{0}{a}{b}{0}{c}{d}(z,\tau)
=\left(
\frac{z + Nr\tau + Ns}{c\tau+d},
\frac{a\tau+b}{c\tau+d}
\right).
\]
Then we can take $\eee ^\circ = H(N)\setminus \ccc \times\hhh$, 
and there is a map
\[
\xymatrix{
\eee ^\circ\ar@{=}[r]\ar[d]&
H(N)\setminus \ccc \times\hhh \ar[d]& 
\supset(N\Z\tau \oplus N\Z)\setminus\ccc \ar[d]
\\
Y(N)\ar@{=}[r]&
\Gamma(N)\setminus \hhh &\ni \tau{\>\>\>}
}
\]

To construct the compactification of $\eee ^\circ$, we define affine pieces
which patch together, and which ``fill in'' something over each cusp.
We will just look at what happens for the cusp at infinity.

We take the quotient in two steps.

%
The stabilizer in $H(N)$ of the fiber over $\infty$ is given by 

\[
H_\infty=\left\{
\mthree{1}{Nr}{Ns}{0}{1}{Nb}{0}{0}{1}
\Bigg\vert b,r,s\in\Z\right\} .
\]
This is an extension:

\[
\xymatrix{
1\ar[r]&
A\ar[r]&
H_\infty\ar[r]&
B\ar[r]&
1},
\]
where
\[
A=\left\{
\mthree{1}{0}{Ns}{0}{1}{Nb}{0}{0}{1}
\Bigg\vert s,b\in\Z\right\},
\]
and
\[
B=\left\{
\mthree{1}{0}{Nr}{0}{1}{0}{0}{0}{1}
\Bigg\vert r\in\Z\right\}.
\]
We take the quotient in two steps.
$A$ acts on a neigborhood $W$ of $\ccc\times\infty$ by
\[
\mthree{1}{0}{Ns}{0}{1}{Nb}{0}{0}{1}:(z,\tau)\mapsto (z+Ns,\tau+Nb),
\]
and the quotient is a subset of $(\ccc^\times )^2$---a torus.
This quotient map is given by
\[
(z,\tau)\mapsto(w=e^{2\pi iz/N},t=e^{2\pi i\tau/N}).
\]
Then $B$ acts on this quotient by
\[
\mthree{1}{0}{Nr}{0}{1}{0}{0}{0}{1}\ :\ (w,t)\mapsto (t^{Nr}w,t).
\]
We define an infinite fan $\Delta$
in $\Z^2$, with cones spanned by $(n,1),(n+1,1)$.
Then $X(\Delta)\setminus (\ccc^\times )^2$ is
given by an infinite number of copies of
$\P^1$, and we still have an action of $B$ on $X(\Delta)$, which
extends the action of $B$ on $\ccc^2\subset X(\Delta)$.  
This action on 
$X(\Delta)\setminus (\ccc^\times )^2$ is given by sending a $\P^1$ to a $\P^1$
$N$ steps further along. 

The variety $X(\Delta)$ is smooth, and the action of $B$
is fixed point free, and so the quotient $X(\Delta)/B$ is smooth,
and gives a compactification of $\eee ^\circ$ in the neighborbood of
$\infty$.  The fiber over $\infty$ is given by $N$ copies of $\P_1$.

\vskip .2in

{\it Acknowledgements}: We'd like to thank
Jim Fennell, Caroline Melles, Will Traves, and David Zelinski for 
useful discussions.

\printindex

\end{document}